\documentclass[10pt]{articleHJ}
\usepackage[english]{babel}
\usepackage{amsmath, amsthm, amsfonts, mathrsfs, amsfonts,bm,array}
\usepackage{mathtools}
\usepackage{bbm}
\usepackage{amssymb}
\usepackage{graphicx}
\DeclareMathAlphabet\bmcal{OMS}{cmsy}{b}{n}

\newcolumntype{L}{>{$}l<{$} }

\usepackage[text={6.0in,8.6in},centering,letterpaper]{geometry}
\usepackage[numbers,comma,square,sort&compress]{natbib}

\usepackage{tabularx}
\usepackage{enumitem}
\usepackage{subcaption}
\newlength{\defbaselineskip}
\newcommand{\setlinespacing}[1]%
           {\setlength{\baselineskip}{#1 \defbaselineskip}}

\newcommand{\singlespacing}{\setlength{\baselineskip}{\defbaselineskip}}

\newcommand{\D}{\ensuremath{\mathcal{D}}}

\newcommand{\N}{\ensuremath{\mathbb{N}}}
\newcommand{\Z}{\ensuremath{\mathbb{Z}}}
\newcommand{\Q}{\ensuremath{\mathbb{Q}}}
\newcommand{\R}{\ensuremath{\mathbb{R}}}
\newcommand{\C}{\ensuremath{\mathbb{C}}}

\newcommand{\J}{\ensuremath{\mathcal{J}}}
\newcommand{\M}{\ensuremath{\mathcal{M}}}
\newcommand{\I}{\ensuremath{\mathcal{I}}}
\newcommand{\V}{\ensuremath{\mathcal{V}}}

\newcommand{\Y}{\ensuremath{\mathcal{Y}}}
\newcommand{\ttheta}{\ensuremath{{\tilde{\theta}}}}

\renewcommand{\epsilon}{\ensuremath{\varepsilon}}
\renewcommand{\Re}{\ensuremath{\mathrm{Re}}}

\renewcommand{\ker}{\ensuremath{\mathrm{Ker}}}
\newcommand{\Range}{\ensuremath{\mathrm{Range}}}

\renewcommand{\span}{\ensuremath{\mathrm{span}}}

\renewcommand{\min}{\ensuremath{\mathrm{min}}}
\renewcommand{\H}{\ensuremath{\mathcal{H}}}
\newcommand{\K}{\ensuremath{\mathcal{K}}}

\DeclarePairedDelimiter{\ip}\langle\rangle
\DeclarePairedDelimiter{\nrm}\lVert\rVert
\theoremstyle{plain}
\newtheorem{theorem}{Theorem}[section]
\newtheorem{proposition}[theorem]{Proposition}
\newtheorem{lemma}[theorem]{Lemma}
\newtheorem{corollary}[theorem]{Corollary}

\theoremstyle{definition}

\newtheorem*{assumption*}{\assumptionnumber}
\providecommand{\assumptionnumber}{}
\makeatletter
\newenvironment{assumption}[2]
 {%
  \renewcommand{\assumptionnumber}{Assumption ($#2$#1)}%
  \begin{assumption*}%
  \protected@edef\@currentlabel{#1}%
 }
 {%
  \end{assumption*}
 }
\makeatother
\newcommand{\asref}[2]{$#2$\ref{#1}}

\numberwithin{equation}{section}

\begin{document}
\begin{frontmatter}

\title{Travelling wave solutions for fully discrete FitzHugh-Nagumo type equations with infinite-range interactions}
\journal{...}
\author[LDA]{W. M. Schouten-Straatman\corauthref{coraut}},
\corauth[coraut]{Corresponding author. }
\author[LDB]{H. J. Hupkes}
\address[LDA]{
  Mathematisch Instituut - Universiteit Leiden \\
  P.O. Box 9512; 2300 RA Leiden; The Netherlands \\ Email:  {\normalfont{\texttt{w.m.schouten@math.leidenuniv.nl}}}
}
\address[LDB]{
  Mathematisch Instituut - Universiteit Leiden \\
  P.O. Box 9512; 2300 RA Leiden; The Netherlands \\ Email:  {\normalfont{\texttt{hhupkes@math.leidenuniv.nl}}}
}

\date{\today}

\begin{abstract}
\singlespacing
We investigate the impact of spatial-temporal discretisation schemes on the dynamics of a class of reaction-diffusion equations that includes the FitzHugh-Nagumo system. For the temporal discretisation we consider the family of six backward differential formula (BDF) methods, which includes the well-known backward-Euler scheme. The spatial discretisations can feature infinite-range interactions, allowing us to consider neural field models. We construct travelling wave solutions to these fully discrete systems in the small time-step limit by viewing them as singular perturbations of the corresponding spatially discrete system. In particular, we refine the previous approach by Hupkes and Van Vleck for scalar fully discretised systems, which is based on a spectral convergence technique that was developed by Bates, Chen and Chmaj.
\end{abstract}

\begin{keyword}
\singlespacing
Travelling waves, FitzHugh-Nagumo system, singular perturbation, spatial-temporal discretisation.
\MSC 34K26, 39A06, 39A12, 47B39
\end{keyword}

\end{frontmatter}

\section{Introduction}\label{fdfhn:sec:introduction}

In this paper, we consider spatial-temporal discretisations of
a class of reaction-diffusion systems that contains the FitzHugh-Nagumo partial differential equation (PDE). This PDE is given by
\begin{equation}\label{fdfhn:FHNPDE}\begin{array}{lcl}u_t&=&u_{xx}+g(u;r)-w\\[0.2cm]
w_t&=&\rho(u-\gamma w).\end{array}\end{equation}
Here $g$ is the bistable, cubic nonlinearity $g(u;r)=u(1-u)(u-r)$ with $r\in(0,1)$, while $\rho>0$ and $\gamma>0$ are positive constants. In particular, our goal is to show that travelling waves for the system (\ref{fdfhn:FHNPDE}) persist under these spatial-temporal discretisations. As such, we contribute to the broad study of numerical schemes and their impact on the solutions under consideration, which has produced an immense quantity of literature. The main distinguishing feature is that we are interested in structures that persist for all time, while almost all of the studies in this area focus on finite time estimates.

\paragraph{Pulse propagation}
The system (\ref{fdfhn:FHNPDE}) was introduced in the 1960s \cite{FITZHUGH19662,FITZHUGH1966} as a simplification of the Hodgkin-Huxley equations, which were used to describe the propagation of spike signals through the nerve fibers of giant squids \cite{HODHUX1952}. After observing similar pulse solutions for the system (\ref{fdfhn:FHNPDE}) numerically \cite{FITZHUGH1968}, a more rigorous, analytical approach to understanding these pulse solutions turned out to be rather delicate. Indeed, many new tools have been developed, some even very recently, to construct these pulses and analyse their stability in various settings. These techniques include geometric singular perturbation theory \cite{CARP1977,HAST1976,JONESKOPLAN1991,JONES1984},
the variational principle \cite{chen2015traveling},
Lin's method \cite{KRUSANSZM1997,CARTER2015,CARTER2016},
and the Maslov index \cite{cornwell2017opening,cornwell2018existence}. Pulse solutions for the system (\ref{fdfhn:FHNPDE}) take the form
\begin{equation}
    \begin{array}{lcl}
         (u,w)(x,t)&=&(u_0,w_0)(x+c_0t) 
    \end{array}
\end{equation}
for some wavespeed $c_0$ and smooth wave profiles $u_0,w_0$ that satisfy the limits
\begin{equation}\label{fdfhn:pdelimits}
\begin{array}{lcl}
\lim\limits_{|\xi|\rightarrow\infty}(u_0,w_0)(\xi)&=&0.
\end{array}
\end{equation}

\paragraph{Spatially discrete systems} It is well-known that electrical pulses can only move through nerve fibres at appropriate speeds if the nerves are insulated with a myelin coating. This coating admits regularly spaced gaps at the so-called nodes of Ranvier \cite{RANVIER1878}. In fact, through a process called saltatory conduction, excitations of these nerves appear to jump from one node to the next \cite{LILLIE1925}. Since the FitzHugh-Nagumo PDE (\ref{fdfhn:FHNPDE}) does not take this discrete structure into account directly, it has been proposed \cite{EVVPD18} to, instead, model these phenenomena using a so-called lattice differential equation (LDE). For example, by applying a nearest-neighbour spatial discretisation to (\ref{fdfhn:FHNPDE}), we arrive at
\begin{equation}\label{fdfhn:finiterangeversion}\begin{array}{lcl}\dot{u}_j&=&\tau(u_{j+1}+u_{j-1}-2u_j)+g(u_j;r)-w_j\\[0.2cm]
\dot{w}_j&=&\rho[u_j-\gamma w_j],\end{array}\end{equation}
where the variable $j$ ranges over the lattice $\Z$. In the system (\ref{fdfhn:finiterangeversion}), the variable $u_j$ represents the potential at the $j$th node of the nerve fibre, while the variable $w_j$ describes a recovery component. Finally, we have $\tau \sim h^{-2}$, where $h>0$ is the distance between subsequent nodes. We emphasize that the time variable remains continuous.\\

Spatialy discrete travelling pulses for the system (\ref{fdfhn:finiterangeversion}) take the form
\begin{equation}\label{fdfhn:ldeansatz}
\begin{array}{lcl}
(u,w)_j(t)&=&(\overline{u}_0,\overline{w}_0)(j+\overline{c}_0t),
\end{array}
\end{equation}
for some wavespeed $\overline{c}_0$, again with the limits (\ref{fdfhn:pdelimits}). Plugging the Ansatz (\ref{fdfhn:ldeansatz}) into the LDE (\ref{fdfhn:finiterangeversion}) yields the functional differential equation of mixed type (MFDE)
\begin{equation}\begin{array}{lcl}\overline{c}_0\overline{u}_0'(\xi)&=&\tau[\overline{u}_0(\xi +1 )+\overline{u}_0(\xi -1)-2\overline{u}_0(\xi)]+g(\overline{u}_0(\xi);r)-\overline{w}_0(\xi)\\[0.2cm]
\overline{c}_0\overline{w}_0'(\xi)&=&\rho[\overline{u}_0(\xi)-\gamma \overline{w}_0(\xi)]\end{array}\end{equation}
in which $\xi=j+\overline{c}_0t$. In \cite{HJHFZHNGM,HJHSTBFHN}, Hupkes and Sandstede developed an infinite dimensional version of the exchange lemma to show that the system (\ref{fdfhn:finiterangeversion}) admits nonlinearly stable travelling pulse solutions. They relied heavily on the existence of exponential dichotomies for MFDEs, which were established in \cite{MPVL,harterich}. In addition, we established the existence and nonlinear stability of pulse solutions for a spatially periodic version of (\ref{fdfhn:finiterangeversion}) \cite{SCHPerExtensions} by building on a spectral convergence method developed by Bates, Chen and Chmaj \cite{BatesInfRange}. The spectral convergence method plays an important role in this paper as well and will be treated in more detail later on.\\

\paragraph{Infinite-range interactions} 
Neural field models aim to describe the dynamic behaviour of large networks of neurons. In neural networks, neurons interact with each other over large distances through their interconnecting nerve axons \cite{BRESS2014,BRESS2011,PINTO2001,SNEYD2005}. It has been proposed \cite[Eq. (3.31)]{BRESS2011} to capture these long distance interactions using an infinite-range version of the system (\ref{fdfhn:finiterangeversion}). To be concrete, we focus our discussion on the prototype system
\begin{equation}\begin{array}{lcl}\label{fdfhn:infiniterangeversion}\dot{u}_j&=&\tau\sum\limits_{m\in\Z_{>0}}e^{-m^2}[u_{j+m}+u_{j-m}-2u_j]+g(u_j;r)-w_j\\[0.2cm]
\dot{w}_j&=&\rho[u_j-\gamma w_j].\end{array}\end{equation}
This system can also be obtained directly from the PDE (\ref{fdfhn:FHNPDE}) by using an infinite-range spatial discretisation.\\

We emphasize that infinite-range interactions also arise naturally when considering discretisations of fractional Laplacians \cite{Ciaurri}. Indeed, such operators are intrinsically nonlocal and are used in many physical systems that feature nonstandard diffusion processes, such as amorphous semiconductors \cite{GU1996} and liquid crystals \cite{Ciuchi}.\\

Substituting the travelling pulse Ansatz (\ref{fdfhn:ldeansatz}) into (\ref{fdfhn:infiniterangeversion}) now yields the MFDE 
\begin{equation}\label{fdfhn:MFDE}\begin{array}{lcl}\overline{c}_0\overline{u}_0'(\xi)&=&\tau\sum\limits_{m\in\Z_{>0}}e^{-m^2}[\overline{u}_0(\xi +m )+\overline{u}_0(\xi -m)-2\overline{u}_0(\xi)]+g(\overline{u}_0(\xi);r)-\overline{w}_0(\xi)\\[0.2cm]
\overline{c}_0\overline{w}_0'(\xi)&=&\rho[\overline{u}_0(\xi)-\gamma \overline{w}_0(\xi)],\end{array}\end{equation}
which features infinitely many shifts. Since exponential dichotomies for MFDEs with infinitely many shifts have only been established very recently \cite{SCHExpDich}, the techniques used by Hupkes and Sandstede for the LDE (\ref{fdfhn:finiterangeversion}) have not yet been fully developed for the system (\ref{fdfhn:MFDE}). Instead, Faye and Scheel \cite{Faye2015} used a functional analytic approach to construct pulse solutions for the system (\ref{fdfhn:infiniterangeversion}). In addition, by applying the previously mentioned spectral convergence method, we were able to show that these pulses are nonlinearly stable \cite{HJHFHNINFRANGE} for $\tau\gg 1$, which corresponds to fine discretisations of the PDE (\ref{fdfhn:FHNPDE}). As of now, no comprehensive result has been found for the system (\ref{fdfhn:infiniterangeversion}).\\

\paragraph{Spatial-temporal discretisations} 
Our main goal here is to understand the impact of temporal discretisation schemes on the behaviour of travelling wave solutions of the system (\ref{fdfhn:infiniterangeversion}). This is a relatively novel area of study, although a handful of results have been established for scalar problems. For example, Bambusi, Faou, Gre\'ebert and J\'ez\'equel constructed solutions to fully discrete Schr\"odinger equations with Dirichlet or periodic spatial boundary conditions in \cite{bambusi2013existence,faou2015resonant}. Most other studies have focused on spatial-temporal discretisations of the Nagumo PDE
\begin{equation}\label{fdfhn:nagumopde}
\begin{array}{lcl}
u_t&=&u_{xx}+g(u;r),
\end{array}
\end{equation}
or, equivalently, temporal discretisations of the Nagumo LDE
\begin{equation}\label{fdfhn:nagumoLDE}\begin{array}{lcl}\dot{u}_j&=&\tau(u_{j+1}+u_{j-1}-2u_j)+g(u_j;r).\end{array}\end{equation}
The PDE (\ref{fdfhn:nagumopde}) and the LDE (\ref{fdfhn:nagumoLDE}) can be seen as scalar versions of the FitzHugh-Nagumo PDE (\ref{fdfhn:FHNPDE}) and LDE (\ref{fdfhn:finiterangeversion}) respectively.\\

The early works by Elmer and Van Vleck \cite{Elmer2006,EVV03,EVVJCP03} provided ad-hoc techniques to understand the impact of spatial-, temporal- and spatial-temporal discretisations of the PDE (\ref{fdfhn:nagumopde}) on the dynamics of travelling waves. In addition, Chow, Mallet-Paret and Shen \cite{chow1996dynamics} established the existence of travelling wave solutions to temporal discretisations of the LDE (\ref{fdfhn:nagumoLDE}) by considering Poincare return maps for the dynamics of this LDE. These results were later expanded by Hupkes and Van Vleck \cite{HJHBDF}, whose methods allowed them to address issues of uniqueness and parameter-dependence. Let us also mention the recent series of papers \cite{HJHADPGRID1,HJHADPGRID2,HJHADPGRID3} by Hupkes and Van Vleck, who study spatial discretisation schemes with an adaptive grid. That is, the authors consider a time dependent moving mesh method which aims to equidistribute the arclength of the solution under consideration.\\

In order to introduce the temporal discretisation schemes that we study in this paper, we briefly discuss the test problem 
\begin{equation}\label{fdfhn:testproblem}
    \begin{array}{lcl}
         \dot{v}&=&\lambda v
    \end{array}
\end{equation}
 with $\lambda<0$. Applying the forward-Euler discretisation scheme with time-step $\Delta t>0$ yields
 \begin{equation}
\begin{array}{lclcl}
v_{n+1}&=&v_n+\lambda\Delta tv_n&=&(1+\lambda \Delta t)v_n,
\end{array}
\end{equation}
where $n\in\Z$. Since a nontrivial solution of the test problem (\ref{fdfhn:testproblem}) converges to zero as $t\rightarrow\infty$, the convergence $v_n\rightarrow 0$ should also be enforced. However, this yields the restriction $0<\Delta t<2|\lambda|^{-1}$, which cannot be satisfied for all $\lambda<0$ for a fixed time-step $\Delta t>0$. In contrast, these issues do not occur for the backward-Euler discretisation scheme. For the test problem (\ref{fdfhn:testproblem}), this scheme yields
\begin{equation}
\begin{array}{lcl}
v_{n+1}&=&v_n+\lambda \Delta tv_{n+1},
\end{array}
\end{equation}
or equivalently
\begin{equation}
\begin{array}{lcl}
v_{n+1}&=&(1-\lambda\Delta t)^{-1}v_n.
\end{array}
\end{equation}
In particular, we see that $v_n\rightarrow 0$ for any value of $\lambda<0$ and time-step $\Delta t>0$. A numerical scheme is called $A(\alpha)$ stable if this property holds for all $\lambda$ in the wedge $\{z\in\C\setminus\{0\}:\mathrm{Arg}(-z)<\alpha\}$. We note that the backward-Euler discretisation is $A(\frac{\pi}{2})$ stable.\\

In fact, the backward-Euler discretisation scheme is one of six so-called backwards differentiation formula (BDF) methods. These BDF methods are all $A(\alpha)$ stable for various coefficients $0<\alpha\leq\frac{\pi}{2}$ and have several convenient analytical properties. For this reason, we have to chosen to focus on these temporal discretisation schemes in this paper. We do, however, emphasize that there are other stable discretisation schemes which we could have used, see for example \cite{HAIRER1991}. \\

Applied to the Nagumo system, the backward-Euler discretisation scheme yields the evolution
\begin{equation}\label{fdfhn:eq:ngm:BDF1}\begin{array}{lcl}\frac{1}{\Delta t}\Big[U_j(n\Delta t)-U_j\big((n-1)\Delta t\big)\Big]&=&\tau\big[U_{j+1}+U_{j-1}-2U_j\big](n\Delta t)+g\big(U_j(n\Delta t);r\big).\end{array}\end{equation} 
A travelling wave solution for the system (\ref{fdfhn:eq:ngm:BDF1}) with wavespeed $c$ takes the form
\begin{equation}
\begin{array}{lcl}
U_j(n\Delta t)&=&\Phi(j+nc\Delta t),
\end{array}
\end{equation}
with the limits
\begin{equation}
    \begin{array}{lclcl}
         \lim\limits_{\xi\rightarrow-\infty}\Phi(\xi)&=&0,\qquad \lim\limits_{\xi\rightarrow\infty}\Phi(\xi)&=&1.
    \end{array}
\end{equation}
As such, the travelling waves need to satisfy the system
\begin{equation}
\begin{array}{lcl}\label{fdfhn:discretenagumowaveequation}
\frac{1}{\Delta t}\big[\Phi(\xi)-\Phi(\xi-c\Delta t)\Big]&=&\tau\big[\Phi(\xi+1)+\Phi(\xi-1)-2\Phi(\xi)\big]+g(\Phi(\xi);r).
\end{array}
\end{equation}
Hupkes and Van Vleck showed \cite{HJHBDF} that, for sufficiently large, rational values of $M=(c \Delta t)^{-1}$, the system (\ref{fdfhn:eq:ngm:BDF1}) admits travelling wave solutions with wavespeed $c$. These travelling waves are constructed as perturbations of travelling wave solutions of the LDE (\ref{fdfhn:nagumoLDE}). The corresponding transition from the semi-discrete setting to the fully discrete setting is highly singular, since a derivative is replaced by a difference. The rationality of $M$ plays a key role here, as it ensures that the domain of the variable $\xi$ in the system (\ref{fdfhn:discretenagumowaveequation}) is a discrete subset of the real line. This restriction arises naturally in the analysis, since it ensures we can use finitely many interpolations to go from a fully discrete to a spatially discrete setting. \\

\paragraph{Spectral convergence} In order to analyse this singular perturbation, Hupkes and Van Vleck relied heavily on the previously mentioned spectral convergence method, which also plays an important role in \cite{SCHPerExtensions,HJHFHNINFRANGE,HJHADPGRID1,HJHADPGRID2,HJHADPGRID3}. This method was introduced in \cite{BatesInfRange} to construct travelling wave solutions to an infinite-range version of the Nagumo LDE (\ref{fdfhn:nagumoLDE}) in the continuum limit, i.e. when the discretisation distance $h\sim\tau^{-\frac{1}{2}}$ is sufficiently small. A key role in \cite{BatesInfRange} is reserved for the family of operators 
\begin{equation}\label{fdfhn:scalaroperator}\begin{array}{lcl}\mathcal{L}_hv(\xi)&=&c_0v'(\xi)-\frac{1}{h^2}\big[v(\xi+h)+v(\xi-h)-2v(\xi)\big]-g_U(\overline{u}_0(\xi);r)v(\xi),\end{array}\end{equation}
which arise as the linearization of the travelling wave MFDE corresponding to the LDE (\ref{fdfhn:nagumoLDE}) around the travelling wave solution $(c_0,\overline{u}_0)$ to the PDE (\ref{fdfhn:nagumopde}). The main question is what properties these operators inherit from their continuous counterpart
\begin{equation}\label{fdfhn:eq:op:L0}\begin{array}{lcl}\mathcal{L}_0v(\xi)&=&c_0v'(\xi)-v''(\xi)-g_U(\overline{u}_0(\xi);r)v(\xi).\end{array}\end{equation}
In particular, the authors in \cite{BatesInfRange} fixed a constant $\delta>0$ and used the invertibility of the operator $\mathcal{L}_0+\delta$ to establish the invertibility of the operator $\mathcal{L}_h+\delta$ for $h>0$ sufficiently small. Indeed, they considered weakly converging sequences $\{v_n\}$ and $\{w_n\}$ with $\mathcal{L}_hv_n+\delta v_n=w_n$ and tried to find a uniform (in $h$ and $\delta$) lower bound on the norm of $v_n'$ in terms of the norm of $w_n$. Such a lower bound prevents the limitless transfer of energy into oscillatory modes, a common concern when dealing with weakly converging sequences. The bistable nature of the nonlinearity $g$ was used to control the behaviour at $\pm\infty$, while the local $L^2$-norm can be bounded on the remaining compact set. We emphasize that this method requires a detailed understanding of the limiting operator $\mathcal{L}_0$.\\

In \cite{HJHBDF}, this method was lifted to the fully discrete Nagumo equation (\ref{fdfhn:discretenagumowaveequation}). Writing $M=\frac{p}{q}$ with $\mathrm{gcd}(p,q)=1$, the corresponding limiting operator resembles a $q$ times coupled version of the operator $\mathcal{L}_h$ given by (\ref{fdfhn:scalaroperator}). For $q=2$, this limiting operator takes the form
\begin{equation}\label{fdfhn:coupledoperator}
    \begin{array}{lcl}
         \mathcal{K}_q v(\zeta,\xi)&=&\overline{c}v'(\zeta,\xi)-\tau\big[v(\zeta+\frac{1}{2},\xi+1)+v(\zeta-\frac{1}{2},\xi-1)-2v(\zeta,\xi)\big]-g_U(\overline{u}(\xi);r)v(\zeta,\xi),
    \end{array}
\end{equation}
where $\overline{u}$ is the travelling wave solution of the LDE (\ref{fdfhn:nagumoLDE}) with wavespeed $\overline{c}$. Here the domain of the variables $\zeta$ and $\xi$ is given by $\zeta\in\{0,\frac{1}{2}\}$ and $\xi\in \R$, with the convention that $v(\zeta+1,\xi)=v(\zeta,\xi)$. Since the MFDE corresponding to (\ref{fdfhn:coupledoperator}) admits a comparison principle, the Fredholm properties of the operator $\mathcal{K}_q$ follow directly from the general results in \cite{HJHNEGDIF}. Hupkes and Van Vleck generalized the spectral convergence method to lift the Fredholm properties of the operator $\mathcal{K}_q$ to the operator
\begin{equation}
    \begin{array}{lcl}
         \mathcal{K}_Mv(\zeta,\xi)&=&\overline{c}M\big[v(\zeta,\xi)-v(\zeta,\xi-M^{-1})\big]\\[0.2cm]
         &&\qquad -\tau\big[v(\zeta+\frac{1}{2},\xi+1-\frac{1}{2}M^{-1})+v(\zeta-\frac{1}{2},\xi-1+\frac{1}{2}M^{-1})-2v(\zeta,\xi)\big]\\[0.2cm]
         &&\qquad -g_U\big(\overline{u}(\xi);r\big)v(\zeta,\xi),
    \end{array}
\end{equation}
in the regime $M\gg 1$,
again with $\zeta\in\{0,\frac{1}{2}\}$ and $\xi\in\frac{1}{2}M^{-1}\Z$. The operator $\mathcal{K}_M$ arises as the linearisation of the fully discrete system (\ref{fdfhn:discretenagumowaveequation}) around the travelling wave $\overline{u}$, using the additional $\zeta$ variable to ensure that all $\xi$-shifted arguments are multiples of $M^{-1}$ .\\

\paragraph{Results} In this paper, we consider reaction-diffusion LDEs such as (\ref{fdfhn:infiniterangeversion}) and replace the temporal derivative by one of the six BDF discretisation schems. For example, applying the backward-Euler method to (\ref{fdfhn:infiniterangeversion}), we arrive at the prototype system
\begin{equation}\label{fdfhn:infspatialtemporal}\begin{array}{lcl}\frac{1}{\Delta t}[U_j(n\Delta t)-U_j((n-1)\Delta t)]&=&\tau\sum\limits_{m=1}^\infty e^{-m^2}\big[U_{j+m}+U_{j-m}-2U_j\big](n\Delta t)\\[0.2cm]
&&\qquad +g(U_j(n\Delta t);r)-W_j(n\Delta t)\\[0.2cm]
\frac{1}{\Delta t}[W_j(n\Delta t)-W_j((n-1)\Delta t)]&=&\rho[U_j(n\Delta t)-\gamma W_j(n\Delta t)].\end{array}\end{equation}
Our main result states that systems such as (\ref{fdfhn:infspatialtemporal}) admit travelling wave solutions. To achieve this, we extend the spectral convergence method that was developed in \cite{HJHBDF} for scalar LDEs with finite-range spatial interactions to the current setting, which features multi-component systems with infinite-range interactions. This generalisation is far from trivial and requires several technical obstructions to be resolved.\\

 The first main obstacle is that the spectral convergence method hinges on the understanding of the corresponding limiting operator. Indeed, the analog of the operator $\mathcal{K}_q$ from (\ref{fdfhn:coupledoperator}) for our system (\ref{fdfhn:infspatialtemporal}) does not admit a comparison principle, since this is not available for FitzHugh-Nagumo type systems. As such, very limited  a-priori knowledge is available for this limiting operator, which forces us to prove many of its properties from scratch. For this, we mainly employ techniques from harmonic analysis.\\
 
 The second main obstacle is that the system setting introduces several cross-terms that need to be controlled. Several key techniques from our earlier works \cite{HJHFHNINFRANGE,SCHPerExtensions} concerning spatially discrete systems can be adjusted to handle these cross-terms in the present fully-discrete setting. However, several crucial points in the analysis still require 
 these terms to be handled
 with 
 special care. \\ 
 %these terms to be handled carefully. consideration we still need to be very careful when dealing with these cross-terms. \\
 
 The remaining obstacles are directly related to the infinite-range interactions, which introduce several convergence issues that need to be overcome. It also requires us to establish more refined estimates on the decay rates of solutions to our limiting MFDE. We achieve this by employing an explicit representation of the corresponding inverse linear operator that was first introduced in \cite{HJHFHNINFRANGE}.\\

\paragraph{Loss of uniqueness}
In \cite{HJHBDF}, Hupkes and Van Vleck extensively studied the uniqueness and parameter-dependence of the travelling wave solutions of (\ref{fdfhn:eq:ngm:BDF1}). The key observation is that the rationality of the variable $M=(c\Delta t)^{-1}$ breaks the translational symmetry in the travelling wave problem, potentially allowing a \textit{family} of solutions to exist. For example, one can apply an irrational phase shift to the continuous waveprofiles for (\ref{fdfhn:nagumoLDE}) that underlies the perturbation argument discussed above. In this fashion, one could construct a different fully discrete wave for the same detuning parameter value $r$ in the nonlinearity $g(\cdot;r)$. However, this is a very delicate issue. In particular, $M=(c\Delta t)^{-1}$ is fixed in the analysis, so additional work is required to obtain results for fixed time-steps $\Delta t>0$.\\

For the backward-Euler discretisation scheme, this nonuniqueness can be made fully rigorous. In particular, Hupkes and Van Vleck showed that, for a fixed time step $\Delta t>0$ both the $r(c)$ relation and the $c(r)$ relation can be multi-valued. In particular, for a fixed value of $c$ there can be multiple values of $r$ for which a solution to the system (\ref{fdfhn:eq:ngm:BDF1}) exists and vice-versa. This can be achieved by embedding the system (\ref{fdfhn:discretenagumowaveequation}) into an MFDE that admits a comparison principle, allowing it to be analysed using the techniques developed by Keener \cite{VL28} and Mallet-Paret \cite{MPB}.\\

By contrast, the $c(r)$ relation for travelling wave solutions to the PDE (\ref{fdfhn:nagumopde}) and the LDE (\ref{fdfhn:nagumoLDE}) are both single-valued. The same holds for the $r(c)$ relation, with the single exception that it can be multi-valued for (\ref{fdfhn:nagumoLDE}) in the special case $c=0$ \cite{EVV2005AppMath,HOFFMPcrys}. This reflects the well-known wave-pinning phenomenon caused by the broken translational symmetry of the lattice \cite{VL28,HOFFMPcrys,BEYN1990,mallet2001crystallographic,32,EVV2005Nonlin}.\\

In this paper we study the $r(c)$ and the $c(r)$ relation for a fully-discrete version of the FitzHugh-Nagumo system.
%corresponding to the system (\ref{fdfhn:infspatialtemporal}). 
For the corresponding %FitzHugh-Nagumo 
PDE (\ref{fdfhn:FHNPDE}) and LDE (\ref{fdfhn:finiterangeversion}), numerical evidence \cite{carter2018unpeeling,BachKwok} suggests that both these relations are at most $2$-valued. In addition, theoretical results \cite{CARTER2016} for this PDE usually yield a locally unique $r(c)$ relation. For the system (\ref{fdfhn:infspatialtemporal}) a comparison principle is not available, rendering a direct analysis similar to the one in \cite{HJHBDF} infeasible. Instead, we run several numerical simulations to investigate these issues. These computations indicate that both the $r(c)$ and the $c(r)$ relation are typically multi-valued. Indeed, the points $(r,c)$ points at which we were able to find solutions
%fully-discrete travelling waves 
appear to map onto a surface instead of a curve. That is, there exists an entire spectrum of travelling wave solutions with different wavespeeds to the same fully discrete system. \\

\paragraph{Acknowledgements.}
Both authors acknowledge support from the Netherlands Organization for Scientific Research (NWO) (grant 639.032.612). 

\section{Main result}
Our main goal is to study the impact of several important temporal discretisation schemes on travelling wave solutions of reaction-diffusion LDEs of the form
\begin{equation}\label{fdfhn:ditishetalgemeneprobleem}\begin{array}{lcl}\dot{U}_j&=&\tau\sum\limits_{m>0} \alpha_m[U_{j+m}+U_{j-m}-2U_j]+\mathcal{G}(U_j;r).\end{array}\end{equation}
This LDE is posed on the one-dimensional lattice $j\in\Z$, but may have multiple components in the sense that $U_j\in\R^d$ for some integer $d\geq 1$. We start by discussing the structural conditions that we impose on the LDE (\ref{fdfhn:ditishetalgemeneprobleem}) and its travelling wave solutions in \S \ref{fdfhn:subsec:LDE} respectively \S \ref{fdfhn:subsec:trvwavesLDE}. In \S \ref{fdfhn:subsec:fullydiscrete} we introduce the appropriate temporal discretisation schemes and formulate our main result. Finally, we discuss some numerical results concerning the nonuniqueness of the fully discrete travelling waves in \S \ref{fdfhn:subsec:numerical}.

\subsection{The spatially discrete system}\label{fdfhn:subsec:LDE}
Besides a handful of exceptions \cite{Faye2014,Faye2015,HJHFHNINFRANGE,BatesInfRange,SCHExpDich,gui2015traveling}, almost all results concerning LDEs of the form (\ref{fdfhn:ditishetalgemeneprobleem}) assume that only finitely many of the coefficients $\alpha_m$ in (\ref{fdfhn:ditishetalgemeneprobleem}) are nonzero. However, following \cite{HJHFHNINFRANGE,BatesInfRange}, we will impose the following much weaker conditions.

\begin{assumption}{}{\text{HS}1}\label{fdfhn:aannames} The coefficients $\{\alpha_m\}_{m\in\Z_{>0}}$ are diagonal $d\times d$ matrices and $\tau>0$ is a positive constant. There exists $1\leq d_{\mathrm{diff}}\leq d$ so that for each $1\leq i\leq d_{\mathrm{diff}}$ we have $\alpha_m^{(i,i)}\neq 0$ for some $m\in\Z_{>0}$, while $\alpha_n^{(j,j)}=0$ for all $n\in\Z_{>0}$ and all $d_{\mathrm{diff}}<j\leq d$. The coefficients $\{\alpha_m\}_{m\in\Z_{>0}}$ satisfy the bound 
\begin{equation}\begin{array}{lcl}\sum\limits_{m>0}|\alpha_m|e^{m\nu}&<&\infty\end{array}\end{equation}
for some constant $\nu>0$, as well as the identity
\begin{equation}\begin{array}{lcllcl}\sum\limits_{m>0}\alpha_m^{(i,i)} m^2&=&1\end{array}\end{equation}
for each $1\leq i\leq d_{\mathrm{diff}}$. Finally, the inequality
\begin{equation}\label{fdfhn:conditionAz}\begin{array}{lcl}A_i(z):=\sum\limits_{m>0}\alpha_m^{(i,i)}\Big(1-\cos(mz)\Big)&>& 0\end{array}\end{equation}
holds for all $z\in(0,2\pi)$ and all $1\leq i\leq d_{\mathrm{diff}}$. \end{assumption}

In particular, the diffusion matrices $\{\alpha_m\}_{m\in\Z_{>0}}$ only act directly on the first $d_{\mathrm{diff}}$ components of $U_j$. For example, for the FitzHugh-Nagumo LDE
\begin{equation}\label{fdfhn:fhnLDE}\begin{array}{lcl}\dot{u}_j&=&\tau\sum\limits_{m>0} \alpha_m[u_{j+m}+u_{j-m}-2u_j]+u_j(1-u_j)(u_j-r)-w_j\\[0.2cm]
\dot{w}_j&=&\rho\big[u_j-\gamma w_j\big],\end{array}\end{equation}
we have $d=2$ and $d_{\mathrm{diff}}=1$, while for the Nagumo LDE
\begin{equation}\label{fdfhn:ngmLDE}\begin{array}{lcl}\dot{u}_j&=&\tau\sum\limits_{m>0} \alpha_m[u_{j+m}+u_{j-m}-2u_j]+u_j(1-u_j)(u_j-r)\end{array}\end{equation}
we have $d=d_{\mathrm{diff}}=1$.\\

We note that (\ref{fdfhn:conditionAz}) is automatically satisfied if $\alpha_m^{(i,i)}\geq 0$ for all $m\in\Z_{>0}$ and $\alpha_1^{(i,i)}\neq 0$. The conditions in (\asref{fdfhn:aannames}{\text{HS}1}) ensure that for $\phi\in L^\infty(\R;\R)$ with $\phi''\in L^2(\R;\R)$ and $1\leq i\leq d_{\mathrm{diff}}$, we have the limit
\begin{equation}\begin{array}{lcl}\lim\limits_{h\downarrow 0}\enskip\nrm{\frac{1}{h^2}\sum\limits_{m>0}\alpha_m^{(i,i)}\big[\phi(\cdot+h m)+\phi(\cdot-hm)-2\phi(\cdot)\big]-\phi''}_{L^2(\R;\R)}=0;\end{array}\end{equation}
see \cite[Lem. 1]{BatesInfRange}. In particular, (\asref{fdfhn:aannames}{\text{HS}1}) ensures that (\ref{fdfhn:fhnLDE}) can be interpreted as the spatial discretisation of the FitzHugh-Nagumo PDE (\ref{fdfhn:FHNPDE}) on a grid with distance $h$, where $\tau=\frac{1}{h^2}$. Additional remarks concerning this assumption in the scalar case $d=1$ can be found in \cite[\S 1]{BatesInfRange}.\\

We now turn to the spatially homogeneous equilibrium solutions to (\ref{fdfhn:ditishetalgemeneprobleem}), which are roots of the nonlinearity $\mathcal{G}$. We will assume that there are two $r$-independent equilibria $P^{\pm}$, but emphasize that they are allowed to be identical.
\begin{assumption}{}{\text{HS}2}\label{fdfhn:aannames2} The parameter dependent nonlinearity $\mathcal{G}:\R^d\times (0,1)\rightarrow \R^d$ is $C^2$-smooth. There exist $P^{\pm}\in\R^d$ so that $\mathcal{G}(P^{\pm};r)=0$ holds for all $r\in(0,1)$. \end{assumption}

The temporal stability of these two equilibria $P^{\pm}$ plays an essential and delicate role in our analysis. Indeed, it does not suffice to simply require that the eigenvalues of $DG(P^\pm)$ have strictly negative real parts, see the proof of \cite[Lem. 4.6]{SCHPerExtensions} for details. Following \cite{SCHPerExtensions}, we consider two auxiliary assumptions on the triplet $(\mathcal{G},P^-,P^+)$ to address this issue. Recalling the constant $1\leq d_{\mathrm{diff}}\leq d$ from (\asref{fdfhn:aannames}{\text{HS}1}), we first write $D\mathcal{G}(U;r)$ in the block form
%Writing $DG(U)$ in the block form
\begin{equation}\label{fdfhn:blockstructure}
\begin{array}{lcl}
D\mathcal{G}(U;r)&=&\left(\begin{array}{ll}\mathcal{G}^{[1,1]}(U;r)&\mathcal{G}^{[1,2]}(U;r)\\ \mathcal{G}^{[2,1]}(U;r)& \mathcal{G}^{[2,2]}(U;r)
\end{array}\right)
\end{array}
\end{equation}
for any $U \in \R^{d}$ and $r\in(0,1)$, taking $D\mathcal{G}^{[1,1]}(U;r) \in \R^{d_{\mathrm{diff}} \times d_{\mathrm{diff}}}$.

\begin{assumption}{$3_{\overline{r}}$}{\text{HS}}\label{fdfhn:aannamesconstanten2} The triplet $(\mathcal{G},P^-,P^+)$ satisfies at least one of the following conditions.
\begin{enumerate}[label=(\alph*)]
    \item\label{fdfhn:aannamesconstantena} The matrices $-D\mathcal{G}(P^-;\overline{r})$ and $-D\mathcal{G}(P^+;\overline{r})$ are positive definite.
    \item\label{fdfhn:aannamesconstantenb} The matrices
$-\mathcal{G}^{[1,1]}(P^-;\overline{r}),-\mathcal{G}^{[1,1]}(P^+;\overline{r}),-\mathcal{G}^{[2,2]}(P^-;\overline{r})$ and $-\mathcal{G}^{[2,2]}(P^+;\overline{r})$ are positive definite.
In addition, there exists a constant $\Gamma > 0$
so that $\mathcal{G}^{[1,2]}(U;\overline{r})=-\Gamma \mathcal{G}^{[2,1]}(U;\overline{r})^T$ holds for all $U \in \R^{d}$. 
\end{enumerate}\end{assumption}

To illustrate these assumptions, we consider the nonlinearity
\begin{equation}
\label{fdfhn:eq:mr:def:g:fhn}
\begin{array}{lcl}
G_{\mathrm{fhn}}(u,w;r)  &=&\left(\begin{array}{l}
u(1-u)(u-r)-w\\[0.2cm]
\rho\big[u-\gamma w\big]
\end{array}\right)
\end{array}
\end{equation}
corresponding to the FitzHugh-Nagumo LDE (\ref{fdfhn:fhnLDE}). The triplet $(G_{\mathrm{fhn}},0,0)$ can easily be seen to satisfy (\asref{fdfhn:aannamesconstantenb}{\text{HS}3_{\overline{r}}}) with $\Gamma=\frac{1}{\rho}$. However, when $r>0$ is sufficiently small the Jacobian $DG_{\mathrm{fhn}}(0;r)$ has a pair of complex eigenvalues with negative real part. In this case, the condition (\asref{fdfhn:aannamesconstantena}{\text{HS}3_{\overline{r}}}) may fail to hold.\\

\subsection{Spatially discrete travelling waves}\label{fdfhn:subsec:trvwavesLDE}
Our final two assumptions for (\ref{fdfhn:ditishetalgemeneprobleem}) concern the existence and stability of travelling wave solutions that connect the equilibria $P^-$ and $P^+$. These solutions take the form
\begin{equation}\label{fdfhn:oudeansatz}\begin{array}{lcl}U_j(t)&=&\overline{U}_0(j+\overline{c}_0t)\end{array}\end{equation}
for some smooth profile $\overline{U}_0$ and nonzero wavespeed $\overline{c}_0$.  Substituting the Ansatz (\ref{fdfhn:oudeansatz}) into (\ref{fdfhn:ditishetalgemeneprobleem}) and writing $\xi=j+\overline{c}_0t$, we see that the pair $(\overline{c}_0,\overline{U}_0)$ must satisfy the travelling wave MFDE
\begin{equation}\label{fdfhn:oudetravellingwaveeq}\begin{array}{lcl}\overline{c}_0\overline{U}_0'(\xi)&=&\tau\sum\limits_{m>0}\alpha_m\Big[\overline{U}_0(\xi + m)+\overline{U}_0(\xi -m)-2\overline{U}_0(\xi)\Big]+\mathcal{G}\big(\overline{U}_0(\xi);r\big),\end{array}\end{equation}
together with the boundary conditions
\begin{equation}\label{fdfhn:oudeinitialconditions}\begin{array}{lcl}\lim\limits_{\xi\rightarrow\pm\infty}\overline{U}_0(\xi)&=&P^{\pm}.\end{array}\end{equation}

\begin{assumption}{$1_{\overline{r}}$}{\text{HW}}\label{fdfhn:existentie} There exists a waveprofile $\overline{U}_0$ and a wavespeed $\overline{c}_0\neq 0$ that solve the travelling wave MFDE (\ref{fdfhn:oudetravellingwaveeq}) for $r=\overline{r}$, together with the boundary conditions (\ref{fdfhn:oudeinitialconditions}).
 \end{assumption}

We now turn to the spectral stability of these travelling wave solutions. To this end, we introduce the operator $L_0:H^1(\R;\R^d)\rightarrow L^2(\R;\R^d)$ for the linearisation of (\ref{fdfhn:oudetravellingwaveeq}) around the travelling wave $\overline{U}_0$, which acts as
\begin{equation}\label{fdfhn:definitieLh}\begin{array}{lcl}L_0 &=&\overline{c}_0\partial_\xi-\Delta_0  -D_U\mathcal{G}\big(\overline{U}_0;\overline{r}\big).\end{array}\end{equation}
Here the operator $\Delta_0:L^2(\R;\R^d)\rightarrow L^2(\R;\R^d)$ is given by
\begin{equation}\label{fdfhn:definitieDelta0}\begin{array}{lcl}\Delta_0&=&\tau\sum\limits_{m>0}\alpha_m\Big[T_0^m+T_0^{-m}-2\Big],\end{array}\end{equation}
where
\begin{equation}\label{fdfhn:definitieT0}
    \begin{array}{lcl}
         (T_0\Phi)(\xi)&=&\Phi(\xi+1). 
    \end{array}
\end{equation}
In addition, we introduce the formal adjoint $L_0^*:H^1(\R;\R^d)\rightarrow L^2(\R;\R^d)$ of $L_0$ that acts as
\begin{equation}\begin{array}{lcl}L_0^* &=&-\overline{c}_0\partial_\xi-\Delta_0  -D_U\mathcal{G}\big(\overline{u}_0;\overline{r}\big)^T.\end{array}\end{equation}
We remark that the spectrum of $L_0$
is $2 \pi i c_0$-periodic on account of the identity
\begin{equation}\begin{array}{lcl}
\big(L_0 + \lambda \big) e^{2 \pi i \cdot }
&= &e^{2 \pi i \cdot} \big(L_0 + \lambda + 2 \pi i c_0 \big),\end{array}
\end{equation}
see \cite[Lem. 5.1]{HJHFHNINFRANGE}. We impose the following condition on the spectral properties of this operator $L_0$.
\begin{assumption}{$2_{\overline{r}}$}{\text{HW}}\label{fdfhn:spectralstability} There exist functions $\Phi_0^\pm\in H^1(\R;\R^d)$, together with a constant $\tilde{\lambda}>0$ so that the following properties hold for the LDE (\ref{fdfhn:ditishetalgemeneprobleem}) with $r=\overline{r}$.
\begin{enumerate}[label=(\roman*)]
\item We have the identity
\begin{equation}
\begin{array}{lcl}
\Phi_0^+&=&\overline{U}_0',
\end{array}
\end{equation}
together with the normalisation
\begin{equation}\label{fdfhn:eq:normalisation}
\begin{array}{lcl}
\ip{\Phi_0^+,\Phi_0^-}_{L^2(\R;\R^d)}&=&1.
\end{array}
\end{equation}

\item The spectrum of the operator $-L_0$ in the half-plane $\{z\in\C:\Re\,  z\geq-\tilde{\lambda}\}$ consists precisely of the points $2\pi im\overline{c}_0$ with $m\in\Z$, which are all eigenvalues of $L_0$. Moreover, we have the identities
\begin{equation}\begin{array}{lcl}\ker (L_0)&=&\span\{\Phi_0^+\}
\\[0.2cm]
&=&\{g\in L^2(\R;\R^d):\ip{g,\Psi}_{L^2(\R;\R^d)}=0\text{ for all }\Psi\in\Range(L_0^*)\}

\end{array}\end{equation}
and 
\begin{equation}\begin{array}{lcl}\ker (L_0^*)&=&\span\{\Phi_0^-\}
\\[0.2cm]
&=&\{g\in L^2(\R;\R^d):\ip{g,\Psi}_{L^2(\R;\R^d)}=0\text{ for all }\Psi\in\Range(L_0)\}.

\end{array}\end{equation}\end{enumerate}\end{assumption}

Recall that an eigenvalue $\lambda$ of a Fredholm operator $L$ is said to be \textit{simple} if the kernel of $L-\lambda$ is spanned by one vector $v$ and the equation $(L-\lambda)w=v$ does not have a solution $w$. Note that if $L$ has a formal adjoint $L^*$, this is equivalent to the condition that $\ip{v,w}\neq 0$ for all nontrivial $w\in \ker(L^*-\overline{\lambda})$. In particular, the normalisation (\ref{fdfhn:eq:normalisation}) implies that the eigenvalues $2\pi i\overline{c}_0\Z$ are all simple eigenvalues of $-L_0$. \\

For the FitzHugh-Nagumo system (\ref{fdfhn:fhnLDE}), the assumptions (\asref{fdfhn:existentie}{\text{HW}}) and (\asref{fdfhn:spectralstability}{\text{HW}}) are both satisfied for all sufficiently small discretisation distances $h>0$ and sufficiently small $\rho>0$, see \cite[Thm. 2.1, Thm. 2.2, Prop. 4.2]{HJHFHNINFRANGE}. If the shifts have finite-range, i.e. $\alpha_m=0$ for all sufficiently large $m$, then these assumptions are satisfied \cite[Thm. 1]{HJHFZHNGM}-\cite[Prop. 5.1]{HJHSTBFHN} for sufficiently small $\rho>0$ without any restriction on the discretisation distance $h$. There are, however, conditions on $r$ and $\gamma$ in both cases.\\

\subsection{The fully discrete system}\label{fdfhn:subsec:fullydiscrete}

We aim to approximate solutions to (\ref{fdfhn:ditishetalgemeneprobleem}) at discrete time intervals $t=n\Delta t$ by
\begin{equation}
\begin{array}{lcl}
U_j(n\Delta t)\sim W_j(n\Delta t).
\end{array}
\end{equation} 
We need to apply an appropriate discretisation scheme to the temporal derivative in (\ref{fdfhn:ditishetalgemeneprobleem}). Although there are many different approximation schemes available, we mainly focus on the six so-called BDF methods. These methods are based on interpolation polynomials of different degrees. In particular, the BDF method of order $k\in\{1,2,...,6\}$ 
approximates $U'$ in (\ref{fdfhn:ditishetalgemeneprobleem}) at $t=n\Delta t$ by first constructing an interpolating polynomial of degree
$k$ through the $k+1$ points $\{W((n-n')\Delta t)\}_{n'=0}^k$ and then computing the derivative of this polynomial at $W(n\Delta t)$. As such, the temporal discretisations of the LDE (\ref{fdfhn:ditishetalgemeneprobleem}) under consideration are of the form
\begin{equation}\label{fdfhn:ditishetnieuwealgemeneprobleem}
\begin{array}{lcl}
\beta_k^{-1}\frac{1}{\Delta t}\sum\limits_{n'=0}^k \mu_{n';k}W_j\big(n\Delta t-(k-n')\Delta t\big)&=&\tau\sum\limits_{m>0} \alpha_m[W_{j+m}(n\Delta t)+W_{j-m}(n\Delta t)-2W_j(n\Delta t)]\\[0.2cm]
&&\qquad +\mathcal{G}\big(W_j(n\Delta t);r\big).
\end{array}
\end{equation}
The coefficients $\beta_k$ and $\{\mu_{n;k}\}$ in (\ref{fdfhn:ditishetnieuwealgemeneprobleem}) are given implicitly by the identities 
\begin{equation}\label{fdfhn:definitiemun}
\begin{array}{rcl}
\sum\limits_{n=0}^k \mu_{n;k}v\big((n-k)\Delta t\big)&=&\sum\limits_{n'=1}^k[\partial^{n'}v](0),\\[0.2cm]
\beta_k&=& \sum\limits_{n=0}^k \mu_{n;k}(n-k),
\end{array}
\end{equation}
which must hold for any scalar function $v$. Here we have introduced the notation
\begin{equation}
\begin{array}{lcl}
[\partial v](n\Delta t)&=&v\big(n\Delta t\big)-v\big((n-1)\Delta t\big).
\end{array}
\end{equation}
This definition yields that $\sum\limits_{n=0}^k \mu_{n;k}=0$, which allows us to identify 
\begin{equation}
\begin{array}{lclcl}\beta_k&=& \sum\limits_{n=0}^k \mu_{n;k}(n-m)&=&\sum\limits_{n=1}^k\mu_{n;k}n.
\end{array}
\end{equation}
For convenience, the values of the coefficients $\beta_k$ and $\mu_{n;k}$ can be found in Table \ref{fdfhn:Tabelcoefficienten}. We note that the BDF method of order 1 is the well-known backward-Euler method. \\

\begin{table}\begin{center}
\begin{tabular}{|L|L|L|L|L|L|L| }
\hline 
\rule[-1ex]{0pt}{2.5ex} \mu_{n;k} & k=1 & k=2 & k=3 & k=4 & k=5 & k=6 \\ 
\hline 
\rule[-1.5ex]{0pt}{4.5ex} n=0 & -1 & \frac{1}{3} & -\frac{2}{11} & \frac{3}{25} & -\frac{12}{137} & \frac{10}{147} \\
\hline 
\rule[-1.5ex]{0pt}{4.5ex} n=1 & 1 & -\frac{4}{3} & \frac{9}{11} & -\frac{16}{25} & \frac{75}{137} & -\frac{72}{147} \\
\hline 
\rule[-1.5ex]{0pt}{4.5ex} n=2 & • & 1 & -\frac{18}{11} & \frac{36}{25} & -\frac{200}{137} & \frac{225}{147} \\
\hline 
\rule[-1.5ex]{0pt}{4.5ex} n=3 & • & • & 1 & -\frac{48}{25} & \frac{300}{137} & -\frac{400}{147} \\ 
\hline 
\rule[-1.5ex]{0pt}{4.5ex} n=4 & • & • & • & 1 & -\frac{300}{137} & \frac{450}{147} \\
\hline 
\rule[-1.5ex]{0pt}{4.5ex} n=5 & • & • & • & • & 1 & -\frac{360}{147} \\
\hline 
\rule[-1.5ex]{0pt}{4.5ex} n=6 & • & • & • & • & • & 1 \\
\hline 
\rule[-1.5ex]{0pt}{4.5ex} \beta_k & 1 & \frac{2}{3} & \frac{6}{11} & \frac{12}{25} & \frac{60}{137} & \frac{60}{147} \\[0.3cm]
\hline 
\end{tabular} \end{center}
\caption{The coefficients $\mu_{n;k}$ and $\beta_k$ associated to the BDF discretisation schemes as given by (\ref{fdfhn:definitiemun}).}
\label{fdfhn:Tabelcoefficienten}
\end{table}

Our main goal is to study travelling wave solutions to the fully discrete system (\ref{fdfhn:ditishetnieuwealgemeneprobleem}), utilizing our assumptions for the spatially discrete system (\ref{fdfhn:ditishetalgemeneprobleem}). Such solutions are given by the Ansatz
\begin{equation}\label{fdfhn:nieuweansatz}
\begin{array}{lcl}
W_j(n\Delta t)&=&\Phi(j+nc\Delta t),
\end{array}
\end{equation}
for some wave speed $c$ and profile $\Phi$ with the boundary conditions
\begin{equation}\label{fdfhn:nieuweboundaryvalues}
\begin{array}{lcl}
\Phi(\pm\infty)&=&P^{\pm},
\end{array}
\end{equation}
in a sense that we make precise below.\\

For notational convenience, we introduce the quantity $M=(c\Delta t)^{-1}$. Substituting the Ansatz (\ref{fdfhn:nieuweansatz}) into (\ref{fdfhn:ditishetnieuwealgemeneprobleem}) yields the system
\begin{equation}\label{fdfhn:ditishetnieuweprobleem}
\begin{array}{lcl}
c[\D_{k,M}\Phi](\xi)&=&\tau\sum\limits_{m>0} \alpha_m[\Phi(\xi+m)+\Phi(\xi-m)-2\Phi(\xi)] +\mathcal{G}\big(\Phi(\xi);r\big),
\end{array}
\end{equation}
for all $\xi$ that can be written as $\xi=n+jM^{-1}$ for $(j,n)\in\Z^2$. Here we have introduced the discrete derivatives
\begin{equation}\label{fdfhn:definitiondiscretederivatives}
\begin{array}{lcl}
[\D_{k,M}\Phi](\xi)&=&\beta_k^{-1}M\sum\limits_{n'=0}^k \mu_{n';k}\Phi\big(\xi-(k-n')M^{-1}\big),
\end{array}
\end{equation}
for $k\in\{1,2,...,6\}$. From \cite[eq. (2.13)]{HJHBDF} we obtain the useful estimate
\begin{equation}\label{fdfhn:DmMconvergeert}
\begin{array}{lcl}
|[\D_{k,M}\Phi](\xi)-\Phi'(\xi)|&\leq & C_l M^{-l}\sup_{-kM^{-1}\leq \theta\leq 0}|\Phi^{(l+1)}(\xi+\theta)|,
\end{array}
\end{equation}
for all integers $1\leq l\leq k$ and all $\Phi\in C^{l+1}(\R;\R^d)$, in which the constant $C_l\geq 1$ is independent of $k$, $\Phi$ and $M$. Indeed, this estimate shows that the regular derivative can be approximated by the discrete derivatives as the time step $\Delta t$ shrinks to zero. We emphasize that BDF discretisation schemes of order $k\geq 2$ do not allow for a comparison principle, even when the original LDE does allow for one. This is a consequence of the existence of coefficients $\mu_{n;k}>0$ that have $n<k$.\\

Most of our results, including our main theorem, require a restriction on the values of $M$ that are allowed. In particular, upon fixing an integer $q\geq 1$, we introduce the set
\begin{equation}
\begin{array}{lcl}
\mathcal{M}_{q}&=&\{\frac{p}{q}:p\in\N\text{ has }\mathrm{gcd}(p,q)=1\text{ and }p\geq q\}.
\end{array}
\end{equation}
Often, we introduce $M=\frac{p}{q}\in\mathcal{M}_{q}$, which implicitly defines the integer $p=p(M)=qM$. Moreover, we see that the natural domain for the values of $\xi$ in the system (\ref{fdfhn:ditishetnieuweprobleem}), as well as in the boundary conditions (\ref{fdfhn:nieuweboundaryvalues}), is precisely the set $p^{-1}\Z$.

\begin{theorem}\label{fdfhn:maintheorem} Assume that (\asref{fdfhn:aannames}{\text{HS}1}) and (\asref{fdfhn:aannames2}{\text{HS}2}) are satisfied and pick $\overline{r}$ in such a way that (\asref{fdfhn:aannamesconstanten2}{\text{HS}}), (\asref{fdfhn:existentie}{\text{HW}}) and (\asref{fdfhn:spectralstability}{\text{HW}}) are satisfied. Fix a pair of integers $1\leq k\leq 6$ and $q\geq 1$. Then there exist constants $M_*\gg 1$ and $\delta_r>0$ so that for any $M=\frac{p}{q}\in\M_{q}$ with $M\geq M_*$, there exist continuous functions
\begin{equation}
\begin{array}{lcl}
c_{M}:\R\times [\overline{r}-\delta_r,\overline{r}+\delta_r]&\rightarrow&\R,\\[0.2cm]
\overline{U}_{M}:\R\times [\overline{r}-\delta_r,\overline{r}+\delta_r]&\rightarrow&\ell^\infty(p^{-1}\Z;\R^d)
\end{array}
\end{equation}
that satisfy the following properties.
\begin{enumerate}[label=(\roman*)]
\item\label{fdfhn:item:thm1} For any $(\theta,r)\in\R\times [\overline{r}-\delta_r,\overline{r}+\delta_r]$, the pair $c=c_{M}(\theta,r)$ and $\overline{U}=\overline{U}_{M}(\theta,r)$ satisfies the system
\begin{equation}\label{fdfhn:ditishetprobleeminhoofdstelling}
\begin{array}{lcl}
c[\D_{k,M}\overline{U}](\xi)&=&\tau\sum\limits_{m>0} \alpha_m[\overline{U}(\xi+m)+\overline{U}(\xi-m)-2\overline{U}(\xi)]+\mathcal{G}\big(\overline{U}(\xi);r\big)
\end{array}
\end{equation}
for $\xi\in p^{-1}\Z$, together with the boundary conditions
\begin{equation}
\label{fdfhn:boundaryconditionsinhoofdstelling}
\begin{array}{lcl}
\lim\limits_{\xi\rightarrow\pm\infty,\xi\in p^{-1}\Z}\overline{U}(\xi)&=&P^{\pm}.
\end{array}
\end{equation}
\item\label{fdfhn:item:thm2} For any $(\theta,r)\in\R\times [\overline{r}-\delta_r,\overline{r}+\delta_r]$, the solution $\overline{U}=\overline{U}_{M}(\theta,r)$ admits the normalisation
\begin{equation}
\label{fdfhn:normalisationinhoofdstelling}
\begin{array}{lcl}
\sum\limits_{\xi\in p^{-1}\Z}\Big[\Big\langle\Phi_0^-(\xi+\theta),\overline{U}(\xi)-\overline{U}_0(\xi+\theta)\Big\rangle_{\R^d}\Big]&=&0.
\end{array}\end{equation}
\item\label{fdfhn:item:thm3} For any $(\theta,r)\in\R\times [\overline{r}-\delta_r,\overline{r}+\delta_r]$, we have the shift-periodicity
\begin{equation}
\label{fdfhn:shiftperiodicityinhoofdstelling}
\begin{array}{lcl}
c_{M}(\theta+p^{-1},r)&=&c_{M}(\theta,r),\\[0.2cm]
\overline{U}_{M}(\theta+p^{-1},r)(\xi)&=&\overline{U}_{M}(\theta,r)(\xi+p^{-1}).
\end{array}
\end{equation}
\end{enumerate}
In  addition,  there  exists $\delta>0$ such  that  the  following  holds  true.  Any  triplet $(c,\overline{U},\theta)\in \R\times \ell^\infty(p^{-1}\Z;\R^d)\times \R$ that satisfies (\ref{fdfhn:ditishetprobleeminhoofdstelling}) for some pair $(r,M)\in (0,1)\times \M_{q}$ with 
\begin{equation}
\begin{array}{lclclclcl}
|r-\overline{r}|&<&\delta,\qquad M&=&\frac{p}{q}&>&\delta^{-1}&\geq &M_*
\end{array}
\end{equation}
and also enjoys the estimate
\begin{equation}\label{fdfhn:uniquenessestimateinhoofdstelling}
\begin{array}{lcl}
p^{-1}\sum\limits_{\xi\in p^{-1}\Z}\Big[|\overline{U}(\xi)-\overline{U}_0(\xi+\theta)|^2+|\D_{k,M}\overline{U}(\xi)-\D_{k,M}\overline{U}_0(\xi+\theta)|^2\Big]&<&\delta^2,
\end{array}
\end{equation}
must actually satisfy $c=c_{M}(\ttheta,r)$ and $\overline{U}=\overline{U}_{M}(\ttheta,r)$ for some $\ttheta\in\R$.
\end{theorem}

The factor $p^{-1}$ in (\ref{fdfhn:uniquenessestimateinhoofdstelling}) is used to compensate the growing number of terms as $p\rightarrow \infty$. In particular, we can view this as a uniqueness result with respect to a scaled $L^2$-norm that will be specified later.

\subsection{Nonuniqueness and numerical examples}\label{fdfhn:subsec:numerical}
Fixing $r\in[\overline{r}-\delta_r,\overline{r}+\delta_r]$, $M=\frac{p}{q}\geq M_*$ and $\theta\in\R$, the travelling wave $(c_M(\theta,r),\overline{U}_M(\theta,r))$ is constructed as a perturbation of the travelling wave $(\overline{c}_0,\overline{U}_0(\cdot+\theta))$ on the domain $p^{-1}\Z$. Since the wave profiles $\overline{U}_0(\cdot+\theta)$ and $\overline{U}_0(\cdot+\theta+p^{-1})$ are simply translates of each other on this domain, the shift-periodicity (\ref{fdfhn:shiftperiodicityinhoofdstelling}) follows easily. However, it is not clear how, specifically, the travelling wave depends on $\theta$. Indeed, in \cite[\S 5]{HJHBDF}, Hupkes and Van Vleck show that it is reasonable to expect that the derivative $\partial_\theta c_M(\theta,\overline{r})$ is exponentially small in $M$. As such, it is unclear how to further analyse this dependence.\\

We emphasize that in general the travelling wave solution will not necessarily be unique, even up to translation. In particular, fixing $\theta\in(0,p^{-1})$, we note that the waves $\overline{U}_0$ and $\overline{U}_0(\cdot+\theta)$ are different on the domain $p^{-1}\Z$. One might be tempted to conclude that if $M$ is sufficiently large, the wave profiles $\overline{U}_M(0,r)$ and $\overline{U}_M(\theta,r)$ are different as well. However, a larger value of $M$ means that the grid $p^{-1}\Z$ becomes finer. In particular, since the travelling waves $\overline{U}_M(0,r)$ and $\overline{U}_M(\theta,r)$ are perturbations of the waves $\overline{U}_0$ and $\overline{U}_0(\cdot+\theta)$, it could be that these perturbations cancel out the difference between $\overline{U}_0$ and $\overline{U}_0(\cdot+\theta)$.\\

In addition, since the constant $M=(c\Delta t)^{-1}$ is fixed in the statement of Theorem \ref{fdfhn:maintheorem}, fluctuations in $c$ automatically lead to changes in $\Delta t$. This complicates our understanding of the fully discrete system for a fixed timestep $\Delta t>0$. Our main goal here is to show that the wavespeed $c$ and the detuning parameter $r$ do not depend on each other in a locally unique fashion, which is in major contrast to the corresponding continuous and semi-discrete systems. \\

However, the lack of a comparison principle for FitzHugh-Nagumo systems heavily complicates a direct analysis. As such, we have chosen to, instead, use numerical simulations to illustrate these phenomena. In particular, we focus on the backward-Euler discretisation of the FitzHugh-Nagumo MFDE, which takes the form
\begin{equation}\label{fdfhn:specificfhnLDE}\begin{array}{lcl}(h\Delta t)^{-1}[u(\xi)-u(\xi-c\Delta t)]&=& h^{-2}[u(\xi+1)+u(\xi-1)-2u(\xi)]+g(u(\xi);r)-w(\xi)\\[0.2cm]
(h\Delta t)^{-1}[u(\xi)-u(\xi-c\Delta t)]&=&\rho\big[u(\xi)-\gamma w(\xi)\big].\end{array}\end{equation}
Here we fix $\rho=0.01$, $\gamma=5$, $h=\frac{5}{8}$ and we let $g$ be the bistable nonlinearity
\begin{equation}
    \begin{array}{lcl}
         g(u;r)&=&u(1-u)(u-r). 
    \end{array}
\end{equation}
Upon fixing the timestep $\Delta t=2$, we repeatedly solved the system (\ref{fdfhn:specificfhnLDE}) with Neumann boundary conditions on the interval $[-80,80]$ for different values of the parameters $(c,r)\in \Q\times(0,1)$.\\

These simulations turned out to be rather delicate, since the quality of the initial condition heavily influenced whether a solution could be found. In many cases, the simulation returned the zero solution. Simply augmenting an extra nontriviality condition often produced no solution at all. In addition, the value of $c$ greatly determines the number of points $\xi\in\R$ for which the values $(u,v)(\xi)$ need to be determined. In particular, upon writing 
\begin{equation}\label{fdfhn:eq:cinterval}
    \begin{array}{lcl}
         c&=&\frac{q\Delta t}{p},
    \end{array}
\end{equation}
we needed to consider the points in the set $p^{-1}\Z\cap[-80,80]$, which rapidly grows in number as $p$ increases. We considered values of $c$ of the form (\ref{fdfhn:eq:cinterval}) for values of $p\in\{1,2,...,8\}$ and $q\in\{1,2,...,2p\}$ with $\mathrm{gcd}(p,q)=1$, while the values of $r$ were taken in $\frac{1}{100}\Z\cap(0,\frac{1}{5})$.\\
\begin{figure}
\centering
\subcaptionbox{\label{fdfhn:fig:num:cvsr}}
[.47\textwidth]{\includegraphics[trim={2cm 3cm 2cm 2.5cm},scale=0.28,clip]{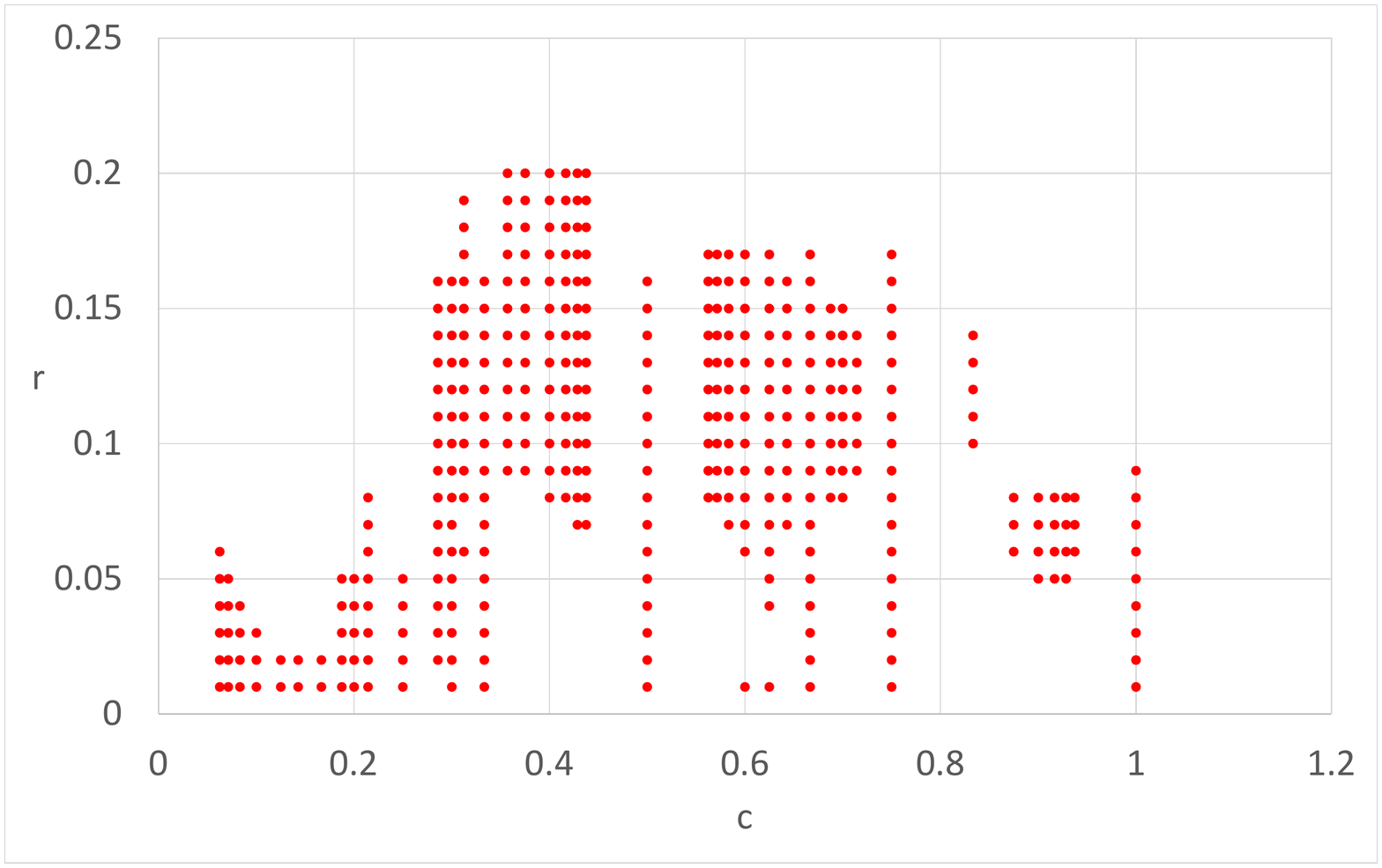}}
\subcaptionbox{\label{fdfhn:fig:num:plot}}
[.47\textwidth]{\includegraphics[trim={1cm 0cm 0 0cm},scale=0.23,clip]{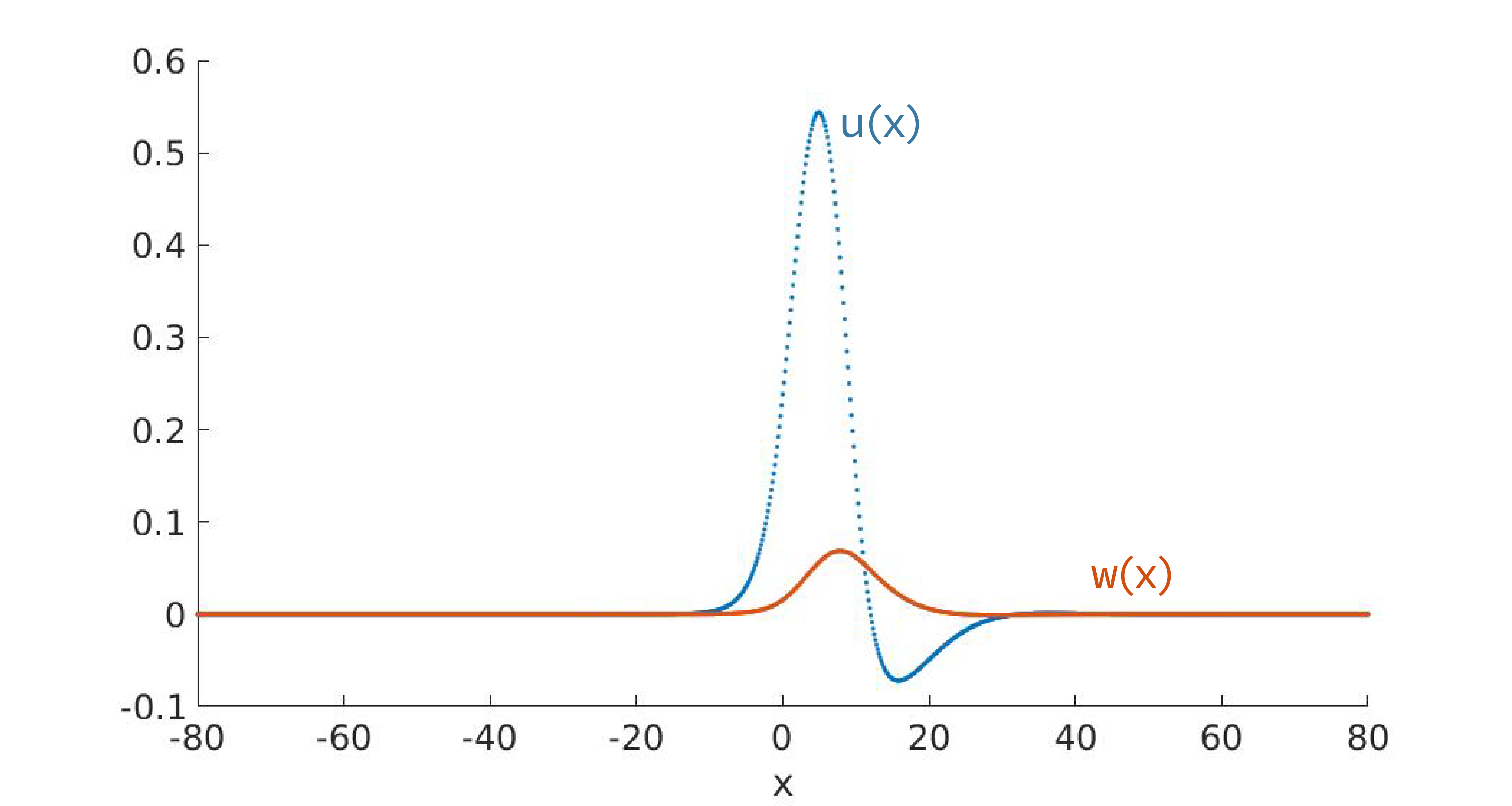}}
\caption{(a) Numerical computations of the pairs $(c,r)$ for which travelling wave solutions to the system (\ref{fdfhn:specificfhnLDE}) exist. We emphasize that there may be parameter values where we could not find a solution, but where a solution exists nonetheless. These simulations clearly show that the relationship $r(c)$ is multi-valued.
(b) A plot of one of the travelling waves found in this numerical procedure with $r=0.11$ and $c=0.3125$.}\label{fdfhn:fig:num}
\end{figure}

Figure \ref{fdfhn:fig:num:cvsr} depicts the pairs $(c,r)$ for which such a numerical solution could be found. It is highly likely that a solution still exists at some of the other parameter values that we investigated. In any case, our simulations clearly show that the parameters $c$ and $r$ depend on each other in an intricate fashion. In particular, our results suggest that travelling wave solutions to the system (\ref{fdfhn:specificfhnLDE}) are not unique, since we were able to find solutions with a range of different wavespeeds at the same value for $r$. We refer to \cite{carter2018unpeeling} and \cite{BachKwok} for the corresponding dependence for the FitzHugh-Nagumo PDE and LDE respectively. In both cases, this dependence is given by a curve in the $(c,r)$-plane that resembles the symbol $\cap$.

\section{Setup}\label{fdfhn:sectionsetup}

The fully discrete travelling wave equation (\ref{fdfhn:ditishetnieuweprobleem}) is a highly singular perturbation of the semi-discrete travelling wave MFDE (\ref{fdfhn:oudetravellingwaveeq}), which is the key complication for our analysis. In order to tackle this issue, we start by studying the linear operators that arise when linearizing the fully discrete travelling wave equation (\ref{fdfhn:ditishetnieuweprobleem}) around the semi-discrete travelling wave $(\overline{c}_0,\overline{U}_0)$. In particular, we define the linear expressions
\begin{equation}\label{fdfhn:definitieLhmM}\begin{array}{lcl}L_{k,M}\Phi(\xi)&=&\overline{c}_0[\D_{k,M}\Phi](\xi)-\Delta_0 \Phi(\xi) -D_U\mathcal{G}\big(\overline{U}_0(\xi)\big)\Phi(\xi).\end{array}\end{equation}
Our aim is to establish that the operators $L_{k,M}$ inherit several useful properties from the operator $L_0$ defined in (\ref{fdfhn:definitieLh}) in the small timestep limit $\Delta t\downarrow 0$.\\

In this section we summarize and adept the setup from \cite{HJHBDF}, sticking to the same notation as much as possible. In order to formulate our results, we need to define several function spaces. For any $\eta\in\R$, we write 
\begin{equation}\label{fdfhn:eq:bcetadef}
\begin{array}{lcl}
BC_\eta(\R;\R^d)&=&\{F\in C(\R;\R^d)\, |\, \sup_{\xi\in\R}e^{-\eta|\xi|}|F(\xi)|<\infty\},\\[0.2cm]
BC_\eta^1(\R;\R^d)&=&\{F\in C^1(\R;\R^d)\, |\, \sup_{\xi\in\R}e^{-\eta|\xi|}[|F(\xi)|+|F'(\xi)|]<\infty\}.
\end{array}
\end{equation}
In addition, given a Hilbert space $H$ and any $\mu>0$, we define the corresponding sequence space
\begin{equation}
\begin{array}{lcl}
\ell_\mu^2(H)&=&\{v:\mu^{-1}\Z\rightarrow H\, |\, \nrm{v}_{\ell_\mu^2(H)}:=\ip{v,v}_{\ell_\mu^2(H)}^\frac{1}{2}<\infty\},
\end{array}
\end{equation}
which is a Hilbert space equipped with the inner product
\begin{equation}\begin{array}{lcl}
\ip{v,v}_{\ell_\mu^2(H)}&=&\mu^{-1}\sum\limits_{\xi\in \mu^{-1}\Z}\ip{v(\xi),w(\xi)}_H.
\end{array}
\end{equation}

For now, we fix two integers $q\geq 1$ and $1\leq k\leq 6$, together with a constant $M=\frac{p}{q}\in \mathcal{M}_{q}$. To streamline our notation, we write $\mathcal{Y}_M$ to refer to the space $\ell_{p}^2(\R^d)$, i.e.,
\begin{equation}\label{fdfhn:eq:defYM}
\begin{array}{lclcl}
\mathcal{Y}_M&=&\ell_{p}^2(\R^d),\qquad \ip{\Phi,\Psi}_{\mathcal{Y}_M}&=&\ip{\Phi,\Psi}_{\ell_{p}^2(\R^d)}.
\end{array}
\end{equation}
Moreover, we introduce the space $\mathcal{Y}_{k,M}^1$, which differs from $\Y_M$ only by its inner product. To be more precise, we write
\begin{equation}\label{fdfhn:eq:defYkm1}
\begin{array}{lcl}
\mathcal{Y}_{k,M}^1&=&\ell_{p}^2(\R^d),\\[0.2cm]
\ip{\Phi,\Psi}_{\mathcal{Y}_{k,M}^1}&=&\ip{\Phi,\Psi}_{\ell_{p}^2(\R^d)}+\ip{\D_{k,M}\Phi,\D_{k,M}\Psi}_{\ell_{p}^2(\R^d)}.
\end{array}
\end{equation}

In addition, for $f\in BC_{-\eta}(\R;\R^d)$ with $\eta>0$, we write $\pi_{\mathcal{Y}_M}$ for the sequence
\begin{equation}\label{fdfhn:definitionpiym}
\begin{array}{lcl}
\big[\pi_{\mathcal{Y}_M}f\big](\xi)&=&f(\xi),\qquad \xi\in p^{-1}\Z.
\end{array}
\end{equation}
If moreover $f\in BC_{-\eta}^1(\R;\R^d)$ and we wish to be explicit, we often write $\pi_{\Y_{k,M}^1}f$ to refer to the restriction (\ref{fdfhn:definitionpiym}). The restriction operators $\pi_{\mathcal{Y}_M}$ and $\pi_{\Y_{k,M}^1}$ are bounded, see Lemma \ref{fdfhn:Lemma31BDF}.\\

We can now consider the operators $L_{k,M}$ appearing in (\ref{fdfhn:definitieLhmM}) as bounded linear maps 
\begin{equation}
\begin{array}{llcl}
L_{k,M}:& {\Y}_{k,M}^1&\rightarrow &  {\Y}_M.
\end{array}
\end{equation}
Our goal is to define new sequence spaces, which allow us to pass to the limit $M\rightarrow \infty$ in a controlled fashion. The basic idea is to use $L^2$-interpolants for functions in $\Y_M$ and $H^1$-interpolants for functions in $\Y_{k,M}^1$, so that the sequences in these spaces can be compared regardless of the different values of $M$. The main difficulty is to control terms of the form $v(\xi+p^{-1})-v(\xi)$ for $v\in \Y_{k,M}^1$ with $M=\frac{p}{q}$, which is impossible to extract solely from the behaviour of $\D_{k,M}v$. \\

To tackle this issue, we need to perform $q$ separate interpolations. Each of these interpolations must bridge a gap of size $M^{-1}=\frac{q}{p}$. In particular, upon fixing an integer $q\geq 1$ and writing
\begin{equation}
\begin{array}{lcl}
\Z_q&=&\{0,1,2,...,q\},\\[0.2cm]
\Z_q^\circ &=&\{1,2,...,q-1\},
\end{array}
\end{equation}
we introduce the space
\begin{equation}
\begin{array}{lcl}
\ell_{q,\perp}^2&=&\{\Phi:q^{-1}\Z_{q}\rightarrow \R^d\},
\end{array}
\end{equation}
equipped with the inner product
\begin{equation}
\begin{array}{lcl}
\ip{\Phi,\Psi}_{\ell_{q,\perp}^2}&=&q^{-1}\Big[\frac{1}{2}\Phi(0)\Psi(0)+\frac{1}{2}\Phi(1)\Psi(1)+\sum\limits_{\zeta\in q^{-1}\Z_q^\circ}\Phi(\zeta)\Psi(\zeta)\Big].
\end{array}
\end{equation}
Upon introducing the notation $\Phi(\zeta,\xi)=[\Phi(\xi)](\zeta)$ for $\Phi\in \ell_M^2(\ell_{q,\perp}^2)$ with $\zeta\in q^{-1}\Z_q$ and $\xi\in M^{-1}\Z$, we define the space
\begin{equation}\label{fdfhn:eq:defHM}
\begin{array}{lcl}
\mathcal{H}_M&=&\{\phi\in \ell_M^2(\ell_{q,\perp}^2):\Phi(1,\xi)=\Phi(0,\xi+M^{-1})\text{ for all }\xi\in M^{-1}\Z\},
\end{array}
\end{equation}
equipped with the inner product
\begin{equation}
\begin{array}{lcl}
\ip{\Phi,\Psi}_{\mathcal{H}_M}&=&M^{-1}\sum\limits_{\xi\in M^{-1}\Z}\ip{\Phi(\cdot,\xi),\Psi(\cdot,\xi)}_{\ell_{q,\perp}^2}.
\end{array}
\end{equation}

For any $\eta>0$ and any $f\in BC_{-\eta}(\R;\R^d)$, we now write $\pi_{\mathcal{H}_M}f\in\mathcal{H}_M$ for the function
\begin{equation}
\begin{array}{lcl}
[\pi_{\mathcal{H}_M}f](\zeta,\xi)&=&f(\xi+\zeta M^{-1}),\qquad \zeta\in q^{-1}\Z_q,\qquad \xi\in M^{-1}\Z.
\end{array}
\end{equation}
We extend the operators $\D_{k,M}$ to $\mathcal{H}_M$ by writing
\begin{equation}\label{fdfhn:eq:defHkM1}
\begin{array}{lcl}
[\D_{k,M}\Phi](\zeta,\xi)&=&[\D_{k,M}\Phi(\zeta,\cdot)](\xi).
\end{array}
\end{equation}
Note that these operators act only on the second component of $\Phi$. This allows us to define our final space 
\begin{equation}
    \begin{array}{lcl}
         \mathcal{H}_{k,M}^1&=&\mathcal{H}_M,
    \end{array}
\end{equation}
equipped with the inner product
\begin{equation}
\begin{array}{lcl}
\ip{\Phi,\Psi}_{\mathcal{H}_{k,M}^1}&=&\ip{\Phi,\Psi}_{\mathcal{H}_M}+\ip{\D_{k,M}\Phi,\D_{k,M}\Psi}_{\mathcal{H}_M}.
\end{array}
\end{equation}

In fact, we can relate the spaces $\H_M$ and $\H_{k,M}^1$ to the spaces defined earlier. To see this, we define the isometries
\begin{equation}
\begin{array}{lclcl}
\mathcal{J}_M:\mathcal{Y}_M&\rightarrow &\mathcal{H}_M,\qquad 
\mathcal{J}_{k,M}^1:\Y_{k,M}^1&\rightarrow &\mathcal{H}_{k,M}^1,
\end{array}
\end{equation}
for $M=\frac{p}{q}\in\mathcal{M}_{q}$, which both act as
\begin{equation}\label{fdfhn:definitieJM}
\begin{array}{lclcl}
[\mathcal{J}_M\Phi](\zeta,\xi)&=&[\mathcal{J}_{k,M}^1\Phi](\zeta,\xi)&=&\Phi(\xi+M^{-1}\zeta),
\end{array}
\end{equation}
for $\zeta\in q^{-1}\Z_q$ and $\xi\in M^{-1}\Z$, see Lemma \ref{fdfhn:Lemma35BDF}. Note that $\pi_{\mathcal{H}_M}=\mathcal{J}_M\pi_{\mathcal{Y}_M}$.\\

Our goal is to interpret $L_{k,M}$ as a map from $ {\H}_{k,M}^1$ into $ {\H}_M$. To this end, we pick $n\in\Z$ and $0<\vartheta\leq 1$ in such a way that
\begin{equation}\label{fdfhn:definitionrotation}
\begin{array}{lcl}
1&=&(n+\vartheta)M^{-1}.
\end{array}
\end{equation}
Since $M=\frac{p}{q}\in\mathcal{M}_{q}$, we see that $\vartheta=\frac{p-nq}{q}$, which yields
\begin{equation}
\begin{array}{lcl}
nM^{-1}&=&1-\vartheta M^{-1},\qquad \vartheta\in q^{-1} \Z_q\setminus\{0\}.
\end{array}
\end{equation}
In fact, because $\mathrm{gcd}(p,q)=1$, it follows that $\mathrm{gcd}(p,\vartheta q)=1$.\\

With these preparations in hand, we now write $\mathcal{K}_{k,M}: {\H}_{k,M}^1\rightarrow {\H}_M$ for the linear operator that acts as
\begin{equation}
\begin{array}{lcl}
\left[\mathcal{K}_{k,M}\Phi\right](\zeta,\xi)&=& \overline{c}_0[\D_{k,M}\Phi](\zeta,\xi)-\big[\Delta_{M}\Phi\big](\zeta,\xi)-D_U\mathcal{G}\big(\overline{U}_0(\xi+\zeta M^{-1});\overline{r}\big)\Phi(\zeta,\xi),
\end{array}
\end{equation}
for $\zeta\in q^{-1}\Z_q$ and $\xi\in M^{-1}\Z$. Here the operator $\Delta_M$ is given by
\begin{equation}\label{fdfhn:definitionDeltahM2}
\begin{array}{lcl}
\Delta_{M}&=&\tau\sum\limits_{m>0}\alpha_m\Big[  {T}_M^m+ {T}_M^{-m}-2\Big],
\end{array}
\end{equation}
where we have introduced the twist operator $T_M:\H_M\rightarrow\H_M$ that acts as
\begin{equation}\label{fdfhn:twistoperator}
\begin{array}{lcl}
[T_M\Phi](\zeta,\xi)&=&\Phi(\zeta+\vartheta,\xi+nM^{-1}),
\end{array}
\end{equation}
taking into account the convention
\begin{equation}\label{fdfhn:convention}
\begin{array}{lcl}
\Phi(\zeta\pm 1,\xi)&=&\Phi(\zeta,\xi\pm M^{-1}).
\end{array}
\end{equation}
In particular, we see that the shift $\vartheta$ acts as a rotation number, connecting the different components of $\Phi$ in the $\zeta$-direction. The inequality
\begin{equation}
\begin{array}{lcl}
\ip{\Delta_{M}\Phi,\Phi}_{ {\H}_M}&\leq &0
\end{array}
\end{equation}
for $\Phi\in {\H}_M$ is almost trivial to verify in the finite-range setting, but turns out to be much harder to establish when dealing with infinite-range interactions; see Lemma \ref{fdfhn:Lemma313BDF}.\\

Finally, we introduce the notation
\begin{equation}
\begin{array}{lcl}
D\mathcal{G}\big(\pi_{\H_M}\overline{U}_0;\overline{r}\big): {\H}_M&\rightarrow &  {\H}_M
\end{array}
\end{equation}
to refer to the multiplication operator
\begin{equation}
\begin{array}{lcl}
[D\mathcal{G}\big(\pi_{\H_M}\overline{U}_0;\overline{r}\big)\Phi](\zeta,\xi)&=&D_U\mathcal{G}\big(\overline{U}_0(\xi+\zeta M^{-1});\overline{r}\big)\Phi(\zeta,\xi).
\end{array}
\end{equation}
In fact, it is easy to see that
\begin{equation}\label{fdfhn:Eq348BDF}
\begin{array}{lcl}
\mathcal{K}_{k,M} {\J}_{k,M}^1&=& {\J}_ML_{k,M},
\end{array}
\end{equation}
which shows that $\mathcal{K}_{k,M}$ and $L_{k,M}$ are equivalent.\\

Since the operator $\K_{k,M}$ is not self-adjoint, we need to introduce the formal adjoint $\K_{k,M}^*: {\H}_{k,M}^1\rightarrow {\H}_M$ of $\K_{k,M}$ by writing
\begin{equation}
\begin{array}{lcl}
\mathcal{K}_{k,M}^*\Phi&=&\overline{c}_0[\D_{k,M}^*\Phi]-\Delta_{M}\Phi-D\mathcal{G}\big(\pi_{\H_M}\overline{U}_0;\overline{r}\big) ^T\Phi,
\end{array}\end{equation}
in which we have defined
\begin{equation}
\begin{array}{lcl}
[\D_{k,M}^*\Phi](\zeta,\xi)&=&\beta_k^{-1}M\sum\limits_{n'=0}^k \mu_{n';k}\Phi(\xi+(k-n')M^{-1}).
\end{array}
\end{equation}

Moreover, we introduce the space
\begin{equation}
\begin{array}{lcl}
\ell_{q,\perp;\infty}^2&=&\{\phi\in\ell_{q,\perp}^2:\phi(1)=\phi(0)\},
\end{array}
\end{equation}
together with the map
\begin{equation}
\begin{array}{lcl}
[\pi_{\perp}f](\zeta,\xi)&=&f(\xi),\qquad \zeta\in q^{-1}\Z_q,\qquad \xi\in\R,
\end{array}
\end{equation}
which constructs a function $\pi_{\perp}f\in L^2(\R,\ell_{q,\perp;\infty}^2)$ from a function $f\in L^2(\R;\R^d)$.\\

Taking the limit $M\rightarrow \infty$, while keeping $\vartheta$ and $q$ fixed as in (\ref{fdfhn:definitionrotation}), we see that $\K_{k,M}$ and $\K_{k,M}^*$ formally approach the limiting operators
\begin{equation}
\begin{array}{lcl}
\overline{\K}_{q,\vartheta}: {H}^1(\R,\ell_{q,\perp;\infty}^2)&\rightarrow & {L}^2(\R,\ell_{q,\perp;\infty}^2),\\[0.2cm]
\overline{\K}_{q,\vartheta}^*: {H}^1(\R,\ell_{q,\perp;\infty}^2)&\rightarrow & {L}^2(\R,\ell_{q,\perp;\infty}^2),
\end{array}
\end{equation}
that act as
\begin{equation}
\begin{array}{lcl}
\overline{\K}_{q,\vartheta}\Theta&=& \overline{c}_0\partial_\xi \Theta-\Delta_{q,\vartheta}\Theta-D\mathcal{G}\big(\pi_{\H_M}\overline{U}_0;\overline{r}\big) \Theta,\\[0.4cm]

\overline{\K}_{q,\vartheta}^*\Theta&=& -\overline{c}_0\partial_\xi \Theta-\Delta_{q,\vartheta}\Theta-D\mathcal{G}\big(\pi_{\H_M}\overline{U}_0;\overline{r}\big)^T \Theta.
\end{array}\end{equation}
Here the operator $\Delta_{q,\vartheta}$ is given by
\begin{equation}\label{fdfhn:def:deltaqtheta}
\begin{array}{lcl}
\Delta_{q,\vartheta}&=&\tau\sum\limits_{m>0}\alpha_m\Big[T_{q,\vartheta}^m +T_{q,\vartheta}^{-m}-2 \Big],
\end{array}
\end{equation}
in which we have introduced the twist operator
\begin{equation}\label{fdfhn:def:Ttheta}
    \begin{array}{lcl}
         \big[T_{q,\vartheta}\Theta\big](\zeta,\xi)&=&\Theta(\zeta+\vartheta,\xi+1),
    \end{array}
\end{equation}
for $\zeta\in q^{-1}\Z_q$ and $\xi\in\R$. In the same spirit as (\ref{fdfhn:convention}), we here make the convention $\Phi(\zeta+1,\xi)=\Phi(\zeta,\xi)$. Notice that the limiting operator $\overline{\K}_{q,\vartheta}$ reduces to the operator $L_0$ defined in (\ref{fdfhn:definitieLh}) for $\zeta$-independent functions.\\

\section{The limiting system}\label{fdfhn:section:limiting}

Our goal here is to exploit our understanding of the operator $L_0$ in order to determine the Fredholm properties of the limiting operator $\overline{\K}_{q,\vartheta}$. Due to the lack of a comparison principle we cannot immediately appeal to a general Frobenius-Peron-type result as was possible in \cite{HJHBDF}. The theory in this section aims to fill these gaps and can be considered the key technical contribution of this paper. We collect the main results in the following Proposition, which plays an essential role in Lemma \ref{fdfhn:Lemma316BDF} below.
\begin{proposition}[{cf. \cite[Lem. 3.6]{HJHBDF}}]\label{fdfhn:fredholmeigenschappenkh} Assume that (\asref{fdfhn:aannames}{\text{HS}1}) and (\asref{fdfhn:aannames2}{\text{HS}2}) are satisfied and pick $\overline{r}$ in such a way  that (\asref{fdfhn:aannamesconstanten2}{\text{HS}}), (\asref{fdfhn:existentie}{\text{HW}}) and (\asref{fdfhn:spectralstability}{\text{HW}}) are satisfied.  Fix an integer $q\geq 1$, together with a constant $\vartheta\in q^{-1}\Z_q$ that has $\mathrm{gcd}(\vartheta q,q)=1$. Then the operators $\overline{\K}_{q,\vartheta}$ and $\overline{\K}_{q,\vartheta}^*$ are both Fredholm operators with index $0$ and we have the identities
\begin{equation}
\begin{array}{lclccl}
\ker(\overline{\K}_{q,\vartheta})&=&\span\{\pi_{\perp}\Phi_0^+\}
,\qquad 
\ker(\overline{\K}_{q,\vartheta}^*)&=&\span\{\pi_{\perp}\Phi_0^-\}

.\end{array}
\end{equation}
Moreover, recalling the constant $\tilde{\lambda}$ appearing in (\asref{fdfhn:spectralstability}{\text{HW}}), the operator $\overline{\K}_{q,\vartheta}+\lambda$ is invertible for each $\lambda\in\C$ that has $\Re\, \lambda\geq-\tilde{\lambda}$
and $\lambda\notin2\pi i\overline{c}_0 q^{-1}\Z$. Finally, there exists constants $C>0$ and $\delta_0>0$ so that for each $0<\delta<\delta_0$ and each $\Theta\in L^2(\R,\ell_{q,\perp;\infty}^2)$ we have the bound
\begin{equation}\label{fdfhn:batesafschattingKh}
\begin{array}{lcl}
\nrm{[\overline{\K}_{q,\vartheta}+\delta]^{-1}\Theta}_{H^1(\R,\ell_{q,\perp;\infty}^2)}&\leq &C\Big[\nrm{\Theta}_{L^2(\R,\ell_{q,\perp;\infty}^2)}+\frac{1}{\delta}|\ip{\Theta,\pi_{\perp}\Phi_0^-}_{L^2(\R,\ell_{q,\perp;\infty}^2)}|\Big].
\end{array}
\end{equation}
\end{proposition}

The first step towards proving Proposition \ref{fdfhn:fredholmeigenschappenkh} is to find the eigenvalues of the operator $\overline{\K}_{q,\vartheta}$. After that, we will focus on the essential spectrum of this operator. The idea behind the proof of Lemma \ref{fdfhn:kernuitrekenen} below can best be illustrated by considering the case $q=2$. In this case, we have $\vartheta=\frac{1}{2}$, together with
\begin{equation}
    \begin{array}{lcl}
         \big[T_{2,\frac{1}{2}}\Theta\big](\zeta,\xi)&=&\Theta\big(\zeta+\frac{1}{2},\xi+1\big).
    \end{array}
\end{equation}
Upon writing
\begin{equation}
    \begin{array}{lclcl}
         [\Pi_0\Theta](\xi)&:=&\Theta(0,\xi)+\Theta(\frac{1}{2},\xi),\qquad [\Pi_1\Theta](\xi)&:=&\Theta(0,\xi)-\Theta(\frac{1}{2},\xi),
    \end{array}
\end{equation}
one may verify the commutation relations
\begin{equation}
    \begin{array}{lclcl}
         \big[T_0\Pi_0\Theta\big](\xi)&=&\big[\Pi_0T_{2,\frac{1}{2}}\Theta\big](\xi),\qquad 

         \big[T_0\Pi_1\Theta\big](\xi)&=&-\big[\Pi_1T_{2,\frac{1}{2}}\Theta\big](\xi) .
    \end{array}
\end{equation}
In particular, if $\Theta$ is in the kernel of $\overline{\K}_{2,\frac{1}{2}}+\lambda$, the functions
\begin{equation}\label{fdfhn:eq:qis2kernuitrekenen}
    \begin{array}{lclcl}
         X_0(\xi)&=&[\Pi_0\Theta](\xi),\qquad X_1(\xi)&=&e^{-\pi i\xi}[\Pi_1\Theta](\xi)
    \end{array}
\end{equation}
are eigenfunctions of the operator $L_0$ with eigenvalues $-\lambda$ and $-\lambda-\overline{c}_0\pi i$ respectively. Since $-\lambda$ and $-\lambda-\overline{c}_0\pi i$ cannot both be eigenvalues of $L_0$ at the same time in view of (\asref{fdfhn:spectralstability}{\text{HW}}), this means that at least one of the functions $X_0$ or $X_1$ is identically $0$.\\

Without loss, we assume that $X_0=0$. In this case, the function $\Theta$ can explicitly be identified as
\begin{equation}
    \begin{array}{lclcl}
        \Theta(0,\xi)&=&\frac{1}{2}e^{\pi i\xi}X_1(\xi),\qquad  \Theta(\frac{1}{2},\xi)&=&-\frac{1}{2}e^{\pi i\xi}X_1(\xi).
    \end{array}
\end{equation} 
As such, the eigenfunctions of $\overline{\K}_{q,\vartheta}$ can be expressed in terms of those of $L_0$, thus providing an upper bound on the dimension of the corresponding eigenspace. 

\begin{lemma}\label{fdfhn:kernuitrekenen} Consider the setting of Proposition \ref{fdfhn:fredholmeigenschappenkh}. Then for any $\lambda\in\C$ with $\Re\, \lambda\geq-\tilde{\lambda}$ and $\lambda\notin\overline{c}_02\pi i q^{-1}\Z$, we have the identity
\begin{equation}\label{fdfhn:tweedebeperkingkern}
\begin{array}{lcl}
\ker(\overline{\K}_{q,\vartheta}+\lambda)&= &\{0\}.
\end{array}
\end{equation}
In addition, we have the identity
\begin{equation}\label{fdfhn:eq:kerKqtheta}
\begin{array}{lcl}
\ker(\overline{\K}_{q,\vartheta})&=&\span\{\pi_{\perp}\Phi_0^+\}.
\end{array}
\end{equation}\end{lemma}
\textit{Proof.} Fix $\lambda\in\C$ with $\Re\, \lambda\geq-\tilde{\lambda}$. Suppose that $\Theta$ is in the kernel of the operator $\overline{\K}_{q,\vartheta}+\lambda$. For $n\in\{0,...,q-1\}$ we set
\begin{equation}\label{fdfhn:eq:harmonicproj}
    \begin{array}{lcl}
         [\Pi_n\Theta](\xi)&=&\sum\limits_{n'=0}^{q-1}\zeta_q^{n\cdot n'}\Theta\big(n'\vartheta ,\xi\big),
    \end{array}
\end{equation}
together with
\begin{equation}
\begin{array}{lclcl}
X_n(\xi)&=&e^{-\frac{2\pi in}{q}\xi}[\Pi_n\Theta](\xi)
&=&\zeta_q^{-n\xi}[\Pi_n\Theta](\xi),
\end{array}
\end{equation}
with $\zeta_q=\mathrm{exp}[2\pi i/q]$ the $q$-th root of unity. Recalling that $\mathrm{gcd}(\vartheta q,q)=1$, it follows that this sum contains each of the functions $\Theta\big(0,\xi\big),...,\Theta\big((q-1)q^{-1},\xi\big)$ exactly once. Recalling the definitions of the operators $T_0$ and $T_{q,\vartheta}$ from (\ref{fdfhn:definitieT0}) and (\ref{fdfhn:def:Ttheta}), we can compute
\begin{equation}\label{fdfhn:deltahoverzetten}
\begin{array}{lcl}
[T_0 \Pi_n\Theta](\xi)&=&[\Pi_n\Theta](\xi+1)\\[0.2cm]
&=&\sum\limits_{n'=0}^{q-1}\zeta_q^{nn'}\Theta\big(n'\vartheta ,\xi+1\big)\\[0.2cm]
&=&\sum\limits_{n'=0}^{q-1}\zeta_q^{nn'}(T_{q,\vartheta}\Theta)\big((n'-1)\vartheta ,\xi\big)\\[0.2cm]
&=&\zeta_q^n\sum\limits_{n'=0}^{q-1}\zeta_q^{n(n'-1)}(T_{q,\vartheta}\Theta)\big((n'-1)\vartheta ,\xi\big)\\[0.2cm]
&=&\zeta_q^n[\Pi_nT_{q,\vartheta}\Theta](\xi),
\end{array}
\end{equation}
which implies
\begin{equation}\label{fdfhn:deltahoverzetten2}
\begin{array}{lcl}
T_0 X_n(\xi)&=&\zeta_q^{-n(\xi+1)}[T_0\Pi_n\Theta](\xi+1)\\[0.2cm]
&=&\zeta_q^{-n(\xi+1)}\zeta_q^n[\Pi_nT_{q,\vartheta}\Theta](\xi)\\[0.2cm]
&=&\zeta_q^{-n\xi}[\Pi_nT_{q,\vartheta}\Theta](\xi).
\end{array}
\end{equation}
This allows us to obtain the identity 
\begin{equation}\label{fdfhn:Lhoverzetten}
\begin{array}{lcl}
(L_0+\lambda) X_n(\xi)&=&\overline{c}_0 X_n'(\xi)-\Delta_0 X_n(\xi)-D_U\mathcal{G}\big(\overline{U}_0(\xi);\overline{r}\big)X_n(\xi)+\lambda X_n(\xi)\\[0.2cm]
&=&\overline{c}_0 \zeta_q^{-n\xi}[\Pi_n\Theta]'(\xi)-\overline{c}_0\frac{2\pi in}{q}X_n(\xi)-\zeta_q^{-n\xi}[\Pi_n\Delta_{q,\vartheta}\Theta](\xi)\\[0.4cm]
&&\qquad -\zeta_q^{-n\xi}D_U\mathcal{G}\big(\overline{U}_0(\xi);\overline{r}\big)[\Pi_n\Theta](\xi)+\zeta_q^{-n\xi}\lambda[\Pi_n\Theta](\xi)\\[0.2cm]
&=&\zeta_q^{-n\xi} \Big[\Pi_n\big(\overline{\K}_{q,\vartheta}+\lambda\big)\Theta\Big](\xi)-\overline{c}_0\frac{2\pi in }{q}X_n(\xi)\\[0.3cm]
&=&-\overline{c}_0\frac{2\pi in }{q} X_n(\xi).
\end{array}
\end{equation}

Suppose first that $\lambda\notin 2\overline{c}_0\pi i q^{-1}\Z$. Then it follows from (\asref{fdfhn:spectralstability}{\text{HW}}) that $-2\overline{c}_0\pi inq^{-1}-\lambda$ is no eigenvalue of $L_0$ for all $0\leq n\leq q-1$. In particular, we must have $X_n=0$ for all $0\leq n\leq q-1$. This means that the functions $\Pi_n\Theta$ for $0\leq n\leq q-1$ are also identically $0$. Since the $q\times q$ Vandermonde matrix $Z$ given by $Z_{n,n'}=\zeta_q^{n\cdot n'}$ is invertible, we obtain $\Theta(n\vartheta ,\cdot)=0$ for all $0\leq n\leq q-1$ from which (\ref{fdfhn:tweedebeperkingkern}) follows.\\

Turning to the case $\lambda=0$, we see that $-2\overline{c}_0\pi inq^{-1}-\lambda=-2\overline{c}_0\pi inq^{-1}$ can only be an eigenvalue of $L_0$ when $nq^{-1}\in\Z$ on account of (\asref{fdfhn:spectralstability}{\text{HW}}). Since $nq^{-1}\notin \Z$ for $1\leq n\leq q-1$, we have $X_n=0$ for those values of $n$. In addition, we have $X_0=\mu \Phi_0^+$ for some $\mu\in\C$. Recalling the invertible matrix $Z$ given by $Z_{n,n'}=\zeta_q^{n\cdot n'}$, we obtain the identity
\begin{equation}
    \begin{array}{lcl}
         \Big(\Theta\big(0,\cdot\big),\Theta\big(\vartheta ,\cdot\big),....,\Theta\big((q-1)\vartheta ,\cdot\big)\Big)^T&=&Z^{-1}\big(\mu \Phi_0^+,0,...,0\big)^T.
    \end{array}
\end{equation}
In particular, the kernel $\ker(\overline{\K}_{q,\vartheta})$ is one-dimensional. Since $L_0\Phi_0^+=0$ by (\asref{fdfhn:spectralstability}{\text{HW}}), it follows immediately that $\overline{\K}_{q,\vartheta}\pi_{\perp}\Phi_0^+=0$, which implies (\ref{fdfhn:eq:kerKqtheta}).\qed\\

We now shift our attention to the Fredholm properties of $\overline{\K}_{q,\vartheta}$, which we aim to extract from those of $L_0$ in a similar fashion. The results in \cite{Faye2014,MPA} show that it suffices to consider the limiting operators
\begin{equation}\begin{array}{lcl}\overline{\K}_{q,\vartheta,\pm\infty}\Theta&=&\overline{c}_0\partial_\xi \Theta-\Delta_{q,\vartheta}\Theta-D\mathcal{G}\big(P^{\pm};\overline{r}\big) \Theta,\\[0.2cm]
L_{\pm\infty} \Theta&=&\overline{c}_0\partial_\xi \Theta-\Delta_{0}\Theta-D\mathcal{G}\big(P^{\pm};\overline{r}\big) \Theta,\end{array}\end{equation}
which have constant coefficients. For $\lambda\in\C$ and $0\leq \rho\leq 1$ we introduce the notation
\begin{equation}
\begin{array}{lcl}
\overline{\K}_{q,\vartheta,\rho;\lambda}&=&\rho\overline{\K}_{q,\vartheta,-\infty}+(1-\rho)\overline{\K}_{q,\vartheta,\infty}+\lambda,\\[0.2cm]
L_{\rho;\lambda}&=&\rho L_{-\infty}+(1-\rho)L_{\infty}+\lambda.
\end{array}
\end{equation} 
We set out to show that for $\lambda$ in a suitable right half-plane and $0\leq\rho\leq 1$, the operators $\overline{\K}_{q,\vartheta,\rho;\lambda}$ and $L_{\rho;\lambda}$ are hyperbolic in the sense of \cite{Faye2014,MPA}. In particular, we write
\begin{equation}\label{fdfhn:definitiedeltalh}\begin{array}{lclcl}\Delta_{q, \vartheta, \rho;\lambda}(z)&=&\Big[\overline{\K}_{q,\vartheta,\rho;\lambda}e^{z\xi}\Big](0),\qquad 
\Delta_{\rho;\lambda}(z)&=&\Big[L_{\rho;\lambda}e^{z\xi}\Big](0)
\end{array}\end{equation}
and establish that $\det\big(\Delta_{q, \vartheta, \rho;\lambda}(iy)\big)\neq 0$ for all $y\in\R$ by first showing that $\det\big(\Delta_{\rho;\lambda}(iy)\big)\neq 0$. We can subsequently use the spectral flow principle to compute the Fredholm index of $\overline{\K}_{q,\vartheta}+\lambda$.\\

We start by considering the characteristic function $\Delta_{\rho;\lambda}$ from (\ref{fdfhn:definitiedeltalh}). For notational convenience we set 
\begin{equation}
\begin{array}{lcl}
D\mathcal{G}_\rho&=&\rho D\mathcal{G}\big(P^-;\overline{r}\big)+(1-\rho)D\mathcal{G}\big(P^+;\overline{r}\big)
\end{array}
\end{equation}
for $0\leq \rho\leq 1$ and use the definition (\ref{fdfhn:conditionAz}) to write 
\begin{equation}
\begin{array}{lcl}
\Delta_{\rho;\lambda}(iy)&=&\overline{c}_0 iy -\tau\sum\limits_{m>0}\alpha_m\Big[e^{m iy}+e^{-miy}-2\Big]-D\mathcal{G}_\rho+\lambda\\[0.2cm]
&=&\overline{c}_0 iy+\tau\sum\limits_{m>0}\alpha_m\Big[2-2\cos(my)\Big]-D\mathcal{G}_\rho+\lambda\\[0.2cm]
&=&\overline{c}_0 iy+2\tau A(y)-D\mathcal{G}_\rho+\lambda.
\end{array}
\end{equation}
For any $V=(v_1,...,v_d)\in\C^d$ we may exploit the inequality (\ref{fdfhn:conditionAz}) to obtain
\begin{equation}
\begin{array}{lclcl}
\tau V^\dagger A(y)V&=&2\tau\sum\limits_{j=1}^d |v_j|^2A_j(y)

&\geq &0.\end{array}
\end{equation}
Here we introduced $\dagger$ for the conjugate transpose.\\

In order to prove that $L_{\pm\infty}+\lambda$ is hyperbolic, we need to distinguish between the setting where the triplet $(\mathcal{G},P^-,P^+)$ satisfies (\asref{fdfhn:aannamesconstantena}{\text{HS}3_{\overline{r}}}) and where it satisfies (\asref{fdfhn:aannamesconstantenb}{\text{HS}3_{\overline{r}}}). A similar computation was performed in \cite[Lem. 4.6]{SCHPerExtensions}. 

\begin{lemma}\label{fdfhn:hyperbolicalpha}
Assume that (\asref{fdfhn:aannames}{\text{HS}1}) and (\asref{fdfhn:aannames2}{\text{HS}2}) are satisfied and pick $\overline{r}$ in such a way that (\asref{fdfhn:existentie}{\text{HW}}) and (\asref{fdfhn:spectralstability}{\text{HW}}) are satisfied. Assume that the triplet $(\mathcal{G},P^-,P^+)$ satisfies (\asref{fdfhn:aannamesconstantena}{\text{HS}3_{\overline{r}}}). Pick $\lambda\in\C$ with $\Re\, \lambda>-\tilde{\lambda}$ and $0\leq \rho\leq 1$. Then we have $\det\big(\Delta_{\rho;\lambda}(iy)\big)\neq 0$ for all $y\in\R$. \end{lemma}
\textit{Proof.} For fixed $y\in\R$ we introduce the matrix
\begin{equation}
\begin{array}{lcl}
X&=&\frac{1}{2}\big[\Delta_{\rho;\lambda}(iy)+\Delta_{\rho;\lambda}(iy)^\dagger\big]\\[0.2cm]
&=&\tau A(y)-D\mathcal{G}_\rho-D\mathcal{G}_\rho^T+\Re\, \lambda.
\end{array}
\end{equation} 
By decreasing $\tilde{\lambda}$ if necessary, we can assume that $-D\mathcal{G}_\rho-D\mathcal{G}_\rho^T+\Re\, \lambda$ is positive definite. It follows that $X$ is the sum of a positive semi-definite matrix and a positive definite matrix and as such, it is positive definite itself. As a consequence, $\Delta_{\rho;\lambda}$ is positive definite as well and hence we obtain $\det\big(\Delta_{\rho;\lambda}(iy)\big)\neq 0$.\qed\\

\begin{lemma}\label{fdfhn:hyperbolicbeta}
Assume that (\asref{fdfhn:aannames}{\text{HS}1}) and (\asref{fdfhn:aannames2}{\text{HS}2}) are satisfied and pick $\overline{r}$ in such a way that (\asref{fdfhn:existentie}{\text{HW}}) and (\asref{fdfhn:spectralstability}{\text{HW}}) are satisfied. Assume that the triplet $(\mathcal{G},P^-,P^+)$ satisfies (\asref{fdfhn:aannamesconstantenb}{\text{HS}3_{\overline{r}}}). Pick $\lambda\in\C$ with $\Re\, \lambda>-\tilde{\lambda}$ and $0\leq \rho\leq 1$. Then we have $\det\big(\Delta_{\rho;\lambda}(iy)\big)\neq 0$ for all $y\in\R$.\end{lemma}
\textit{Proof.} We recall the proportionality constant $\Gamma>0$ from (\asref{fdfhn:aannamesconstantenb}{\text{HS}3_{\overline{r}}}). In particular, upon writing
\begin{equation}
\begin{array}{lcl}
D\mathcal{G}_\rho&=&\left(\begin{array}{lcl}
D\mathcal{G}_{\rho}^{[1,1]}& D\mathcal{G}_{\rho}^{[1,2]}\\
D\mathcal{G}_{\rho}^{[2,1]}& D\mathcal{G}_{\rho}^{[2,2]}
\end{array}\right),
\end{array}
\end{equation} 
we have $D\mathcal{G}_{\rho}^{[1,2]}=-\Gamma (D\mathcal{G}_{\rho}^{[2,1]})^T$. Suppose that $\Delta_{\rho;\lambda}(iy)V=0$ for some $V\in\C^d$. Write $V=(u,w)$ where $u$ contains the first $d_{\mathrm{diff}}$ components of $V$. Then we can compute

\begin{equation}\label{fdfhn:eq:VYVis0}
\begin{array}{lcl}
0&=&\Re\,  V^\dagger \Delta_{\rho;\lambda}(iy)V\\[0.2cm]
&=&\Re\, \Big[-\tau V^\dagger A(y)V- V^\dagger D\mathcal{G}_\rho V+\lambda |V|^2\Big]\\[0.2cm]
&=&\Re\, \Big[-\tau V^\dagger A(y)V-  u^\dagger D\mathcal{G}_{\rho}^{[1,1]}u- u^\dagger D\mathcal{G}_{\rho}^{[1,2]}w -  w^\dagger D\mathcal{G}_{\rho}^{[2,1]} u- w^\dagger D\mathcal{G}_{\rho}^{[2,2]}w+\lambda |u|^2+  \lambda |w|^2\Big].
\end{array}
\end{equation}
The second component of the equation $\Delta_{\rho;\lambda}(iy)V=0$ is equivalent to
\begin{equation}
    \begin{array}{lcl}
          D\mathcal{G}_{\rho}^{[2,1]}u&=&  - D\mathcal{G}_{\rho}^{[2,2]}w+\lambda w.
    \end{array}
\end{equation}
As such, we can rewrite the cross-terms in (\ref{fdfhn:eq:VYVis0}) to obtain
\begin{equation}
\begin{array}{lcl}
\Re\, \big[- u^\dagger  D\mathcal{G}_{\rho}^{[1,2]}w-  w^\dagger D\mathcal{G}_{\rho}^{[2,1]} u\big]
&=&\Re\, (1-\Gamma)\Big[-  w^\dagger D\mathcal{G}_{\rho}^{[2,1]}u\Big]\\[0.2cm]
&=&\Re\, (\Gamma-1)\Big[-  w^\dagger D\mathcal{G}_{\rho}^{[2,2]}w+\lambda |w|^2\Big].
\end{array}
\end{equation}
As a consequence, (\ref{fdfhn:eq:VYVis0}) reduces to
\begin{equation}
\begin{array}{lcl}
0&=&\Re\, \Big[-\tau V^\dagger A(y)V-  u^\dagger D\mathcal{G}_{\rho}^{[1,1]}u+\lambda |u|^2 -\Gamma  w^\dagger D\mathcal{G}_{\rho}^{[2,2]}w+ \Gamma \lambda |w|^2 \Big].
\end{array}
\end{equation}
By decreasing $\tilde{\lambda}$ if necessary, we can assume that $-D\mathcal{G}_{\rho}^{[1,1]}+\Re\, \lambda$ and $-\Gamma D\mathcal{G}_{\rho}^{[2,2]}+\Gamma\Re\, \lambda$ are positive definite. Therefore, we must have $V=0$, from which it follows that $\det\big(\Delta_{\rho;\lambda}(iy)\big)\neq 0$.\qed\\

\begin{lemma}\label{fdfhn:hyperbolic2} Consider the setting of Proposition \ref{fdfhn:fredholmeigenschappenkh}. Pick $\lambda\in\C$ with $\Re\, \lambda>-\tilde{\lambda}$ and $0\leq \rho\leq 1$. Then we have $\det\big(\Delta_{q,\vartheta,\rho;\lambda}(iy)\big)\neq 0$ for all $y\in\R$.\end{lemma}
\textit{Proof.} Suppose there exists $V\in  {\ell}_{q,\perp;\infty}^2$ and $y\in\R$ for which
\begin{equation}\begin{array}{lcl}\Delta_{q,\vartheta,\rho;\lambda}(iy)V&=&0.\end{array}\end{equation}
We then write
\begin{equation}
\begin{array}{lcl}
W\big(\frac{n}{q},\xi\big)&=&e^{iy\xi}V\big(\frac{n}{q}\big)
\end{array}
\end{equation} 
for $0\leq n\leq q-1$. The definition of the characteristic function yields
\begin{equation}
\begin{array}{lcl}
\overline{\K}_{q,\vartheta,\rho;\lambda}W&=&e^{iy\xi}\left[\overline{\K}_{q,\vartheta,\rho;\lambda}e^{iy\xi}V\right](0)\\[0.6cm]
&=&e^{iy\xi}\Delta_{q,\vartheta,\rho;\lambda}(iy)V\\[0.2cm]
&=&0.
\end{array}
\end{equation}
Recalling the projections (\ref{fdfhn:eq:harmonicproj}),
we write
\begin{equation}
\begin{array}{lcl}
X_n(\xi)&=&e^{-\frac{2\pi i n}{q}\xi}[\Pi_nW](\xi)
\end{array}
\end{equation} 
and use a computation similar to (\ref{fdfhn:Lhoverzetten}) to find
\begin{equation}
\begin{array}{lcl}
L_{\rho;\lambda}X_n(\xi)&=&e^{-\frac{2\pi i n}{q}\xi}\big[\Pi_n\overline{\K}_{q,\vartheta,\rho;\lambda} W\big](\xi)-\overline{c}_0\frac{2\pi i n}{q}X_n(\xi)\\[0.2cm]
&=&-\overline{c}_0\frac{2\pi i n}{q} X_n(\xi).
\end{array}
\end{equation}
On account of Lemmas \ref{fdfhn:hyperbolicalpha}-\ref{fdfhn:hyperbolicbeta}, it follows from the spectral flow theorem \cite[Thm. 1.6]{Faye2014} and \cite[Thm. 1.7]{Faye2014} that $L_{\rho;\lambda-\overline{c}_02\pi i nq^{-1} }$ is hyperbolic. Applying \cite[Lem. 6.3]{HJHFHNINFRANGE}, which is a generalization of \cite[Thm. 4.1]{MPA}, yields that $L_{\rho;\lambda-\overline{c}_02\pi i n\vartheta } $ is invertible as a map from $W^{1,\infty}(\R;\R^d)$ to $L^\infty(\R;\R^d)$. Therefore, we must have $X_n=0$ for all $0\leq n\leq q-1$. This implies that $W(\frac{n}{q},\xi)=0$ for all $0\leq n\leq q-1$ and thus that $V=0$, which yields the desired result.\qed\\

\textit{Proof of Proposition \ref{fdfhn:fredholmeigenschappenkh}.} These results, except the bound (\ref{fdfhn:batesafschattingKh}), follow from combining Lemma \ref{fdfhn:kernuitrekenen}, Lemma \ref{fdfhn:hyperbolic2} and the spectral flow theorem \cite[Thm. 1.6-1.7]{Faye2014}. The bound (\ref{fdfhn:batesafschattingKh}) can be obtained by following the proof of \cite[Lem. 5]{BatesInfRange}.\qed\\

\section{Linear theory for $\Delta t\rightarrow 0$}\label{fdfhn:section:lineartheory}

In this section, we apply the spectral convergence method to lift the Fredholm properties of the semi-discrete system to the fully discrete system in the small timestep limit $\Delta t\rightarrow 0$. In particular, we establish the main result below, which gives a quasi-inverse for the operators $L_{k,M}$. This turns out to be the key ingredient in the construction of the discrete waves, which can subsequently be proved by means of a standard fixed point argument.
\begin{proposition}[{cf. \cite[Prop. 3.2]{HJHBDF}}]\label{fdfhn:Prop32BDF} Assume that (\asref{fdfhn:aannames}{\text{HS}1}) and (\asref{fdfhn:aannames2}{\text{HS}2}) are satisfied and pick $\overline{r}$ in such a way that (\asref{fdfhn:aannamesconstanten2}{\text{HS}}), (\asref{fdfhn:existentie}{\text{HW}}) and (\asref{fdfhn:spectralstability}{\text{HW}}) are satisfied.  Fix a pair of integers $1\leq k\leq 6$ and $q\geq 1$, together with a sufficiently small $\eta>0$ and sufficiently large constants $M_*\in\M_{q}$ and $C>0$. Then for each $M\in \M_{q}$ with $M\geq M_*$ there exist linear maps
\begin{equation}\begin{array}{lclcl}\gamma_{k,M}^*: {\Y}_M&\rightarrow&\R,\qquad \V_{k,M}^*: {\Y}_M&\rightarrow & {\Y}_{k,M}^1,\end{array}\end{equation}
so that for all $\Psi\in  {\Y}_M$ the pair
\begin{equation}\begin{array}{lcl}(\gamma,V)&=&(\gamma_{k,M}^*\Psi,\V_{k,M}^*\Psi)\end{array}\end{equation}
is the unique solution to the problem
\begin{equation}\label{fdfhn:bdf1}\begin{array}{lcl}L_{k,M}V&=&\Psi+\gamma {\pi}_{\Y_M} {\D}_{k,M}\overline{U}_0\end{array}\end{equation}
that satisfies the normalisation condition
\begin{equation}\label{fdfhn:bdf2}\begin{array}{lcl}\ip{ {\pi}_{\Y_M}\Phi_0^-,V}_{ {\Y}_M}&=&0.\end{array}\end{equation}
In addition, for all $\Psi\in  {\Y}_M$ we have the bound
\begin{equation}\label{fdfhn:eq315BDF}\begin{array}{lcl}|\gamma_{k,M}^*\Psi|+\nrm{\V_{k,M}^*\Psi}_{ {\Y}_{k,M}^1}&\leq &C\nrm{\Psi}_{ {\Y}_M}.\end{array}\end{equation}
\end{proposition}
\textit{Proof of Theorem \ref{fdfhn:maintheorem}.} On account of Proposition \ref{fdfhn:Prop32BDF}, the procedure outlined in \cite[\S 4.1]{HJHBDF} can be followed to arrive at the desired result.\qed\\

In order to facilitate the reading, we first outline our strategy and formulate two intermediate results in \S \ref{fdfhn:subsec:strategy}. This strategy heavily follows the program in \cite{HJHBDF}, allowing us to simply refer to these results in many cases. However, due to the lack of a comparison principle and the many cross-terms we need to control, there are several key points in the analysis that need a fully new approach, which we develop in \S \ref{fdfhn:sectiontransfer}. In addition, the infinite-range setting forces us to obtain an extra order of regularity on the operator $(L_0+\delta)^{-1}$, which we achieve in \S \ref{fdfhn:sectionexpdecay}.

\subsection{Strategy}\label{fdfhn:subsec:strategy}

%\begin{proposition}[{cf. \cite[Prop. 3.3]{HJHBDF}}]\label{fdfhn:Prop33BDF} Assume that (\asref{fdfhn:aannames}{\text{HS}1}) and (\asref{fdfhn:aannames2}{\text{HS}2}) are satisfied and pick $\overline{r}$ in such a way that (\asref{fdfhn:aannamesconstanten2}{\text{HS}}), (\asref{fdfhn:existentie}{\text{HW}}) and (\asref{fdfhn:spectralstability}{\text{HW}}) are satisfied.  Fix a pair of integers $1\leq k\leq 6$ and $q\geq 1$. Then there exists $C_0>0$ together with a map $M_0:(0,\delta_0)\rightarrow [1,\infty)$ so that the following holds true. For any $0<\delta<\delta_0$ and any $M\in\M_{q}$ for which $M\geq M_0(\delta)$, the operator $L_{k,M}+\delta$ is invertible as a map from $ {\Y}_{k,M}^1$ to $ {\Y}_M$, with the bound
%\begin{equation}\label{fdfhn:Eq320BDF}
%\begin{array}{lcl}
%\nrm{(L_{k,M}+\delta)^{-1}\Psi}_{ {\Y}_{k,M}^1}&\leq & C_0\Big[ \nrm{\Psi}_{ {\Y}_M}+\delta^{-1}|\ip{ {\pi}_{\Y_M}\Phi_0^-,\Psi}_{ {\Y}_M}|\Big].
%\end{array}
%\end{equation}\end{proposition}

Recalling the spaces $\H_{M}$ and $\H_{k,M}^1$ from (\ref{fdfhn:eq:defHM}) and (\ref{fdfhn:eq:defHkM1}), we introduce the quantities
\begin{equation}\label{fdfhn:definitiemathcalE}
\begin{array}{lcl}
\mathcal{E}_{k,M}(\delta)&=&\inf_{\nrm{\Phi}_{ {\H}_{k,M}^1}=1}\Big[\nrm{\mathcal{K}_{k,M}\Phi+\delta \Phi}_{ {\H}_M}+\delta^{-1}\Big|\ip{\pi_{\H_M} \Phi_0^-,\mathcal{K}_{k,M}\Phi+\delta \Phi}_{ {\H}_M}\Big|\Big],\\[0.2cm]
\mathcal{E}_{k,M}^*(\delta)&=&\inf_{\nrm{\Phi}_{ {\H}_{k,M}^1}=1}\Big[\nrm{\mathcal{K}_{k,M}^*\Phi+\delta \Phi}_{ {\H}_M}+\delta^{-1}\Big|\ip{\pi_{\H_M} \Phi_0^+,\mathcal{K}^*_{k,M}\Phi+\delta \Phi}_{ {\H}_M}\Big|\Big],
\end{array}
\end{equation}
together with
\begin{equation}\label{fdfhn:definitiekappadelta}
\begin{array}{lcl}
\kappa(\delta)&=&\liminf_{M\rightarrow\infty,M\in\M_{q}}\mathcal{E}_{k,M}(\delta),\\[0.2cm]
\kappa^*(\delta)&=&\liminf_{M\rightarrow\infty,M\in\M_{q}}\mathcal{E}^*_{k,M}(\delta)
\end{array}
\end{equation}
for $\delta\in(0,\delta_0)$.\\[0.2cm]
The key step towards proving Proposition \ref{fdfhn:Prop32BDF} is the establishment of lower bounds for these quantities. This procedure is based on \cite[Lem. 6]{BatesInfRange}. Our strategy to prove it is essentially the same, but some major modifications are needed to incorporate the difficulties arising from the discrete derivatives.
\begin{proposition}[{cf. \cite[Prop. 3.7]{HJHBDF}}]\label{fdfhn:Prop37BDF} Assume that (\asref{fdfhn:aannames}{\text{HS}1}) and (\asref{fdfhn:aannames2}{\text{HS}2}) are satisfied and pick $\overline{r}$ in such a way that (\asref{fdfhn:aannamesconstanten2}{\text{HS}}), (\asref{fdfhn:existentie}{\text{HW}}) and (\asref{fdfhn:spectralstability}{\text{HW}}) are satisfied.  Fix a pair of integers $1\leq k\leq 6$ and $q\geq 1$. Then there exists $\kappa>0$ such that for all $0<\delta<\delta_0$ we have
\begin{equation}\label{fdfhn:Eq359BDF}
\begin{array}{lclcl}
\kappa(\delta)&\geq & \kappa,\qquad \kappa^*(\delta)&\geq & \kappa.
\end{array}
\end{equation}\end{proposition}

We are now ready to start our interpolation procedure. For any $\xi\in\R$, we pick two quantities $\xi_M^\pm(\xi)\in M^{-1}\Z$ in such a way that
\begin{equation}
\begin{array}{lclclcl}
\xi_M^-(\xi)&\leq & \xi & < & \xi_M^+(\xi),\qquad \xi_M^+(\xi)-\xi_M^-(\xi)&=&M^{-1}.
\end{array}
\end{equation}
Using these quantities, we can define two interpolation operators
\begin{equation}\label{fdfhn:eq:interpolationoperator}
\begin{array}{lclcl}
\I_M^0& : &\H_M & \rightarrow & L^2(\R,\ell_{q,\perp;\infty}^2),\\[0.2cm]
\I_{k,M}^1 & : & \H_{k,M}^1 &\rightarrow & H^1(\R,\ell_{q,\perp;\infty}^2),
\end{array}
\end{equation}
that act as
\begin{equation}
\begin{array}{lcl}
[\I_M^0\phi](\zeta,\xi)&=& \phi\Big(\zeta,\xi_M^-(\xi)\Big),\\[0.2cm]
[\I_{k,M}^1\phi](\zeta,\xi)&=& M\Big[\big(\xi-\xi_M^-(\xi)\big)\phi\Big(\zeta,\xi_M^+(\xi)\Big)+(\xi_M^+(\xi)-\xi)\phi\Big(\zeta,\xi_M^-(\xi)\Big)\Big],
\end{array}
\end{equation}
for all $\zeta\in q^{-1}\Z_q$ and all $\xi\in\R$. These operators can be seen as interpolations of order zero and one respectively, both acting only on the second coordinate of $\phi$. We refer to \cite[Lem. 3.10-3.12]{HJHBDF} for some useful estimates involving these interpolations.\\

With these preparations in hand, we start the proof of Proposition \ref{fdfhn:Prop37BDF} using the methods described in the proof of \cite[Lem. 6]{BatesInfRange}. We focus on the quantity $\kappa(\delta)$ defined in (\ref{fdfhn:definitiekappadelta}), noting that $\kappa^*(\delta)$ can be treated in a similar fashion. In particular, we find a lower bound for $\kappa(\delta)$ by constructing sequences that minimize this quantity. At this point it becomes clear why we work on the spaces $ {H}^1(\R,\ell_{q,\perp}^2)$ and $ {L}^2(\R,\ell_{q,\perp}^2)$, as we exploit the fact that bounded closed subsets of these spaces are weakly compact.

\begin{lemma}[{cf. \cite[Lem. 3.16-3.17]{HJHBDF}}]\label{fdfhn:Lemma316BDF} Assume that (\asref{fdfhn:aannames}{\text{HS}1}) and (\asref{fdfhn:aannames2}{\text{HS}2}) are satisfied and pick $\overline{r}$ in such a way that (\asref{fdfhn:aannamesconstanten2}{\text{HS}}), (\asref{fdfhn:existentie}{\text{HW}}) and (\asref{fdfhn:spectralstability}{\text{HW}}) are satisfied.  Fix a pair of integers $1\leq k\leq 6$ and $q\geq 1$, as well as $0<\delta<\delta_0$. Then there exist two functions
\begin{equation}
\begin{array}{lclcl}
\Phi_*&\in & {H}^1(\R,\ell_{q,\perp;\infty}^2),\qquad
\Psi_*&\in & {L}^2(\R,\ell_{q,\perp;\infty}^2),
\end{array}
\end{equation} 
together with three sequences
\begin{equation}\label{fdfhn:threesequences}
\begin{array}{lclclcl}
\{M_j\}_{j\in\N}&\subset& \M_{q},\enskip \{\Phi_j\}_{j\in\N}&\subset &  {\H}_{k,M_j}^1,\enskip \{\Psi_j\}_{j\in\N}&\subset & {\H}_{M_j}
\end{array}
\end{equation}
and two constants $\vartheta\in q^{-1}\Z_q\setminus\{0\}$ and $K_1>0$ that satisfy the following properties.
\begin{enumerate}[label=(\roman*)]
\item We have $\lim_{j\rightarrow\infty}M_j=\infty$ and $\nrm{\Phi_j}_{ {\H}_{k,M_j}^1}=1$ for all $j\in\N$.
\item The identity
\begin{equation}
\begin{array}{lcl}
\Psi_j&=&\mathcal{K}_{k,M_j}\Phi_j+\delta \Phi_j
\end{array}
\end{equation}
holds for all $j\in\N$.
\item Recalling the constant $\kappa(\delta)$ defined in (\ref{fdfhn:definitiekappadelta}), we have the limit
\begin{equation}\label{fdfhn:limietnaarkappadelta}
\begin{array}{lcl}
\kappa(\delta)&=&\lim_{j\rightarrow\infty}\Big[\nrm{\mathcal{K}_{k,M}\Phi_j+\delta \Phi_j}_{ {\H}_{M_j}}+\delta^{-1}\Big|\ip{\pi_{\H_{M_j}}\Phi_0^-,\mathcal{K}_{k,M_j}\Phi_j+\delta \Phi_j}_{ {\H}_{M_j}}\Big|\Big].
\end{array}
\end{equation}
\item As $j\rightarrow\infty$, we have the weak convergences
\begin{equation}
\begin{array}{lcl}
\I_{k,M_j}^1\Phi_j&\rightharpoonup&\Phi_*\in  {H}^1(\R,\ell_{q,\perp}^2),\\[0.2cm]
\I_{M_j}^0\Psi_j&\rightharpoonup&\Psi_*\in  {L}^2(\R,\ell_{q,\perp}^2).
\end{array}
\end{equation}
\item For any compact interval $\I\subset \R$, we have the strong convergences
\begin{equation}\label{fdfhn:eq:strongconve}
\begin{array}{lcl}
(\I_{k,M_j}^1,\I_{k,M_j}^1)\Phi_j&\rightarrow&\Phi_*\in  {L}^2(\I,\ell_{q,\perp}^2),\\[0.2cm]
(\I_{M_j}^0,\I_{M_j}^0)\Psi_j&\rightarrow&\Psi_*\in  {L}^2(\I,\ell_{q,\perp}^2)
\end{array}
\end{equation}
as $j\rightarrow\infty$.
\item The function $\Phi_*$ is a weak solution to $(\overline{\K}_{q,\vartheta}+\delta)\Phi_*=\Psi_*$ and we have the bound
\begin{equation}\label{fdfhn:eq:specconv1big}
\begin{array}{lcl}
\nrm{\Phi_*}_{ {H}^1(\R,\ell_{q,\perp;\infty}^2)}&\leq & K_1\kappa(\delta).
\end{array}
\end{equation}
\end{enumerate}\end{lemma}
\textit{Proof.} In view of Proposition \ref{fdfhn:fredholmeigenschappenkh} and Lemma \ref{fdfhn:lemmaconvergencetestfunction}, we can follow the proof of \cite[Lem. 3.16-3.17]{HJHBDF} almost verbatim.\qed\\

In order to prove Proposition \ref{fdfhn:Prop37BDF}, we need to establish a lower bound on the norm $\nrm{\Phi_*}_{ {H}^1(\R,\ell_{q,\perp;\infty}^2)}$ on account of (\ref{fdfhn:eq:specconv1big}). In Proposition \ref{fdfhn:Lemma318BDF} we follow the approach of \cite[Lem. 3.18]{HJHBDF} in order to obtain this lower bound. Here we have to deal with both the cross-terms arising from the system setting as well as the infinite-range interactions. 

\begin{proposition}[{see \S \ref{fdfhn:sectiontransfer}}]\label{fdfhn:Lemma318BDF} Consider the setting of Lemma \ref{fdfhn:Lemma316BDF}. Then there exist constants $K_2>1$ and $K_3>1$ so that for any $0<\delta<\delta_0$, the function $\Phi_*$ satisfies the bound
\begin{equation}\label{fdfhn:eq:specconv2big}
\begin{array}{lcl}
\nrm{\Phi_*}_{H^1(\R,\ell_{q,\perp;\infty}^2)}^2&\geq & K_2-K_3\kappa(\delta)^2.
\end{array}
\end{equation}\end{proposition}

\textit{Proof of Proposition \ref{fdfhn:Prop37BDF}.} Combining the bounds (\ref{fdfhn:eq:specconv1big}) and (\ref{fdfhn:eq:specconv2big}) immediately yields
\begin{equation}
    \begin{array}{lcl}
         K_2-K_3\kappa(\delta)^2&\leq &K_1^2\kappa(\delta)^2.
    \end{array}
\end{equation}
Solving this quadratic inequality, we obtain
\begin{equation}
    \begin{array}{lclcl}
         \kappa(\delta)&\geq & \sqrt{\frac{K_2}{K_1^2+K_3}}&:=&\kappa.
    \end{array}
\end{equation}
The lower bound on $\kappa^*(\delta)$ follows in a similar fashion.\qed\\

In order to establish Proposition \ref{fdfhn:Prop32BDF}, we need more control on the operator $L_0$ than in \cite{HJHFHNINFRANGE}. In particular, due to the infinite-range interactions it is not immediately clear that this operator preserves the exponential decay properties of the function spaces (\ref{fdfhn:eq:bcetadef}).
\begin{proposition}\label{fdfhn:exponentieelvervalprop} Assume that (\asref{fdfhn:aannames}{\text{HS}1}) and (\asref{fdfhn:aannames2}{\text{HS}2}) are satisfied and pick $\overline{r}$ in such a way that (\asref{fdfhn:aannamesconstanten2}{\text{HS}}), (\asref{fdfhn:existentie}{\text{HW}}) and (\asref{fdfhn:spectralstability}{\text{HW}}) are satisfied. Fix a sufficiently small constant $\eta>0$. Then there exist constants $\delta_*>0$ and $K>0$, so that for each $0<\delta<\delta_*$ and each $G \in BC_{-\eta}^1(\R;\R^d)$ we have the bounds
\begin{equation}
\begin{array}{lcl}
\nrm{(L_0+\delta)^{-1}G}_{BC_{-\eta}(\R;\R^d)}&\leq & K \delta^{-1}\nrm{G}_{BC_{-\eta}(\R;\R^d)}
\\[0.2cm]
\nrm{[(L_0+\delta)^{-1}G]'}_{BC_{-\eta}(\R;\R^d)}&\leq & K \delta^{-1}\nrm{G}_{BC_{-\eta}(\R;\R^d)}
\\[0.2cm]
\nrm{[(L_0+\delta)^{-1}G]''}_{BC_{-\eta}(\R;\R^d)}&\leq & K \delta^{-1}\nrm{G}_{BC_{-\eta}^1(\R;\R^d)}.
\end{array}
\end{equation}
\end{proposition}

\textit{Proof of Proposition \ref{fdfhn:Prop32BDF}.} On account of Proposition \ref{fdfhn:exponentieelvervalprop}, we can follow the procedure developed in \cite[\S 3.3]{HJHBDF} to arrive at the desired result.\qed\\

\subsection{Spectral convergence}\label{fdfhn:sectiontransfer}

In this section we set out to prove Proposition \ref{fdfhn:Lemma318BDF} using the spectral convergence method. The main idea is to derive an upper bound for the discrete derivative $\D_{k,M_j}\Phi_j$, together with a lower bound for $\Phi_j$ restricted to a large---but finite---interval. This prevents the $\H^1_{k,M_j}$-norm of $\Phi_j$ from leaking away into oscillations or tail effects, providing the desired control on the limit (\ref{fdfhn:eq:strongconve}). All constants introduced in Lemmas \ref{fdfhn:Lemma318BDFdeel1}-\ref{fdfhn:Lemma318BDFdeel2b} and Proposition \ref{fdfhn:Lemma318BDF} are independent of $0<\delta<\delta_0$.

\begin{lemma}\label{fdfhn:Lemma318BDFdeel1} Consider the setting of Lemma \ref{fdfhn:Lemma316BDF}. Then there exists a constant $C_1>0$ so that the bound 
\begin{equation}\label{fdfhn:eq:specconvest}
\begin{array}{lcl}
2\nrm{\Psi_j}_{ {\H}_{M_j}}^2+2C_1\nrm{\Phi_j}_{ {\H}_{M_j}}^2&\geq &\overline{c}_0^2\nrm{ {\D}_{k,M_j}\Phi_j}_{ {\H}_{M_j}}^2
\end{array}
\end{equation}
holds for all $j\in\N$.\end{lemma}
\textit{Proof.} We will assume $\overline{c}_0>0$, noting that the case where $\overline{c}_0<0$ can be treated in a similar fashion. In view of the identity
\begin{equation}
\begin{array}{lcl}
\K_{k,M_j}\Phi_j+\delta \Phi_j&=&\Psi_j,
\end{array}
\end{equation}
we can compute
\begin{equation}
\begin{array}{lcl}
\ip{\Psi_j, {\D}_{k,M_j}\Phi_j}_{ {\H}_{M_j}}&=&\overline{c}_0\nrm{ {\D}_{k,M_j}\Phi_j}_{ {\H}_{M_j}}^2-\ip{\Delta_{M_j}\Phi_j, {\D}_{k,M_j}\Phi_j}_{ {\H}_{M_j}}\\[0.2cm]
&&\qquad -\ip{D\mathcal{G}\big(\pi_{\H_{M_j}}\overline{U}_0;\overline{r}\big)\Phi_j, {\D}_{k,M_j}\Phi_j}_{ {\H}_{M_j}}+\delta\ip{\Phi_j, {\D}_{k,M_j}\Phi_j}_{ {\H}_{M_j}}.
\end{array}
\end{equation}
Writing 
\begin{equation}
\begin{array}{lcl}
K&=&\nrm{D\mathcal{G}\big(\overline{U}_0;\overline{r}\big)}_\infty+4\tau\sum\limits_{m>0}|\alpha_m|
\end{array}
\end{equation}
and remembering that $0<\delta<\delta_0<1$, we may use the Cauchy-Schwarz inequality to obtain 
\begin{equation}
\begin{array}{lcl}
K\nrm{\Phi_j}_{ {\H}_{M_j}}\nrm{ {\D}_{k,M_j}\Phi_j}_{ {\H}_{M_j}}&\geq & \ip{\Delta_{M_j}\Phi_j, {\D}_{k,M_j}\Phi_j}_{ {\H}_{M_j}}+\ip{D\mathcal{G}\big(\pi_{\H_{M_j}}\overline{U}_0;\overline{r}\big)\Phi_j, {\D}_{k,M_j}\Phi_j}_{ {\H}_{M_j}}\\[0.2cm]
&&\qquad -\delta\ip{\Phi_j, {\D}_{k,M_j}\Phi_j}_{ {\H}_{M_j}}\\[0.2cm]
&=&\overline{c}_0\nrm{ {\D}_{k,M_j}\Phi_j}_{ {\H}_{M_j}}^2-\ip{\Psi_j, {\D}_{k,M_j}\Phi_j}_{ {\H}_{M_j}}\\[0.2cm]
&\geq &\overline{c}_0\nrm{ {\D}_{k,M_j}\Phi_j}_{ {\H}_{M_j}}^2-\nrm{\Psi_j}_{ {\H}_{M_j}}\nrm{ {\D}_{k,M_j}\Phi_j}_{ {\H}_{M_j}}.
\end{array}
\end{equation}
This yields the bound
\begin{equation}
\begin{array}{lcl}
\nrm{\Psi_j}_{ {\H}_{M_j}}+K\nrm{\Phi_j}_{ {\H}_{M_j}}&\geq &\overline{c}_0\nrm{ {\D}_{k,M_j}\Phi_j}_{ {\H}_{M_j}}.
\end{array}
\end{equation}
Squaring this inequality gives the desired estimate (\ref{fdfhn:eq:specconvest}).\qed\\

\begin{lemma}\label{fdfhn:Lemma318BDFdeel2a}  Consider the setting of Lemma \ref{fdfhn:Lemma316BDF} and assume that the triplet $(\mathcal{G},P^-,P^+)$ satisfies (\asref{fdfhn:aannamesconstantena}{\text{HS}3_{\overline{r}}}). There exist positive constants $\mu$, $C_3$, $C_4$ and $C_5$ so that the bound
\begin{equation}\label{fdfhn:eq137bdfa}
\begin{array}{lcl}
M_j^{-1}\sum\limits_{\xi\in M_j^{-1}\Z:|\xi|\leq \mu}|\Phi_j(\cdot,\xi)|_{\ell_{q,\perp}^2}^2&\geq &C_3\nrm{\Phi_j}_{ {\H}_{M_j}}^2 -C_4\nrm{\Psi_j}_{ {\H}_{M_j}}^2 -C_5M_j^{-1}\nrm{ {\D}_{k,M_j}\Phi_j}_{ {\H}_{M_j}}^2
\end{array}
\end{equation}
holds for all $j\in\N$.\end{lemma}
\textit{Proof.} Invoking Lemma \ref{fdfhn:Cor315BDF} and Lemma \ref{fdfhn:Lemma313BDF}, we can estimate
\begin{equation}\label{fdfhn:eq:specconv1}
\begin{array}{lcl}
\ip{\Psi_j,\Phi_j}_{ {\H}_{M_j}}&=&\ip{[\K_{k,M_j}+\delta]\Phi_j,\Phi_j}_{ {\H}_{M_j}}\\[0.2cm]
&= & \overline{c}_0\ip{ {\D}_{k,M_j}\Phi_j,\Phi_j}_{ {\H}_{M_j}} -\ip{\Delta_{M_j}\Phi_j,\Phi_j}_{\H_{M_j}}\\[0.2cm]
&&\qquad -\ip{D\mathcal{G}\big(\pi_{\H_{M_j}}\overline{U}_0;\overline{r}\big)\Phi_j,\Phi_j}_{\H_{M_j}}+\delta\nrm{\Phi_j}_{ {\H}_{M_j}}^2\\[0.2cm]
&\geq & \overline{c}_0\ip{ {\D}_{k,M_j}\Phi_j,\Phi_j}_{ {\H}_{M_j}}-\ip{D\mathcal{G}\big(\pi_{\H_{M_j}}\overline{U}_0;\overline{r}\big)\Phi_j,\Phi_j}_{\H_{M_j}}\\[0.2cm]
&\geq & -C_2M_j^{-1}\nrm{ {\D}_{k,M_j}\Phi_j}_{ {\H}_{M_j}}^2-\ip{D\mathcal{G}\big(\pi_{\H_{M_j}}\overline{U}_0;\overline{r}\big)\Phi_j,\Phi_j}_{\H_{M_j}}
\end{array}
\end{equation}
for some $C_2>1$. Since $-D\mathcal{G}\big(P^{\pm};\overline{r}\big)$ is positive definite and $-D\mathcal{G}$ is continuous, we can choose $\mu>0$ and $a>0$ in such a way that the matrix
\begin{equation}
    \begin{array}{lcl}
         B(\xi)&=&-D\mathcal{G}\big(\overline{U}_0(\xi);\overline{r}\big)-a
    \end{array}
\end{equation}
is positive definite for all $|\xi|\geq\mu$. Using the definition of this matrix and writing
\begin{equation}\label{fdfhn:eq:mathcalI}
\begin{array}{lcl}
\I&=&(\nrm{D\mathcal{G}\big(\overline{U}_0;\overline{r}\big)}_\infty+a)M_j^{-1}\sum\limits_{\xi\in M_j^{-1}\Z:|\xi|\leq \mu}|\Phi_j(\cdot,\xi)|_{\ell_{q,\perp}^2}^2,
\end{array}
\end{equation}
we can estimate
\begin{equation}\label{eq:fdfhn:eq:dgafsch}
    \begin{array}{lcl}
         -\ip{D\mathcal{G}\big(\pi_{\H_{M_j}}\overline{U}_0;\overline{r}\big)\Phi_j,\Phi_j}_{\H_{M_j}}&=&a\nrm{\Phi_j}_{\H_{M_j}}^2-\ip{B\Phi_j,\Phi_j}_{\H_{M_j}}\\[0.2cm]
         &\geq &a\nrm{\Phi_j}_{\H_{M_j}}^2-M_j^{-1}\sum\limits_{\xi\in M_j^{-1}\Z}|B(\xi)\Phi_j(\cdot,\xi)|_{\ell_{q,\perp}^2}^2\\[0.2cm]
         &\geq &a\nrm{\Phi_j}_{\H_{M_j}}^2-\mathcal{I}.
    \end{array}
\end{equation}
In particular, we can combine (\ref{fdfhn:eq:specconv1}) and (\ref{eq:fdfhn:eq:dgafsch}) to obtain
\begin{equation}\label{fdfhn:eq:specconv12}
    \begin{array}{lcl}
         \ip{\Psi_j,\Phi_j}_{ {\H}_{M_j}}&\geq & a\nrm{\Phi_j}_{\H_{M_j}}^2-\mathcal{I} -C_2M_j^{-1}\nrm{ {\D}_{k,M_j}\Phi_j}_{ {\H}_{M_j}}^2.
    \end{array}
\end{equation}
We can hence rearrange (\ref{fdfhn:eq:specconv12}) and estimate
\begin{equation}
\begin{array}{lcl}
\I&\geq & a\nrm{\Phi_j}_{ {\H}_{M_j}}^2-C_2M_j^{-1}\nrm{ {\D}_{k,M_j}\Phi_j}_{ {\H}_{M_j}}^2-\ip{\Psi_j,\Phi_j}_{ {\H}_{M_j}}\\[0.2cm]
&\geq & \frac{a}{2}\nrm{\Phi_j}_{ {\H}_{M_j}}^2-\frac{2}{a}\nrm{\Psi_j}_{ {\H}_{M_j}}^2-C_2M_j^{-1}\nrm{ {\D}_{k,M_j}\Phi_j}_{ {\H}_{M_j}}^2
,
\end{array}
\end{equation}
which yields the desired bound.\\

\begin{lemma}\label{fdfhn:Lemma318BDFdeel2b}  Consider the setting of Lemma \ref{fdfhn:Lemma316BDF} and assume that the triplet $(\mathcal{G},P^-,P^+)$ satisfies (\asref{fdfhn:aannamesconstantenb}{\text{HS}3_{\overline{r}}}). Then there exist positive constants $\mu$, $C_3$, $C_4$ and $C_5$ so that the bound
\begin{equation}\label{fdfhn:eq137bdfb}
\begin{array}{lcl}
M_j^{-1}\sum\limits_{\xi\in M_j^{-1}\Z:|\xi|\leq \mu}|\Phi_j(\cdot,\xi)|_{\ell_{q,\perp}^2}^2&\geq &C_3\nrm{\Phi_j}_{ {\H}_{M_j}}^2 -C_4\nrm{\Psi_j}_{ {\H}_{M_j}}^2 \\[0.2cm]
&&\qquad-C_5M_j^{-1}\nrm{ {\D}_{k,M_j}\Phi_j}_{ {\H}_{M_j}}^2
\end{array}
\end{equation}
holds for all $j\in\N$.\end{lemma}
\textit{Proof.} Recall the proportionality constant $\Gamma>0$ from (\asref{fdfhn:aannamesconstantenb}{\text{HS}3_{\overline{r}}}). In particular, upon writing
\begin{equation}
\begin{array}{lcl}
D\mathcal{G}&=&\left(\begin{array}{lcl}
D\mathcal{G}^{[1,1]}& D\mathcal{G}^{[1,2]}\\
D\mathcal{G}^{[2,1]}& D\mathcal{G}^{[2,2]}
\end{array}\right),
\end{array}
\end{equation} 
we have $D\mathcal{G}^{[1,2]}=-\Gamma (D\mathcal{G}^{[2,1]})^T$. For each $M\in\mathcal{M}_q$, we introduce the decomposition
\begin{equation}
    \begin{array}{lcl}
         \mathcal{H}_M &= &\mathcal{H}^{[1]}_M \times \mathcal{H}^{[2]}_M,
    \end{array}
\end{equation}
which splits every $\Phi=(\phi,\theta)\in \mathcal{H}_M$ in such a way that $\phi\in \mathcal{H}^{[1]}_M$ contains the first $d_{\mathrm{diff}}$ components of $\Phi$, while $\theta\in \mathcal{H}^{[2]}_M$ contains the other $d-d_{\mathrm{diff}}$ components. For each $j\geq 0$ we write $\Phi_j=(\phi_j,\theta_j)$ and $\Psi_j=(\psi_j,\chi_j)$ with $\phi_j,\psi_j\in\mathcal{H}^{[1]}_{M_j}$ and $\theta_j,\chi_j\in \mathcal{H}^{[2]}_{M_j} $. \\

Using this decomposition, we can expand the inner product as
\begin{equation}\label{fdfhn:eq:matrixmultiply}
    \begin{array}{lcl}
         -\ip{D\mathcal{G}\big(\pi_{\H_{M_j}}\overline{U}_0;\overline{r}\big)\Phi_j,\Phi_j}_{\H_{M_j}}&=&-\ip{D\mathcal{G}^{[1,1]}(\pi_{\H_{M_j}}\overline{U}_0)\phi_j,\phi_j}_{\mathcal{H}^{[1]}_{M_j}}+\mathcal{C}\\[0.2cm] &&\qquad-\ip{D\mathcal{G}^{[2,2]}(\pi_{\H_{M_j}}\overline{U}_0)\theta_j,\theta_j}_{\mathcal{H}^{[2]}_{M_j}},
    \end{array}
\end{equation}
where we have introduced the cross-terms
\begin{equation}
    \begin{array}{lcl}
        \mathcal{C}&:=&-\ip{D\mathcal{G}^{[1,2]}(\pi_{\H_{M_j}}\overline{U}_0)\theta_j,\phi_j}_{\mathcal{H}^{[1]}_{M_j}} -\ip{D\mathcal{G}^{[2,1]}(\pi_{\H_{M_j}}\overline{U}_0)\phi_j,\theta_j}_{\mathcal{H}^{[2]}_{M_j}}.
    \end{array}
\end{equation}
Recalling $D\mathcal{G}^{[1,2]}=-\Gamma (D\mathcal{G}^{[2,1]})^T$ and exploiting the identity
\begin{equation}
\begin{array}{lcl}
\chi_j&=&\overline{c}_0\D_{k,M_j}\theta_j-D\mathcal{G}^{[2,1]}(\pi_{\H_{M_j}}\overline{U}_0)\phi_j-D\mathcal{G}^{[2,2]}(\pi_{\H_{M_j}}\overline{U}_0)\theta_j+\delta \theta_j,
\end{array}
\end{equation}
we can rewrite the cross-terms to obtain 
\begin{equation}\label{fdfhn:eq:crossterms}
\begin{array}{lcl}
\mathcal{C}&=&-\ip{D\mathcal{G}^{[1,2]}(\pi_{\H_{M_j}}\overline{U}_0)\theta_j,\phi_j}_{\mathcal{H}^{[1]}_{M_j}} -\ip{D\mathcal{G}^{[2,1]}(\pi_{\H_{M_j}}\overline{U}_0)\phi_j,\theta_j}_{\mathcal{H}^{[2]}_{M_j}}\\[0.2cm]
&= & -(1-\Gamma)\ip{D\mathcal{G}^{[2,1]}(\pi_{\H_{M_j}}\overline{U}_0)\phi_j,\theta_j}_{\mathcal{H}^{[2]}_{M_j}}\\[0.2cm]
&=&(\Gamma-1)\ip{\overline{c}_0\D_{k,M_j}\theta_j-D\mathcal{G}^{[2,2]}(\pi_{\H_{M_j}}\overline{U}_0)\theta_j+\delta \theta_j-\chi_j,\theta_j}_{\mathcal{H}^{[2]}_{M_j}}.
\end{array}
\end{equation} 
The identities (\ref{fdfhn:eq:matrixmultiply}) and (\ref{fdfhn:eq:crossterms}) allow us to expand the inner product
\begin{equation}
\begin{array}{lcl}
\ip{\Psi_j,\Phi_j}_{ {\H}_{M_j}}&=&\ip{[\K_{k,M_j}+\delta]\Phi_j,\Phi_j}_{ {\H}_{M_j}}\\[0.2cm]
&= & \overline{c}_0\ip{ {\D}_{k,M_j}\Phi_j,\Phi_j}_{ {\H}_{M_j}} -\ip{\Delta_{M_j}\Phi_j,\Phi_j}_{\H_{M_j}}\\[0.2cm]
&&\qquad -\ip{D\mathcal{G}\big(\pi_{\H_{M_j}}\overline{U}_0;\overline{r}\big)\Phi_j,\Phi_j}_{\H_{M_j}}+\delta\nrm{\Phi_j}_{ {\H}_{M_j}}^2\\[0.2cm]
&= & \overline{c}_0\ip{\D_{k,M_j}\phi_j,\phi_j}_{\mathcal{H}^{[1]}_{M_j}}+\Gamma\overline{c}_0\ip{\D_{k,M_j}\theta_j,\theta_j}_{\mathcal{H}^{[2]}_{M_j}} -\ip{\Delta_{M_j}\Phi_j,\Phi_j}_{\H_{M_j}}\\[0.2cm]
&&\qquad -\ip{D\mathcal{G}^{[1,1]}(\pi_{\H_{M_j}}\overline{U}_0)\phi_j,\phi_j}_{\mathcal{H}^{[1]}_{M_j}} -\Gamma\ip{D\mathcal{G}^{[2,2]}(\pi_{\H_{M_j}}\overline{U}_0)\theta_j,\theta_j}_{\mathcal{H}^{[2]}_{M_j}}\\[0.2cm]
&&\qquad-(\Gamma-1)\ip{\chi_j,\theta_j}_{\mathcal{H}^{[2]}_{M_j}}+\delta\nrm{\phi_j}_{\mathcal{H}^{[1]}_{M_j}}^2+\delta\Gamma\nrm{\theta_j}_{\mathcal{H}^{[2]}_{M_j}}.\\[0.2cm]
\end{array}
\end{equation}
As such, we can use Lemma \ref{fdfhn:Cor315BDF} and Lemma \ref{fdfhn:Lemma313BDF} to estimate
\begin{equation}\label{fdfhn:eq:specconvestderiv}
\begin{array}{lcl}
\ip{\Psi_j,\Phi_j}_{ {\H}_{M_j}}
&\geq &  \overline{c}_0\ip{\D_{k,M_j}\phi_j,\phi_j}_{\mathcal{H}^{[1]}_{M_j}}+\Gamma\overline{c}_0\ip{\D_{k,M_j}\theta_j,\theta_j}_{\mathcal{H}^{[2]}_{M_j}} \\[0.2cm]
&&\qquad -\ip{D\mathcal{G}^{[1,1]}(\pi_{\H_{M_j}}\overline{U}_0)\phi_j,\phi_j}_{\mathcal{H}^{[1]}_{M_j}}-\Gamma\ip{D\mathcal{G}^{[2,2]}(\pi_{\H_{M_j}}\overline{U}_0)\theta_j,\theta_j}_{\mathcal{H}^{[2]}_{M_j}}\\[0.2cm]
&&\qquad -(\Gamma+1)\nrm{\chi_j}_{\mathcal{H}^{[2]}_{M_j}}\nrm{\theta_j}_{\mathcal{H}^{[2]}_{M_j}}\\[0.2cm]
&\geq &  -(1+\Gamma)C_2M_j^{-1}\nrm{ {\D}_{k,M_j}\Phi_j}_{ {\H}_{M_j}}^2\\[0.2cm] &&\qquad -\ip{D\mathcal{G}^{[1,1]}(\pi_{\H_{M_j}}\overline{U}_0)\phi_j,\phi_j}_{\mathcal{H}^{[1]}_{M_j}}-\Gamma\ip{D\mathcal{G}^{[2,2]}(\pi_{\H_{M_j}}\overline{U}_0)\theta_j,\theta_j}_{\mathcal{H}^{[2]}_{M_j}}\\[0.2cm]
&&\qquad -(\Gamma+1)\nrm{\chi_j}_{\mathcal{H}^{[2]}_{M_j}}\nrm{\theta_j}_{\mathcal{H}^{[2]}_{M_j}}
\end{array}
\end{equation}
for some $C_2>1$.\\

 Since $-D\mathcal{G}^{[1,1]}(P^{\pm})$ and $-D\mathcal{G}^{[2,2]}(P^{\pm})$ are positive definite and $-D\mathcal{G}$ is continuous, we can choose $\mu>0$ and $a>0$ to be positive constants such that the matrices
\begin{equation}
    \begin{array}{lclcl}
         B_1(\xi)&=&-D\mathcal{G}^{[1,1]}\big(\overline{U}_0(\xi)\big)-a,\qquad B_2(\xi)&=&-\Gamma D\mathcal{G}^{[2,2]}\big(\overline{U}_0(\xi)\big)-a
    \end{array}
\end{equation}
are positive definite for all $|\xi|\geq\mu$. Defining $\I$ as in (\ref{fdfhn:eq:mathcalI}), this allows us to estimate
\begin{equation}\label{fdfhn:eq:dg11afsch}
    \begin{array}{lcl}
         -\ip{D\mathcal{G}^{[1,1]}(\pi_{\H_{M_j}}\overline{U}_0)\phi_j,\phi_j}_{\mathcal{H}^{[1]}_{M_j}}&= &a\nrm{\phi_j}_{\mathcal{H}^{[1]}_{M_j}}^2+\ip{B_1\phi_j,\phi_j}_{\mathcal{H}^{[1]}_{M_j}}\\[0.2cm]
         &\geq &a\nrm{\Phi_j}_{ {\H}_{M_j}}^2-M_j^{-1}\sum\limits_{\xi\in M_j^{-1}\Z}|B_1(\xi)\phi_j(\cdot,\xi)|_{\ell_{q,\perp}^2}^2\\[0.2cm]
         &\geq &a\nrm{\Phi_j}_{ {\H}_{M_j}}^2-\mathcal{I},
    \end{array}
\end{equation}
together with
\begin{equation}\label{fdfhn:eq:dg22afsch}
    \begin{array}{lcl}
         -\Gamma\ip{D\mathcal{G}^{[2,2]}(\pi_{\H_{M_j}}\overline{U}_0)\theta_j,\theta_j}_{\mathcal{H}^{[2]}_{M_j}}&\geq &a\nrm{\theta_j}_{\mathcal{H}^{[2]}_{M_j}}^2-\Gamma\mathcal{I}.
    \end{array}
\end{equation}
Combining the estimates (\ref{fdfhn:eq:specconvestderiv}), (\ref{fdfhn:eq:dg11afsch}) and (\ref{fdfhn:eq:dg22afsch}) yields the bound
\begin{equation}
\begin{array}{lcl}
\ip{\Psi_j,\Phi_j}_{ {\H}_{M_j}}&\geq & a\nrm{\Phi_j}_{ {\H}_{M_j}}^2-(1+\Gamma)\mathcal{I} -(1+\Gamma)C_2M_j^{-1}\nrm{ {\D}_{k,M_j}\Phi_j}_{ {\H}_{M_j}}^2\\[0.2cm] &&\qquad
-(\Gamma+1)\nrm{\chi_j}_{\mathcal{H}^{[2]}_{M_j}}\nrm{\theta_j}_{\mathcal{H}^{[2]}_{M_j}}.\end{array}
\end{equation}
Hence we obtain
\begin{equation}
\begin{array}{lcl}
(1+\Gamma)\I&\geq & a\nrm{\Phi_j}_{ {\H}_{M_j}}^2-(1+\Gamma) C_2M_j^{-1}\nrm{ {\D}_{k,M_j}\Phi_j}_{ {\H}_{M_j}}^2\\[0.2cm] &&\qquad -\ip{\Psi_j,\Phi_j}_{ {\H}_{M_j}}-(\Gamma+1)\nrm{\chi_j}_{\mathcal{H}^{[2]}_{M_j}}\nrm{\theta_j}_{\mathcal{H}^{[2]}_{M_j}}\\[0.2cm]
&\geq & \frac{a}{2}\nrm{\Phi_j}_{ {\H}_{M_j}}^2-\big(\frac{1}{a}+\frac{\Gamma+1}{a}\big)\nrm{\Psi_j}_{ {\H}_{M_j}}^2-(1+\Gamma)C_2M_j^{-1}\nrm{ {\D}_{k,M_j}\Phi_j}_{ {\H}_{M_j}}^2
,
\end{array}
\end{equation}
which yields the desired bound. \qed\\

\textit{Proof of Proposition \ref{fdfhn:Lemma318BDF}.} Rescaling (\ref{fdfhn:eq:specconvest}) yields
\begin{equation}
\begin{array}{lcl}
0&\geq &\frac{C_3}{\overline{c}_0^2+2C_1}\Big[\overline{c}_0^2\nrm{ {\D}_{k,M_j}\Phi_j}_{ {\H}_{M_j}}^2-2C_1\nrm{\Phi_j}_{ {\H}_{M_j}}^2-2\nrm{\Psi_j}_{ {\H}_{M_j}}^2\Big],
\end{array}
\end{equation}
which can be added to (\ref{fdfhn:eq137bdfa}) or (\ref{fdfhn:eq137bdfb}) to obtain
\begin{equation}
\begin{array}{lcl}
M_j^{-1}\sum\limits_{\xi\in M_j^{-1}\Z:|\xi|\leq \mu}|\Phi_j(\cdot,\xi)|_{\ell_{q,\perp}^2}^2&\geq & C_3\nrm{\Phi_j}_{ {\H}_{M_j}}^2-C_4\nrm{\Psi_j}_{ {\H}_{M_j}}^2-C_5M_j^{-1}\nrm{ {\D}_{k,M_j}\Phi_j}_{ {\H}_{M_j}}^2\\[0.2cm]
&&\qquad +\frac{C_3}{\overline{c}_0^2+2C_1}\Big[\overline{c}_0^2\nrm{ {\D}_{k,M_j}\Phi_j}_{ {\H}_{M_j}}^2 -2C_1\nrm{\Phi_j}_{ {\H}_{M_j}}^2-2\nrm{\Psi_j}_{ {\H}_{M_j}}^2\Big]\\[0.2cm]
&=& \frac{\overline{c}_0^2C_3}{\overline{c}_0^2+2C_1}\Big[\nrm{ {\D}_{k,M_j}\Phi_j}_{ {\H}_{M_j}}^2 + \nrm{\Phi_j}_{ {\H}_{M_j}}^2\Big]\\[0.2cm]
&&\qquad -\Big[C_4+\frac{C_3}{\overline{c}_0^2+2C_1}\Big]\nrm{\Psi_j}_{ {\H}_{M_j}}^2-C_5M_j^{-1}\nrm{ {\D}_{k,M_j}\Phi_j}_{ {\H}_{M_j}}^2.
\end{array}\end{equation}
Remembering that $\nrm{\Phi_j}_{ {\H}_{k,M_j}^1}=1$, we can pick constants $C_6>0$, $C_7>0$ and $C_8>0$, which all are independent of $0<\delta<\delta_0$, in such a way that
\begin{equation}
\begin{array}{lcl}
M_j^{-1}\sum\limits_{\xi\in M_j^{-1}\Z:|\xi|\leq \mu}|\Phi_j(\cdot,\xi)|_{\ell_{q,\perp}^2}^2&\geq & C_6-C_7\nrm{\Psi_j}_{ {\H}_{M_j}}^2-C_8M_j^{-1}.
\end{array}
\end{equation}
The strong convergence $\I_{M_j}^0\Phi_j\rightarrow \Phi_*\in  {L}^2([-\mu-1,\mu+1];\ell_{q,\perp}^2)$ now yields the limiting behaviour
\begin{equation}
\begin{array}{lcl}
M_j^{-1}\sum\limits_{\xi\in M_j^{-1}\Z:|\xi|\leq \mu}|\Phi_j(\cdot,\xi)|_{\ell_{q,\perp}^2}^2&=&\int_{-\mu}^{\mu+M_j^{-1}}\Big|[\I_{M_j}^0\Phi_j](\cdot,\xi)\Big|_{\ell_{q,\perp}^2}^2d\xi\\[0.2cm]
&\leq &\int_{-\mu}^{\mu+1}\Big|[\I_{M_j}^0\Phi_j](\cdot,\xi)\Big|_{ {\ell}_{q,\perp}^2}^2d\xi\\[0.2cm]
&\rightarrow &\int_{-\mu}^{\mu+1}\Big|\Phi_*(\cdot,\xi)\Big|_{ {\ell}_{q,\perp}^2}^2d\xi,
\end{array}
\end{equation}
as $j\rightarrow\infty$. In view of the bound $\limsup_{j\rightarrow\infty}\nrm{\Psi_j}_{ {\H}_{M_j}}^2\leq \kappa(\delta)^2$, this gives the desired inequality
\begin{equation}
\begin{array}{lclcl}
\nrm{\Phi_*}_{ {H}^1(\R,\ell_{q,\perp}^2)}^2&\geq & \int_{-\mu}^{\mu+1}\Big|\Phi_*(\cdot,\xi)\Big|_{ {\ell}_{q,\perp}^2}^2d\xi
&\geq & C_6-C_7\kappa(\delta)^2.
\end{array}
\end{equation}
\qed\\

\subsection{Exponential decay}\label{fdfhn:sectionexpdecay}

In this section we set out to prove Proposition \ref{fdfhn:exponentieelvervalprop}. The main ingredient to establish this result is to show that for $0<\delta<\delta_0$ the map $(L_0+\delta)^{-1}$ maps $BC_{-\eta}^1(\R;\R^d)$ into the space
\begin{equation}
\begin{array}{lcl}
BC_{-\eta}^2(\R;\R^d)&=&\{ F\in BC_{-\eta}(\R;\R^d):\sup_{\xi\in\R}e^{-\eta|\xi|}[|F(\xi)|+|F'(\xi)|+|F''(\xi)|]<\infty\}.
\end{array}
\end{equation} 
This is not immediately clear, since if we have
\begin{equation}
    \begin{array}{lcl}
         F&=&(L_0+\delta)^{-1}G
    \end{array}
\end{equation}
with $G\in BC_{-\eta}^1(\R;\R^d)$, it is impossible to express $F$ as a local function of $G$ due to the infinite-range interactions. We first establish this result for the subspaces $H^1(\R;\R^d)$ and $H^2(\R;\R^d)$.
\begin{lemma}\label{fdfhn:exponentieelvervalbewijs1} Assume that (\asref{fdfhn:aannames}{\text{HS}1}) and (\asref{fdfhn:aannames2}{\text{HS}2}) are satisfied and pick $\overline{r}$ in such a way that (\asref{fdfhn:aannamesconstanten2}{\text{HS}}), (\asref{fdfhn:existentie}{\text{HW}}) and (\asref{fdfhn:spectralstability}{\text{HW}}) are satisfied.  Then for each $0<\delta<\delta_0$ and each $G\in H^1(\R;\R^d)$ we have
\begin{equation}
\begin{array}{lcl}
(L_0+\delta)^{-1}G&\in & H^2(\R;\R^d).
\end{array}
\end{equation}\end{lemma}
\textit{Proof.} Fix $0<\delta<\delta_0$ and $G\in H^1(\R;\R^d)$. Write $F=(L_0+\delta)^{-1}G\in H^1(\R;\R^d)$. Then we can rewrite the equation $(L_0+\delta)F=G$ in the form
\begin{equation}\label{fdfhn:inverseomschrijven}
\begin{array}{lcl}
\overline{c}_0F'&=& G+\Delta_0 F+D_U\mathcal{G}\big(\overline{U}_0;\overline{r}\big)F-\delta F.
\end{array}
\end{equation}
From this representation it immediately follows that $F'\in L^\infty(\R;\R^d)$. Differentiating both sides yields
\begin{equation}\label{fdfhn:secondderivativeexponentieel}
\begin{array}{lcl}
\overline{c}_0F''&=&G'+(\Delta_0 F)'+D_U\mathcal{G}\big(\overline{U}_0;\overline{r}\big)F'+D^2\mathcal{G}\big(\overline{U}_0;\overline{r})[\overline{U}_0',F]-\delta F'.
\end{array}
\end{equation}
Writing
\begin{equation}
\begin{array}{lcl}
F_n(x)&=&\tau\sum\limits_{m=1}^n\alpha_m\Big[F(x+m)+F(x-m)-2F(x)\Big]
\end{array}
\end{equation}
for $n\in\Z_{>0}$, we can compute
\begin{equation}
\begin{array}{lcl}
F_n'(x)&=&\tau\sum\limits_{m=1}^n\alpha_m\Big[F'(x+m)+F'(x-m)-2F'(x)\Big].
\end{array}
\end{equation}
This allows us to estimate 
\begin{equation}
\begin{array}{lcl}
|F_n'(x)-(\Delta_0F')(x)|&\leq &4\tau\sum\limits_{m=n+1}^\infty|\alpha_m|\nrm{F'}_\infty.
\end{array}
\end{equation}
In particular, the sequence $\{F_n'\}$ converges uniformly to $\Delta_0F'$, from which it follows that
\begin{equation}
\begin{array}{lclcl}
(\Delta_0 F)'(x)&=&\tau\sum\limits_{m=1}^\infty\alpha_m\Big[F'(x+m)+F'(x-m)-2F'(x)\Big]&=& (\Delta_0 F')(x).
\end{array}
\end{equation}
Since $F,G\in H^1(\R;\R^d)$, these considerations yield that $F''\in L^2(\R;\R^d)$, from which the desired result follows.\qed\\

We now turn to the desired exponential decay. The assumptions (\asref{fdfhn:existentie}{\text{HW}}) and (\asref{fdfhn:spectralstability}{\text{HW}}) yield the following useful properties of the operator $L_0$. 

\begin{lemma}\label{fdfhn:mapsintoitself} Assume that (\asref{fdfhn:aannames}{\text{HS}1}) and (\asref{fdfhn:aannames2}{\text{HS}2}) are satisfied and pick $\overline{r}$ in such a way that (\asref{fdfhn:aannamesconstanten2}{\text{HS}}), (\asref{fdfhn:existentie}{\text{HW}}) and (\asref{fdfhn:spectralstability}{\text{HW}}) are satisfied. Then the following properties hold for the LDE (\ref{fdfhn:ditishetalgemeneprobleem}) with $r=\overline{r}$.
\begin{enumerate}[label=(\roman*)]
\item\label{fdfhn:item1} The functions $\Phi_0^+$ and $\Phi_0^-$ together with their derivatives decay exponentially.

\item\label{fdfhn:item2} Upon introducing the spaces
\begin{equation}\begin{array}{lcl}X_0&=&\{F\in H^1(\R;\R^d): \ip{\Phi_0^-,F}_{ L^2(\R;\R^d)}=0\}\end{array}\end{equation}
and
\begin{equation}\begin{array}{lcl}Y_0&=&\{G\in  L^2(\R;\R^d): \ip{\Phi_0^-,G}_{ L^2(\R;\R^d)}=0\},\end{array}\end{equation}
the operator $L_0:X_0\rightarrow Y_0$ is invertible.\end{enumerate}
In addition, there exists a constant $\tilde{\eta}>0$ in such a way that for each $0<\eta<\tilde{\eta}$ the map $L_0$ maps $BC_{-\eta}^1(\R;\R^d)$ into $BC_{-\eta}(\R;\R^d)$.\end{lemma}
\textit{Proof.} The proof of the statements \ref{fdfhn:item1}-\ref{fdfhn:item2} follows the procedure described \cite[Lem. 4.15, 6.8, 6.9]{HJHFHNINFRANGE} and will hence be omitted. It hence suffices to prove that $\Delta_0$ maps $BC_{-\eta}(\R;\R^{d_\mathrm{diff}})$ into itself for $\eta$ small enough. Upon picking $F\in BC_{-\eta}(\R;\R^{d_\mathrm{diff}})$ and $K\in\R_{>0}$ in such a way that the bound
\begin{equation}
\begin{array}{lcl}
|F(\xi)|&\leq & K e^{-\eta|\xi|}
\end{array}
\end{equation}
holds, we estimate
\begin{equation}\label{fdfhn:Deltahexponentieelbegrensd}
\begin{array}{lcl}
|\Delta_0 F(\xi)|&\leq &\tau\sum\limits_{m>0}|\alpha_m|K\Big(e^{-\eta|\xi+m|}+e^{-\eta|\xi-m|}+2 e^{-\eta|\xi|}\Big)\\[0.2cm]
&\leq &\tau\sum\limits_{m>0}|\alpha_m|Ke^{-\eta|\xi|}\Big(2e^{\eta m}+2\Big).
\end{array}
\end{equation}
We can hence set $\tilde{\eta}= \nu$, where $\nu$ is defined in (\asref{fdfhn:aannames}{\text{HS}1}). A computation similar to the proof of \cite[Lem. 6.5]{HJHFHNINFRANGE} yields the continuity of $\Delta_0f$, from which the desired result follows.\qed\\

We now recall the notation $L_0^{\mathrm{qinv}}G$ that was introduced in \cite[Cor. 4.4]{HJHFHNINFRANGE} for the unique solution $F$ of the equation
\begin{equation}\begin{array}{lcl}L_0F&=&G-\frac{\ip{\Phi_0^-,G}_{L^2(\R;\R^d)}}{\ip{\Phi_0^-,\Phi_0^+}_{L^2(\R;\R^d)}}\Phi_0^+\end{array}\end{equation}
in the space $X_0$, which is given explicitly by
\begin{equation}\label{fdfhn:Lhqinv}\begin{array}{lcl}L_0^{\mathrm{qinv}}G&=&L_0^{-1}\Big[G-\frac{\ip{\Phi_0^-,G}_{L^2(\R;\R^d)}}{\ip{\Phi_0^-,\Phi_0^+}_{L^2(\R;\R^d)}}\Phi_0^+\Big].\end{array}\end{equation}
The proof of \cite[Prop. 5.2]{HJHFHNINFRANGE} provides the representation
\begin{equation}\label{fdfhn:Lhdeltainverse}
\begin{array}{lcl}
(L_0+\delta)^{-1}G&=&\delta^{-1}\frac{\ip{\Phi_0^-,G}_{L^2(\R;\R^d)}}{\ip{\Phi_0^-,\Phi_0^+}_{L^2(\R;\R^d)}}\Phi_0^++[I+\delta L_0^{-1}]^{-1} L_0^{\mathrm{qinv}}G
\end{array}
\end{equation} 
for each $0<\delta<\delta_0$ and each $G\in L^2(\R;\R^d)$. In addition, we can use Lemma \ref{fdfhn:mapsintoitself} to pick constants $\tilde{K}>0$ and $\tilde{\alpha}>0$ in such a way that
\begin{equation}\label{fdfhn:exponentieelvervalPhihplus}
\begin{array}{lcl}
|\Phi_0^+(x)|&\leq & \tilde{K}e^{-\tilde{\alpha}|x|}
\end{array}
\end{equation}
holds for all $x\in\R$.
Let $\tilde{\eta}>0$ be the constant from Lemma \ref{fdfhn:mapsintoitself}. Using \cite[Lem. 6.6]{HJHFHNINFRANGE}, which is a generalization of \cite[Prop. 5.3]{MPA}, we can pick constants $K_1>0$ and $0<\alpha\leq\min\{\tilde{\eta},\tilde{\alpha}\}$ in such a way that
\begin{equation}\label{fdfhn:expvervalquinv}
\begin{array}{lcl}
|L_0^{\mathrm{qinv}}G(x)|&\leq & K_1e^{-\alpha|x|}\nrm{L_0^{\mathrm{qinv}}G}_{\infty}+K_1\int_{-\infty}^\infty e^{-\alpha|x-y|}\Big|G(y)-\frac{\ip{\Phi_0^-,G}_{L^2(\R;\R^d)}}{\ip{\Phi_0^-,\Phi_0^+}_{L^2(\R;\R^d)}}\Phi_0^+(y)\Big|dy
\end{array}
\end{equation}
holds for each $G\in L^2(\R;\R^d)$. The following three results use (\ref{fdfhn:Lhdeltainverse}) and (\ref{fdfhn:expvervalquinv}) to establish the desired pointwise bound for $(L_0+\delta)^{-1}$.

\begin{lemma}\label{fdfhn:exponentieelvervalbewijs2} Assume that (\asref{fdfhn:aannames}{\text{HS}1}) and (\asref{fdfhn:aannames2}{\text{HS}2}) are satisfied and pick $\overline{r}$ in such a way that (\asref{fdfhn:aannamesconstanten2}{\text{HS}}), (\asref{fdfhn:existentie}{\text{HW}}) and (\asref{fdfhn:spectralstability}{\text{HW}}) are satisfied. Recall the constant $\alpha>0$ from (\ref{fdfhn:expvervalquinv}) and fix $0<\eta\leq \alpha$. Then there exists a constant $K>0$ so that for each $G \in BC_{-\eta}^1(\R;\R^d)$ we have the bound
\begin{equation}
\begin{array}{lcl}
|L_0^{\mathrm{qinv}}G(x)|&\leq & K\nrm{G}_{BC_{-\eta}(\R;\R^d)}e^{-\eta|x|}.
\end{array}
\end{equation}\end{lemma}
\textit{Proof.} Pick $0<\eta\leq \alpha$ and $G \in BC_{-\eta}^1(\R;\R^d)$. Recalling (\ref{fdfhn:expvervalquinv}), we can estimate
\begin{equation}\label{fdfhn:inftynormafschatten}
\begin{array}{lcl}
\nrm{L_0^{\mathrm{qinv}}G}_{\infty}&\leq &\nrm{L_0^{\mathrm{qinv}}G}_{H^1(\R;\R^d)}\\[0.2cm]
&\leq &\nrm{L_0^{-1}}_{\mathcal{L}(Y_0,X_0)}\nrm{G-\frac{\ip{\Phi_0^-,G}_{L^2(\R;\R^d)}}{\ip{\Phi_0^-,\Phi_0^+}_{L^2(\R;\R^d)}}\Phi_0^+}_{L^2(\R;\R^d)}\\[0.2cm]
&\leq & \nrm{L_0^{-1}}_{\mathcal{L}(Y_0,X_0)}\Big(1+\frac{1}{|\ip{\Phi_0^-,\Phi_0^+}_{L^2(\R;\R^d)}|}\Big)\nrm{G}_{L^2(\R;\R^d)}.
\end{array}
\end{equation}
Combining these estimates yields the bound
\begin{equation}
\begin{array}{lcl}
|L_0^{\mathrm{qinv}}G(x)|&\leq &K_1e^{-\alpha|x|}\nrm{L_0^{-1}}_{\mathcal{L}(Y_0,X_0)}\Big(1+\frac{1}{|\ip{\Phi_0^-,\Phi_0^+}_{L^2(\R;\R^d)}|}\Big)\nrm{G}_{L^2(\R;\R^d)}\\[0.2cm]
&&\qquad+K_1\Big(\nrm{G}_{BC_{-\eta}(\R,\R^d)}+\tilde{K}\frac{1}{|\ip{\Phi_0^-,\Phi_0^+}_{L^2(\R;\R^d)}|}\nrm{G}_{L^2(\R;\R^d)}\Big)\int_{-\infty}^\infty e^{-\alpha|x-y|}e^{-\eta|x|}dy\\[0.2cm]
&\leq & K_1e^{-\alpha|x|}\nrm{L_0^{-1}}_{\mathcal{L}(Y_0,X_0)}\Big(1+\frac{1}{|\ip{\Phi_0^-,\Phi_0^+}_{L^2(\R;\R^d)}|}\Big)\nrm{G}_{L^2(\R;\R^d)}\\[0.2cm]
&&\qquad+K_1\Big(\nrm{G}_{BC_{-\eta}(\R,\R^d)}+\tilde{K}\frac{1}{|\ip{\Phi_0^-,\Phi_0^+}_{L^2(\R;\R^d)}|}\nrm{G}_{L^2(\R;\R^d)}\Big)e^{-\eta|x|}\\[0.2cm]
&\leq & K_2\Big(\nrm{G}_{BC_{-\eta}(\R;\R^d)}+\nrm{G}_{L^2(\R;\R^d)}\Big)e^{-\eta|x|}.\end{array}\end{equation}
Finally, we note that $\nrm{G}_{L^2(\R;\R^d)}\leq K_3\nrm{G}_{BC_{-\eta}(\R;\R^d)}$ for some constant $K_3>0$, which implies the desired result.\qed\\

\begin{lemma}\label{fdfhn:exponentieelvervalbewijs3} Consider the setting of Lemma \ref{fdfhn:exponentieelvervalbewijs2}. Then there exist constants $0<\delta_*\leq\delta_0$ and $K>0$ so that for each $0<\delta<\delta_*$ and each $G \in BC_{-\eta}^1(\R;\R^d)$ we have the bound
\begin{equation}
\begin{array}{lcl}
|[I+\delta L_0^{-1}]^{-1}L_0^{\mathrm{qinv}}G(x)|
&\leq &K\nrm{G}_{BC_{-\eta}(\R;\R^d)}e^{-\eta|x|}.
\end{array}
\end{equation}
\end{lemma}
\textit{Proof.} Pick $G \in BC_{-\eta}^1(\R;\R^d)$. For $n\in\Z_{>0}$ a calculation similar to (\ref{fdfhn:inftynormafschatten}) yields
\begin{equation}
\begin{array}{lcl}
\nrm{(L_0^{-1})^nL_0^{\mathrm{qinv}}G}_{\infty}&\leq &\nrm{L_0^{-1}}_{\mathcal{L}(Y_0,X_0)}^{n+1}\Big(1+\frac{1}{|\ip{\Phi_0^-,\Phi_0^+}_{L^2(\R;\R^d)}|}\Big)\nrm{G}_{L^2(\R;\R^d)}.
\end{array}
\end{equation} 
Using \cite[Lem. 6.6]{HJHFHNINFRANGE} and Lemma \ref{fdfhn:exponentieelvervalbewijs2}, we obtain
\begin{equation}
\begin{array}{lcl}
|L_0^{-1}L_0^{\mathrm{qinv}}G(x)|&\leq & K_1e^{-\alpha|x|}\nrm{L_0^{-1}L_0^{\mathrm{qinv}}G}_{\infty}+K_1\int_{-\infty}^\infty e^{-\alpha|x-y|}\Big|L_0^{\mathrm{qinv}}G\Big|dy\\[0.2cm]
&\leq & K_1 \Big(1+\nrm{L_0^{-1}}_{\mathcal{L}(Y_0,X_0)}\Big)K\nrm{G}_{BC_{-\eta}(\R;\R^d)}e^{-\eta|x|}.
\end{array}
\end{equation}
Continuing in this fashion we see that the estimate
\begin{equation}
\begin{array}{lcl}
|(L_0^{-1})^nL_0^{\mathrm{qinv}}G(x)|&\leq & K_2^nK\nrm{G}_{BC_{-\eta}(\R;\R^d)}e^{-\eta|x|}
\end{array}
\end{equation}
holds for all $n\in\Z_{>0}$ and for some constant $K_2>0$. If we set
\begin{equation}
\begin{array}{lcl}
\delta_*&=&\min\Big\{\delta_0,\frac{1}{\nrm{L_0^{-1}}_{\mathcal{L}(Y_0,X_0)}\Big(1+\frac{1}{|\ip{\Phi_0^-,\Phi_0^+}_{L^2(\R;\R^d)}|}\Big)},\frac{1}{K_2}\Big\},
\end{array}
\end{equation}
then for each $n\in\Z_{>0}$ and each $0<\delta<\delta_*$ we have
\begin{equation}
\begin{array}{lcl}
\nrm{(-\delta)^n(L_0^{-1})^nL_0^{\mathrm{qinv}}G}_{\infty}&\leq &\frac{1}{2}\nrm{G}_{BC_{-\eta}(\R;\R^d)}.
\end{array}
\end{equation} 
In particular, it follows that
\begin{equation}
\begin{array}{lcl}
\sum\limits_{n=0}^\infty (-\delta)^n(L_0^{-1})^nL_0^{\mathrm{qinv}}G&\rightarrow & [I+\delta L_0^{-1}]^{-1}L_0^{\mathrm{qinv}}G
\end{array}
\end{equation}
in $H^1(\R;\R^d)$. Since $H^1(\R;\R^d)$-convergence implies pointwise convergence we see that
\begin{equation}
\begin{array}{lcl}
|[I+\delta L_0^{-1}]^{-1}L_0^{\mathrm{qinv}}G(x)|&=&|\sum\limits_{n=0}^\infty (-\delta)^n(L_0^{-1})^nL_0^{\mathrm{qinv}}G(x)|\\[0.2cm]
&\leq &\sum\limits_{n=0}^\infty \delta_*^n K K_2^n\nrm{G}_{BC_{-\eta}(\R;\R^d)}e^{-\eta|x|}\\[0.2cm]
&\leq &K_3\nrm{G}_{BC_{-\eta}(\R;\R^d)}e^{-\eta|x|}.
\end{array}
\end{equation}
\qed\\

\begin{corollary}\label{fdfhn:exponentieelvervalbewijs4} Consider the setting of Lemma \ref{fdfhn:exponentieelvervalbewijs3}. There exists a constant $K>0$ so that for each $0<\delta<\delta_*$ and each $G \in BC_{-\eta}^1(\R;\R^d)$ we have the bound
\begin{equation}
\begin{array}{lcl}
|(L_0+\delta)^{-1}G(x)|&\leq &K\delta^{-1}\nrm{G}_{BC_{-\eta}(\R;\R^d)}e^{-\eta|x|},\\[0.2cm]
|[(L_0+\delta)^{-1}G]'(x)|&\leq & K \delta^{-1}\nrm{G}_{BC_{-\eta}(\R;\R^d)}e^{-\eta|x|},
\\[0.2cm]
|[(L_0+\delta)^{-1}G]''(x)|&\leq & K \delta^{-1}\nrm{G}_{BC_{-\eta}^1(\R;\R^d)}e^{-\eta|x|}.
\end{array}
\end{equation}\end{corollary}
\textit{Proof.} Fix $0<\delta<\delta_*$ and $G\in BC_{-\eta}^1(\R;\R^d)$. Write $F=(L_0+\delta)^{-1}G$. The representation (\ref{fdfhn:Lhdeltainverse}) together with Lemma \ref{fdfhn:exponentieelvervalbewijs3} immediately yields the bound
\begin{equation}\label{fdfhn:eq:firstinverseexp}
\begin{array}{lcl}
|F(x)|&\leq &\delta^{-1}\frac{1}{|\ip{\Phi_0^-,\Phi_0^+}_{L^2(\R;\R^d)}|}\nrm{G}_{L^2(\R;\R^d)}\tilde{K}e^{-\tilde{\alpha}|x|}+ K\nrm{G}_{BC_{-\eta}(\R;\R^d)}e^{-\eta|x|}\\[0.2cm]
&\leq & \delta^{-1}K_2\nrm{G}_{BC_{-\eta}(\R;\R^d)}e^{-\eta|x|}.
\end{array}
\end{equation}
In addition, the representation (\ref{fdfhn:inverseomschrijven}) together with the bounds (\ref{fdfhn:Deltahexponentieelbegrensd}) and (\ref{fdfhn:eq:firstinverseexp}) yields that
\begin{equation}
\begin{array}{lcl}
|F'(x)|&\leq & K \delta^{-1}\nrm{G}_{BC_{-\eta}(\R;\R^d)}e^{-\eta|x|}
\end{array}
\end{equation}
for some constant $K>0$. Similarly, the representation (\ref{fdfhn:secondderivativeexponentieel}) yields the bound
\begin{equation}
\begin{array}{lcl}
|F''(x)|&\leq & K \delta^{-1}\nrm{G}_{BC_{-\eta}^1(\R;\R^d)}e^{-\eta|x|}.
\end{array}
\end{equation}
\qed\\
\textit{Proof of Proposition \ref{fdfhn:exponentieelvervalprop}.} Corollary \ref{fdfhn:exponentieelvervalbewijs4} implies the desired result.\qed\\

\appendix
%\appsection
\section{Auxiliary results}\label{fdfhn:section:appendix}
In this section we collect several useful results that we use throughout this paper. The first three results concern the sequence spaces $\Y_{M}$ and $\Y_{k,M}^1$ and their associated inner products (\ref{fdfhn:eq:defYM})-(\ref{fdfhn:eq:defYkm1}).
\begin{lemma}[{\cite[Lem. 3.1]{HJHBDF}}]\label{fdfhn:Lemma31BDF} Fix a pair of integers $1\leq k\leq 6$ and $q\geq 1$, together with a constant $\eta>0$. Then there exists a constant $C\geq 1$ for which the bounds
\begin{equation}
\begin{array}{lcl}
\nrm{\pi_{\mathcal{Y}_M}f}_{\mathcal{Y}_M}&\leq & C\nrm{f}_{BC_{-\eta}},\\[0.2cm]
\nrm{\pi_{\Y_{k,M}^1}g}_{\Y_{k,M}^1}&\leq & C\nrm{g}_{BC_{-\eta}^1}
\end{array}
\end{equation}
hold for all $M\in\mathcal{M}_{q}$ and all functions $f\in BC_{-\eta}(\R;\R)$ and $g\in BC_{-\eta}^1(\R;\R)$. 
\end{lemma}

\begin{lemma}[{\cite[Lem. 3.4]{HJHBDF}}]\label{fdfhn:Lemma34BDF} Fix an integer $q\geq 1$. Then there exists $C>1$ so that the bound
\begin{equation}
\begin{array}{lcl}
\left|\ip{f,g}_{L^2(\R;\R^d)}-\ip{\pi_{\mathcal{Y}_M}f,\pi_{\mathcal{Y}_M}g}_{ {\Y}_M}\right|&\leq & CM^{-1}\nrm{f}_{BC_{-\eta}^1(\R;\R^d)}\nrm{g}_{BC_{-\eta}^1(\R;\R^d)}
\end{array}
\end{equation}
holds for all $M\in \mathcal{M}_{q}$ and all functions $f,g\in BC_{-\eta}^1(\R;\R^d)$.
\end{lemma}

\begin{lemma}[{\cite[Lem. 3.5]{HJHBDF}}]\label{fdfhn:Lemma35BDF} Fix an integer $q\geq 1$. For any $M=\frac{p}{q}\in\mathcal{M}_{q}$, the operators $\mathcal{J}_M$ and $\mathcal{J}_{k,M}^1$ defined in (\ref{fdfhn:definitieJM}) are isometries between $\Y_M$ and $\H_M$ and between $\Y_{k,M}^1$ and $\H_{k,M}^1$ respectively.\end{lemma}

The following results can be seen as the fully discrete generalizations of the well-known facts
\begin{equation}
    \begin{array}{lclcl}
         \ip{u,u'}&=&0,\qquad \ip{u'',u}&\leq &0
    \end{array}
\end{equation}
that hold for smooth, localized functions $u$. When dealing solely with nearest-neighbour interactions as in \cite{HJHBDF} the inequality $\ip{\Delta_{M}\Phi,\Phi}_{ {\H}_M}\leq 0$ follows immediately from the Cauchy-Schwarz inequality. However, in our setting, some of the coefficients $\alpha_k$ may not be positive definite, preventing us from taking them out of the inner product. This motivates the indirect approach that is taken 
%causes major complications 
in the proof of Lemma \ref{fdfhn:Lemma313BDF}.
\begin{lemma}[{\cite[Cor. 3.15]{HJHBDF}}]\label{fdfhn:Cor315BDF} Fix a pair of integers $1\leq k\leq 6$ and $q\geq 1$. There exists a constant $K>1$ so that for all $M\in\M_{q}$ and all $\Phi\in \H_{k,M}^1$ we have the bound
\begin{equation}
\begin{array}{lcl}
\Big|\ip{\Phi,\D_{k,M}\Phi}_{\H_M}\Big|&\leq & KM^{-1}\nrm{\D_{k,M}\Phi}_{\H_M}^2.
\end{array}
\end{equation}
\end{lemma}

\begin{lemma}[{cf. \cite[Lem. 3.13]{HJHBDF}}]\label{fdfhn:Lemma313BDF} Assume that (\asref{fdfhn:aannames}{\text{HS}1}) is satisfied. Fix an integer $q\geq 1$ and pick $M\in\mathcal{M}_{q}$. Then the bound
\begin{equation}
\begin{array}{lcl}
\ip{\Delta_{M}\Phi,\Phi}_{ {\H}_M}&\leq &0
\end{array}
\end{equation}
holds for each $\Phi\in {\H}_M$.
\end{lemma}

\textit{Proof.} Pick $\Phi\in {\H}_M$ and define the stepwise interpolation function $\tilde{\Phi}\in L^2(\R;\R^d)$ by setting
\begin{equation}
\begin{array}{lcl}
\tilde{\Phi}\big(\xi+\zeta M^{-1}+\epsilon\big)&=& \Phi\big(\zeta,\xi\big)
\end{array}
\end{equation}
for $\xi\in M^{-1}\Z$, 
$\zeta\in q^{-1}\Z_q^{\circ}\cup\{0\}$ and $0\leq \epsilon<M^{-1}q^{-1}$. Upon recalling that
\begin{equation}
\begin{array}{lclcl}
\vartheta&=&\frac{p-nq}{q},\qquad 
nM^{-1}&=&1-\vartheta M^{-1}
\end{array}
\end{equation}
and observing that
\begin{equation}
\begin{array}{lclcl}
1&=&p\frac{q}{p}q^{-1}
&=&\big((p-nq)+nq\big)M^{-1}q^{-1},
\end{array}
\end{equation}
we may compute
\begin{equation}
    \begin{array}{lcl}
         \tilde{\Phi}\big(\xi+\zeta M^{-1}+m\big)&=&\tilde{\Phi}\big(\xi+mnM^{-1}+(\zeta+m(p-nq)q^{-1}) M^{-1}\big)\\[0.2cm]
         &=&\Phi\big(\zeta+m(p-nq)q^{-1},\xi+mnM^{-1}\big)\\[0.2cm]
         &=&\Phi\big(\zeta+m\vartheta,\xi+m-m\vartheta M^{-1}\big)
    \end{array}
\end{equation}
for arbitrary  $\xi\in M^{-1}\Z$, $\zeta\in q^{-1}\Z_q^\circ\cup\{0\}$ and $m\in\Z$. In particular, for $m\in\Z$ we obtain the identity 
\begin{equation}
\begin{array}{lcl}
\ip{T_0^m\tilde{\Phi},\tilde{\Phi}}_{L^2(\R;\R^d)}&=&q^{-1}M^{-1}\sum\limits_{\xi\in M^{-1}\Z}\sum\limits_{\zeta\in q^{-1}\Z_q^{\circ}\cup\{0\}}\big\langle\tilde{\Phi}\big(\xi+\zeta M^{-1}+m\big),\tilde{\Phi}\big(\xi+\zeta M^{-1}\big)\big\rangle_{\R^d}\\[0.2cm]
&=&q^{-1}M^{-1}\sum\limits_{\xi\in M^{-1}\Z}\sum\limits_{\zeta\in q^{-1}\Z_q^{\circ}\cup\{0\}}\big\langle\Phi\big(\zeta+m\vartheta,\xi+m-m\vartheta M^{-1}\big),\Phi\big(\zeta,\xi\big)\big\rangle_{\R^d}\\[0.2cm]
&=&\ip{T_{M}^m\Phi,\Phi}_{ {\H}_M}.

\end{array}
\end{equation}
We hence obtain
\begin{equation}
\begin{array}{lcl}
\ip{\Delta_0\tilde{\Phi},\tilde{\Phi}}_{L^2(\R;\R^d)}&=&\tau\sum\limits_{m>0}\alpha_m\big[\ip{T_0^m\tilde{\Phi},\tilde{\Phi}}_{L^2(\R;\R^d)}+\ip{T_0^{-m}\tilde{\Phi},\tilde{\Phi}}_{L^2(\R;\R^d)}-2\ip{\tilde{\Phi},\tilde{\Phi}}_{L^2(\R;\R^d)}\big]\\[0.2cm]
&=&\tau\sum\limits_{m>0}\alpha_m\big[\ip{T_{M}^m\Phi,\Phi}_{ {\H}_M}+\ip{T_{M}^{-m}\Phi,\Phi}_{ {\H}_M}-2\ip{\Phi,\Phi}_{ {\H}_M}\big]\\[0.2cm]
&=&\ip{\Delta_{M}\Phi,\Phi}_{ {\H}_M}.

\end{array}
\end{equation}
The desired result now follows from \cite[Lem. 3]{BatesInfRange}.\qed\\

We now show that $\mathcal{K}_{k,M}$ approaches $\overline{\mathcal{K}}_{q,\vartheta}$ in a more rigorous fashion. The infinite-range interactions cause complications here, because we need to interchange a limit and an infinite sum. For $\tilde{M}\in\M_{q}$ we introduce the notation $\vartheta(\tilde{M})$ to refer to the value of $\vartheta$ in (\ref{fdfhn:definitionrotation}) with $M=\tilde{M}$. 
\begin{lemma}\label{fdfhn:lemmaconvergencetestfunction} Assume that (\asref{fdfhn:aannames}{\text{HS}1}) is satisfied. Fix an integer $q\geq 1$ and consider any sequence $\{M_j\}_{j\in\N}$ in $\M_{q}$ with the property that $\lim_{j\rightarrow\infty}M_j=\infty$ and $\vartheta(M_j)=\vartheta$ for all $j\in\N$ and some $\vartheta\in q^{-1}\Z_q\setminus\{0\}$. Then for any $Z\in  {C}_c^\infty(\R;\ell_{q,\perp;\infty}^2)\subset  {C}_c^\infty(\R;\ell_{q,\perp}^2)$ we have the limit
\begin{equation}
\begin{array}{lcl}
\lim_{j\rightarrow\infty}\nrm{\Delta_{M_j}Z-\Delta_{q,\vartheta}Z}_{ {L}^2(\R,\ell_{q,\perp}^2)}&=&0.
\end{array}
\end{equation}
\end{lemma}
\textit{Proof.} Fix any test function $Z\in  {C}_c^\infty(\R;\ell_{q,\perp;\infty}^2)\subset  {C}_c^\infty(\R;\ell_{q,\perp}^2)$ and pick a sufficiently large $\mu\in\N$ for which $\mathrm{supp}(Z)\subset [-\mu,\mu]$. Without loss of generality we assume that $\nrm{Z}_{ {L}^2(\R,\ell_{q,\perp}^2)}=1$. Pick $\epsilon>0$, together with $K\in\Z_{>\mu}$ in such a way that
\begin{equation}
\begin{array}{lcl}
\tau\sum\limits_{m\geq K-\mu} 4|\alpha_m|&<&\frac{\epsilon}{8}.
\end{array}
\end{equation}
Moreover, by the strong continuity of the shift-semigroup \cite[Example I.5.4]{ENGEL2000}, we can pick $J\in\N$ in such a way that for each $j\geq J$ and each $|m|\leq 4K+4l$ we have the bound
\begin{equation}\label{fdfhn:eq:strongsemi}
\begin{array}{lcl}
\tau|\alpha_m|\nrm{T_{M_j}^mZ-T_{q,\vartheta}^mZ}_{L^2(\R,\ell_{q,\perp}^2)}&=&
\tau|\alpha_m|\nrm{Z(\cdot +m n_jM_j^{-1})-Z(\cdot+m)}_{L^2(\R,\ell_{q,\perp}^2)}\\[0.2cm]
&<&\frac{\epsilon}{16(\mu+K)},
\end{array}
\end{equation}
together with
\begin{equation}
\begin{array}{lcl}
|n_jM_j^{-1}-1|&\leq &\frac{1}{2}.
\end{array}
\end{equation}
Here we introduced $n_j$ for the value of $n$ in (\ref{fdfhn:definitionrotation}) with $M=M_j$. Fix $j\geq J$. Since $\mathrm{supp}(Z)\subset [-\mu,\mu]$, we obtain 
\begin{equation}
\begin{array}{lcl}
\Delta_{M_j}Z(\xi)-\Delta_{q,\vartheta}Z(\xi)&=&\tau\sum\limits_{m\geq K-\mu}\alpha_m\Big[Z(\xi-mn_jM_j^{-1})-Z(\xi-m)\Big]
\end{array}
\end{equation}
for any $\xi>K$, which allows us to estimate
\begin{equation}
\begin{array}{lclcl}
\nrm{\Delta_{M_j}Z-\Delta_{q,\vartheta}Z}_{ {L}^2\big((K,\infty),\ell_{q,\perp}^2\big)}&\leq & \tau\sum\limits_{m\geq K-\mu}2|\alpha_m|\nrm{Z}_{ {L}^2(\R,\ell_{q,\perp}^2)}
&<& \frac{\epsilon}{4}.
\end{array}
\end{equation}
A similar computation yields
\begin{equation}
\begin{array}{lcl}
\nrm{\Delta_{M_j}Z-\Delta_{q,\vartheta}Z}_{ {L}^2\big((-\infty,-K),\ell_{q,\perp}^2\big)}&<& \frac{\epsilon}{4}.
\end{array}
\end{equation}
Finally, for $\xi\in[-K,K]$ we see that
\begin{equation}
\begin{array}{lcl}
\Delta_{M_j}Z(\xi)-\Delta_{q,\vartheta}Z(\xi)&=&\tau\sum\limits_{m=1}^{4l+4K}\alpha_m\Big[Z(\xi+mn_jM_j^{-1})-Z(\xi+m)\\[0.4cm]
&&\qquad \qquad \qquad +Z(\xi-mn_jM_j^{-1})-Z(\xi-m)\Big]\\[0.4cm]
&=&\tau\sum\limits_{m=1}^{4l+4K}\alpha_m\Big[T_{M_j}^mZ(\xi)-T_{q,\vartheta}^mZ(\xi) +T_{M_j}^{-m}Z(\xi)-T_{q,\vartheta}^{-m}Z(\xi)\Big].
\end{array}
\end{equation}
On account of (\ref{fdfhn:eq:strongsemi}), we can, hence, estimate
\begin{equation}
\begin{array}{lcl}
\nrm{\Delta_{M_j}Z-\Delta_{q,\vartheta}Z}_{ {L}^2\big([-K,K],\ell_{q,\perp}^2\big)}&\leq & \tau\sum\limits_{m=1}^{4l+4K}|\alpha_m|\Big[\nrm{T_{M_j}^mZ-T_{q,\vartheta}^mZ}_{L^2(\R,\ell_{q,\perp}^2)}\\[0.4cm]
&&\qquad \qquad +\nrm{T_{M_j}^{-m}Z-T_{q,\vartheta}^{-m}Z}_{L^2(\R,\ell_{q,\perp}^2)}\Big]\\[0.4cm]
&<& 2(4l+4K)\frac{\epsilon}{16(\mu+K)}\\[0.2cm]
&=&\frac{\epsilon}{2}.
\end{array}
\end{equation}
Combining these estimates yields the bound
\begin{equation}\begin{array}{lcl}
\nrm{\Delta_{M_j}Z-\Delta_{q,\vartheta}Z}_{ {L}^2(\R,\ell_{q,\perp}^2)}&<&\epsilon,
\end{array}
\end{equation}
from which the desired limit follows.
\qed\\

\bibliographystyle{plain}

%%%%%%%%%%%%%%%%%%%%%%%%%%%%%%%%%%%%%%%%%%%%%%%%%%%%%%%%%%%%%%%%%%%%%%%%%%%%%%%
\end{document}